\documentclass[english, 11pt]{amsart}

\usepackage{listings}
\usepackage{xcolor}
\usepackage{mathptmx}
\usepackage{tikz}
\usepackage{aliascnt}
\usepackage{mathrsfs}
\usepackage[all,poly,knot]{xy}
\usepackage{hyperref}
\usepackage{comment}  
\usepackage{csquotes}   
\usepackage{amssymb,amsbsy,amsmath,amsfonts,amscd,
	footmisc,fancyhdr,multicol,fancybox,
	graphicx,rotating,ifthen,wasysym}
\usepackage{tikz-cd}
\usepackage{nameref}

\usepackage{pdfpages}
\usepackage{multirow}
\usepackage{nccrules}
\usepackage{textcomp}
\usepackage{framed}
\usetikzlibrary{arrows}
\usepackage{lipsum}
\usepackage{mathtools}

\DeclarePairedDelimiter\floor{\lfloor}{\rfloor}

\usepackage{geometry} % 在 \documentclass 之后、\begin{document} 之前添加这行
\makeatletter
\newcommand{\leqnomode}{\tagsleft@true}
\newcommand{\reqnomode}{\tagsleft@false}
\makeatother

\DeclareMathOperator{\pic}{Pic}

\DeclareMathOperator{\dif}{d}

\DeclareMathOperator{\ord}{ord}

\newcommand{\C}{\mathbb{C}}

\renewcommand{\P}{\mathbb{P}}

\newcommand{\style}[1]{{\sf #1}}
\newcommand{\Jet}{\style{Jet}}

\def\log{\mathrm{log}\,}

\theoremstyle{plain}

\newtheorem{thm}{Theorem}[section]

\newtheorem{lem}[thm]{{Lemma}}

\newtheorem{pro}[thm]{Proposition}

\newtheorem{rem}[thm]{Remark}
\newtheorem{defi}[thm]{Definition}

\newtheorem{KeyObservation}[equation]{Key Observation}

\newtheorem{MiracleFactorizations}[equation]{Miracle Factorizations}

\theoremstyle{remark}
\newtheorem{rmk}[thm]{Remark}

\theoremstyle{definition}

\numberwithin{equation}{section}

\usepackage[hyperpageref]{backref}

\hypersetup{
	colorlinks,
	linkcolor={red!50!black},
	citecolor={blue!62!black},
	urlcolor={blue!80!black}
}
\theoremstyle{plain}
\newcommand{\thistheoremname}{}
\newtheorem*{genericthm*}{\thistheoremname}
\newenvironment{namedthm*}[1]{\renewcommand{\thistheoremname}{#1}%
	\begin{genericthm*}}
	{\end{genericthm*}}

\newtheoremstyle{named}{}{}{\itshape}{}{\bfseries}{.}{.5em}{\thmnote{#3's }#1}
\theoremstyle{named}

\allowdisplaybreaks[4]

\makeatletter
\newcommand\thankssymb[1]{\textsuperscript{\@fnsymbol{#1}}}
\makeatother
\begin{document} 
\title[Second Main Theorem for three conics]{\bf A Second Main Theorem for Entire Curves Intersecting Three Conics}

	\keywords{
	Nevanlinna theory, Second Main Theorem, entire curves, invariant logarithmic $2$-jet differentials, Siu's slanted vector fields,  Demailly-Semple tower}

	\author{Lei Hou}
	\address{Academy of Mathematics and Systems Sciences, Chinese Academy of Sciences, Beijing 100190, China}
	\email{houlei@amss.ac.cn}
	
	\author{Dinh Tuan Huynh}
	
	\address{Department of Mathematics, University of Education, Hue University, Hue City 530000, Vietnam}
	\email{dinhtuanhuynh@hueuni.edu.vn}
	
	\author{
Joël Merker}
\address{ Laboratoire de Math\'ematiques d’Orsay, CNRS, Universit\'e Paris-Saclay, 91405 Orsay
Cedex, France}
\email{joel.merker@universite-paris-saclay.fr}
	
	\author{Song-Yan Xie}
	\address{State Key Laboratory of Mathematical Sciences, Academy of Mathematics and Systems Science, Chinese Academy of Sciences, Beijing 100190, China;  School of Mathematical Sciences, University of Chinese Academy of Sciences, Beijing 100049, China}
	\email{xiesongyan@amss.ac.cn}
	
	\date{\today}
	
\begin{abstract}
We establish a Second Main Theorem for algebraically nondegenerate entire holomorphic curves \( f: \mathbb{C} \to \mathbb{P}^2 \) intersecting a generic configuration of three conics \(\mathcal{C}= \mathcal{C}_1+ \mathcal{C}_2+ \mathcal{C}_3 \) in the complex projective plane $\mathbb{P}^2$. Using invariant logarithmic $2$-jet differentials with negative twists, we prove the estimate
\[
T_f(r) \leqslant 5 \sum_{i=1}^3 N_f^{[1]}(r, \mathcal{C}_i) + o\big(T_f(r)\big)\quad\parallel,
\]
where \( T_f(r) \) is the Nevanlinna characteristic function, and \( N_f^{[1]}(r, \mathcal{C}_i) \) is the $1$-truncated counting function.

We also prove that the logarithmic cotangent bundle of three generic conics
is nef and that all of its negatively twisted symmetric powers have no
nonzero global sections.  Consequently, within the negatively twisted
jet-differential approach, order two is the first level at which an effective
argument can operate.

The key innovation of our approach is establishing new vanishing lemmas of the form
\[
H^0\bigl(\mathbb{P}^2,\, E_{2,m}T_{\mathbb{P}^2}^*(\log \mathcal{C}) \otimes \mathcal{O}_{\mathbb{P}^2}(-t)\bigr) = 0
\]
for specific pairs \((m, t)\), achieved by combining algebro-geometric arguments with computer-assisted computations through a mod-\(p\) reduction technique. This provides a systematic method for proving vanishing results for negatively twisted jet differentials — a key component in complex hyperbolic geometry — thereby filling a missing piece in two-dimensional complex hyperbolicity theory.
\end{abstract}

\thanks{ Appendix~\ref{appendix:results of Pengchao Wang} was authored by Pengchao Wang.}

	\maketitle

\section{\bf Introduction}

\subsection{Motivation}
In the complex projective plane \(\mathbb{P}^2\), the \emph{logarithmic Kobayashi conjecture} asserts that if \(\mathcal{C} \subset \mathbb{P}^2\) is a \emph{generic} algebraic curve of degree \(d \geqslant 5\), then the complement \(\mathbb{P}^2 \setminus \mathcal{C}\) should be hyperbolic \cite{Kobayashi1970}. Various notions of hyperbolicity exist; here we use the most straightforward one: \(\mathbb{P}^2 \setminus \mathcal{C}\) is said to be \emph{hyperbolic} if there is no nonconstant holomorphic map \(f: \mathbb{C} \to \mathbb{P}^2 \setminus \mathcal{C}\). This conjecture generalizes the classical Little Picard Theorem on the hyperbolicity of \(\mathbb{P}^1 \setminus \{\text{3 points}\}\), and is supported by explicit examples constructed by Zaidenberg \cite{Zaidenberg1988}. For generic curves \(\mathcal{C}\) of sufficiently large degree in \(\mathbb{P}^2\), the conjecture was confirmed by Siu and Yeung \cite{siu_yeung1996} by constructing explicit  logarithmic $2$-jet differentials with negative twists. Later, building on the work  \cite{Demailly-Elgoul2000},  El Goul \cite{Elgoul2003} improved the degree bound to \(15\) for \emph{very generic} \(\mathcal{C}\). Combining this with Siu's strategy of using slanted vector fields \cite{Siu2004} and a technique of Păun~\cite{mihaipaun2008}, Rousseau \cite{Rousseau2009} further reduced the degree bound to \(14\).

These results for high-degree curves naturally extend to configurations with multiple components. In the case where the curve \(\mathcal{C} = \cup_{i=1}^q \mathcal{C}_i\) decomposes into the union of \(q \geqslant 2\) generic algebraic curves \(\mathcal{C}_i\) of degree \(d_i\), the hyperbolicity of the complement \(\mathbb{P}^2 \setminus \mathcal{C}\) has been established with smaller degree bounds \(d = \sum_{i=1}^q d_i\) in various works. For instance, a result in \cite{babets1984} 
settles the case where \(q \geqslant 5\). The case \((q = 4, \sum_{i=1}^4 d_i \geqslant 5)\) was treated in \cite{DSW1995-AJM}. The configuration of three conics was studied in \cite{Grauert1989, DSW1995-AJM, Francois-Duval-2001} by different methods.
Dethloff, Schumacher, and Wong treated the case of two conics and one cubic using symmetric differential $1$-forms~\cite{DSW1995-duke}, and settled the case of three generic conics by a clever application of Nevanlinna theory~\cite{DSW1995-AJM}. Rousseau \cite{Rousseau2009} addressed the case \(q = 2\) under the assumptions \((d_1 = 1, d_2 \geqslant 11)\), \((d_1 = 2, d_2 \geqslant 8)\), \((d_1 = 3, d_2 \geqslant 5)\), or \((d_1 = 4, d_2 \geqslant 4)\). Notably, for $q \geqslant 3$, Noguchi, Winkelmann, and Yamanoi \cite{NWA2007} established the \textit{algebraic degeneracy} of every entire curve in $\mathbb{P}^2 \setminus \mathcal{C}$ under the mild condition $\sum_{i=1}^q d_i \geqslant 4$. Their result also extends to the corresponding higher-dimensional setting.

Thus, in the aforementioned settings, every algebraically nondegenerate entire  curve \(f: \mathbb{C} \rightarrow \mathbb{P}^2\) must intersect the curve \(\mathcal{C}\). In this article, we aim to quantify these geometric results by means of Nevanlinna theory. Specifically, we aim to establish a Second Main Theorem type estimate:
\begin{equation}
\label{SMT}
T_f(r) \leqslant c \sum_{i=1}^q N_f^{[k]}(r, \mathcal{C}_i) + o\big(T_f(r)\big) \quad\parallel,
\end{equation}
where
\[
N_f^{[k]}(r, \mathcal{C}_i) := \int_1^r \sum_{|z| < t} \min\{k, \ord_z f^* \mathcal{C}_i\}\, \frac{\dif t}{t} \quad (i = 1, \dots, q; \, r > 1)
\]
are the \textsl{\(k\)-truncated counting functions}, measuring the frequency of intersections \(f(\mathbb{C}) \cap \mathcal{C}_i\); and
\[
T_f(r) := \int_1^r \frac{\dif t}{t} \int_{|z| < t} f^* \omega_{\text{FS}} \quad (r > 1)
\]
is the \textsl{order function}, capturing the area growth rate of  the holomorphic discs \(f(\{|z| < t\})\) with respect to the Fubini-Study metric \(\omega_{\text{FS}}\) on $\mathbb{P}^2$. The constant \(c\) on the right-hand side of~\eqref{SMT} is independent of \(r\), and the notation ``\(\parallel\)'' indicates that the inequality \eqref{SMT} holds for all \(r > 1\) outside a set of finite Lebesgue measure.

When \(\mathcal{C}\) has at least \(4\) irreducible components, such a Second Main Theorem type estimate with truncation at level \(2\) was established by H. Cartan \cite{Cartan1933} in the linear case where \(\mathcal{C}_i\) are lines:
\[
(q-3) T_f(r) \leqslant \sum_{i=1}^q N_f^{[2]}(r, \mathcal{C}_i) + o\big(T_f(r)\big)\quad \parallel.
\]
In the nonlinear case, analogous results (without truncation) were obtained by Eremenko and Sodin \cite{Eremenko-Sodin1992} via pluripotential theory, and by Ru \cite{Minru2009} through a synthesis of Cartan's method and Diophantine geometry. These results also extend to the corresponding higher-dimensional setting (see \cite{Noguchi-Winkelmann2014, Ru21}).

The above estimate becomes vacuous when \(q = 3\), as the factor \((q-3)\) vanishes; thus the case of three curves lies precisely on the borderline where classical methods break down. From the viewpoint of the logarithmic Bloch--Ochiai theorem, the condition \(q > 3\) was originally required to guarantee algebraic degeneracy of every entire curve \(f: \mathbb{C} \to \mathbb{P}^2 \setminus \bigcup_{i=1}^q \mathcal{C}_i\) (see the series of works by Noguchi~\cite{n77,n80,n81}).

%%The problem of establishing \eqref{SMT} remains largely open when \(\mathcal{C}\) has \(q \leqslant 3\) components. In \cite{Huynh-Vu-Xie-2017}, using methods developed for the Green-Griffiths conjecture, we addressed the case \(q = 1\) with a degree bound \(d \geqslant 660\). When the degree of \(\mathcal{C}\) is low, the problem becomes more difficult to track, according to the motto "the lower the degree the harder the hyperbolicity" \cite[p.~3]{Duval2021}.

In the prior work \cite{CHSX2025}, the authors established a Second Main Theorem for $n+1$ hypersurfaces $\{D_i\}_{i=1}^{n+1}$ in $\mathbb{P}^n$ in general position and intersecting transversally, proving the defect relation
\begin{equation}
    \label{CHSX-defect}
\sum_{i=1}^{n+1} \delta_f(D_i) < n+1
\end{equation}
for any algebraically nondegenerate entire curve $f: \mathbb{\mathbb{C}}\rightarrow \mathbb{P}^n$
under the optimal total degree assumption $\sum_{i=1}^{n+1} \deg D_i \geqslant n+2$.
Here,
 each defect $\delta_f(D_i)$ is defined by
\[
\delta_f(D_i) 
\coloneqq 1-\limsup_{r \to \infty}  \frac{N_f(r, D_i)}{\deg(D_i) \, T_f(r)}.
\]
The proof combined parabolic Nevanlinna theory \cite{Siboni-Paun2014} with the Second Main Theorems for semi-Abelian varieties \cite{Noguchi-Yamanoi-Winkelmann2002} and for $n+2$ hypersurfaces~\cite{Ru2004}. 
Subsequently, Ru and Wang~\cite{Ru-Wang2025} extended the result~\eqref{CHSX-defect} to moving targets with truncation level~$1$, while preserving the upper bound $n+1$ on the right-hand side (see also~\cite{Ru-Wang-2026}).

In the case of $n=2$, the defect relation~\eqref{CHSX-defect} might be viewed as the weakest, albeit nontrivial, form of a Second Main Theorem in the vein of~\eqref{SMT} with noneffective $c < \infty$. The primary objective of this article is to establish an effective Second Main Theorem, with truncation level $1$,
for three generic conics with an explicit constant $c = 5$.

\subsection{Two Key Theorems}
The pursuit of a Second Main Theorem encounters a fundamental obstacle: the very nature of entire curves exhibits an inherent and wild flexibility in large-scale geometry, as evidenced by the exotic examples in \cite{HX2021, Xie2024, WX2025, CFX2025}. Our strategy constrains this flexibility by leveraging additional algebro-geometric structures to extract the underlying rigidity of entire curves, building on the following two key results.

\begin{thm}[McQuillan~\cite{Mcquillan1998}]
\label{McQuillan's result}
Let $\mathcal{C}=\sum_{i=1}^q\mathcal{C}_i$ be a normal crossings divisor in $\mathbb{P}^2$, where each $\mathcal{C}_i$ is a smooth algebraic curve of degree $d_i$ with $d=\sum_{i=1}^q d_i\geqslant 4$. Let $f:\mathbb{C}\rightarrow\mathbb{P}^2$ be an algebraically nondegenerate entire holomorphic curve. If the lifting $f_{[1]}:\mathbb{C}\rightarrow X_1=\mathbb{P}\,T_{\mathbb{P}^2}(-\log \mathcal{C})$ is algebraically degenerate, then the following Second Main Theorem type estimate holds:
\[
(d-3)\,T_f(r)\leqslant N_f^{[1]}(r,\mathcal{C}) + o\big(T_f(r)\big)\quad\parallel. 
\]
\end{thm}

\begin{thm}[{\cite[Theorem 3.1]{Huynh-Vu-Xie-2017}}]
\label{smt-form-logarithmic-diff-jet}
Let \(\mathcal{C} = \sum_{i=1}^q \mathcal{C}_i\) be a normal crossings divisor in \(\mathbb{P}^2\). For any logarithmic jet differential
\[
\omega \in H^0\big(\mathbb{P}^2, E_{k,m}^{GG} T_{\mathbb{P}^2}^*(\log \mathcal{C}) \otimes \mathcal{O}_{\mathbb{P}^2}(-t)\big) \quad (k, m, t \geqslant 1),
\]
and any entire holomorphic curve \(f: \mathbb{C} \rightarrow \mathbb{P}^2\) such that \(f^*\omega \not\equiv 0\), the following Second Main Theorem type estimate holds:
\begin{equation}
\label{SMT2}
T_f(r) \leqslant \frac{m}{t} \sum_{i=1}^q N_f^{[1]}(r, \mathcal{C}_i) + O\big(\log T_f(r) + \log r\big)\quad\parallel.
\end{equation}
\end{thm}

Here, \( E_{k,m}^{GG} T_{\mathbb{P}^2}^*(\log \mathcal{C}) \) denotes the logarithmic Green–Griffiths \( k \)-jet bundle of weighted order \( m \) (see Section~\ref{sect: preliminaries}). The use of such higher-order jet differentials in the study of complex hyperbolicity goes back to Andr\'e Bloch~\cite{Bloch1926}. In the present work, we restrict attention to the ``invariant'' subbundle \(E_{k,m} T_{\mathbb{P}^2}^*(\log \mathcal{C}) \subset E_{k,m}^{GG} T_{\mathbb{P}^2}^*(\log \mathcal{C}) \), originally introduced by Demailly~\cite{Demailly1997}. This subbundle possesses stronger positivity properties, making it particularly suitable for constructing sections; further details can be found in Section~\ref{sect: preliminaries}.

Theorem~\ref{smt-form-logarithmic-diff-jet} can be seen as a quantitative refinement of the Fundamental Vanishing Theorem for entire curves (cf.~\cite[pp.~55-59]{Green-Griffiths1980}, \cite[Lemma 1.2]{Siu-Yeung1996} and \cite[Corollary 7.9]{Demailly1997}). Although the original statement in \cite{Huynh-Vu-Xie-2017} is formulated for hypersurfaces in $\mathbb{P}^n$, the proof extends directly to the case of simple normal crossing divisors in projective manifolds (cf.~\cite[Theorem A10.7.2]{Ru21}). 

%Theorem~\ref{smt-form-logarithmic-diff-jet} enables a direct application. Specifically, whenever the existence of sufficiently many nonzero logarithmic jet differentials ensures the complex hyperbolicity of the complement $\mathbb{P}^2 \setminus \mathcal{C}$, it facilitates the derivation of a \emph{quantitative} Second Main Theorem for entire curves in $\mathbb{P}^2$ intersecting $\mathcal{C}$. This passage from qualitative hyperbolicity to precise analytic inequalities represents a key step in our approach.

\subsection{Our First Result}

The first result of this article is an application of Theorem~\ref{smt-form-logarithmic-diff-jet} in the case $k = 1$.

\begin{thm}
\label{smt for three components, case large degrees}
Let $\mathcal{C}=\cup_{i=1}^3\mathcal{C}_i$ be a union of three smooth curves $\mathcal{C}_i$ of degrees $d_i$ with simple normal crossings. Assuming that $d_1\geqslant d_2\geqslant d_3\geqslant 2$ and $d_1\geqslant 3$, so that the total degree $d=d_1+d_2+d_3\geqslant 7$. For a nonconstant algebraically nondegenerate entire holomorphic curve $f\colon \mathbb{C}\rightarrow\mathbb{P}^2$ whose image is not contained in any $\mathcal{C}_i$, one has the estimate
\begin{equation*}
	%%\label{eq:smt-three-components, case large degrees}
	T_f(r) \leqslant c \sum_{i=1}^3 N_f^{[1]}(r, \mathcal{C}_i) + o\big(T_f(r)\big) \quad\parallel,
\end{equation*}
for any positive constant $c$ satisfying 
\[
c > \frac{1}{(d-3)\delta_1}, \quad \text{where } \quad { \delta_1 = \tfrac{1}{2} - \tfrac{\sqrt{9(d-1)^2 - 12(d_1d_2 + d_2d_3 + d_3d_1)}}{6(d-3)}}.
\]
\end{thm}

Here, the numerical conditions on $d_1, d_2, d_3$ are imposed to ensure the positivity of the Euler characteristic of the logarithmic symmetric differential bundles $S^mT_{\mathbb{P}^2}^*(\log \mathscr{C})$. Specifically, we have:
\[
\chi\big(S^mT_{\mathbb{P}^2}^*(\log \mathscr{C})\big)
= L_{d_{1} d_{2} d_{3}} \cdot m^3 
+O(m^2)
> 0, \quad \forall\, m \gg 1,
\]
where $L_{d_{1} d_{2} d_{3}}$ is a positive leading coefficient.

\begin{rem}\rm
In the case where $\mathcal{C}$ consists of two conics and one cubic, a direct computation shows that $\delta_1$ takes the value $\frac{6-\sqrt{33}}{12}$. 
This leads to a Second Main Theorem of the form
\[
T_f(r) \leqslant c \sum_{i=1}^3 N_f^{[1]}(r, \mathcal{C}_i) + o\big(T_f(r)\big)\quad \parallel,
\]
valid for any positive constant $c$ satisfying
$
c > \frac{3}{6-\sqrt{33}} \approx 11.74\ldots$.
This quantifies the hyperbolicity result of $\mathbb{P}^2\setminus \mathcal{C}$ due to Dethloff--Schumacher--Wong \cite{DSW1995-duke}.
\end{rem}

\subsection{Proof  of Theorem~\ref{smt for three components, case large degrees}}
Let \(f\) be algebraically nondegenerate, as in the statement of the theorem.
The key step is to find a nonzero global negatively twisted logarithmic symmetric differential form  
\begin{equation}
    \label{key 1-jet differential existence}
\omega \in H^0\big(\mathbb{P}^2,\, S^m\, T_{\mathbb{P}^2}^*(\log \mathcal{C}) \otimes \mathcal{O}(-t)\big),
\quad m \geqslant 1,\ t \geqslant 1,\ m/t \leqslant c,\ \omega \not\equiv 0.
\end{equation}

\begin{proof}[Proof of Theorem~\ref{smt for three components, case large degrees} assuming~\eqref{key 1-jet differential existence}]
    
If $f^*\omega \not\equiv 0$, then Theorem~\ref{smt-form-logarithmic-diff-jet} yields the desired estimate. On the other hand, if $f^*\omega \equiv 0$, then $f_{[1]}(\mathbb C)$ lies in the proper zero divisor of $\omega$ on $X_1$. We can therefore apply McQuillan's Theorem~\ref{McQuillan's result} to obtain the stronger inequality
\[
(d - 3) \, T_f(r) \leqslant N_f^{[1]}(r, \mathcal{C}) + o(T_f(r))\quad \parallel.
\]
In either case, the conclusion of the theorem holds.
\end{proof}

The proof of~\eqref{key 1-jet differential existence} will be given in Section~\ref{section: The case of three curves with large degrees}, based on  Riemann--Roch computations and a Bogomolov-type vanishing theorem (Theorem~\ref{thm:Bogomolov vanishing theorem}).

\subsection{Our Second, Main  Result: Second Main Theorem for Three Conics.}\label{subsection: main theorem}

Having established an effective Second Main Theorem for configurations of
large degree, we now turn to the more challenging case of three generic
conics $\mathcal C = \mathcal{C}_1 + \mathcal{C}_2 + \mathcal{C}_3$ in $\mathbb P^2$.  
The residue sequence gives $c_1(T_{\mathbb P^2}^*(\log\mathcal C))=3h$ and $c_2(T_{\mathbb P^2}^*(\log\mathcal C))=9h^2$. Here $h$ represents the class of $\mathcal{O}_{\mathbb{P}^2}(1)$, so  satisfies $h^2 = 1$. Consequently,
Hirzebruch--Riemann--Roch yields the exact formula
\begin{equation}\label{eq:one-jet-Euler-characteristic}
 \chi\bigl(S^m T_{\mathbb P^2}^*(\log\mathcal C)\otimes\mathcal{O}_{\mathbb{P}^2}(-t)\bigr)
 =\frac{m+1}{2}\bigl((t-1)(t-2)-3mt\bigr).
\end{equation}
In particular, $\chi(S^m T_{\mathbb P^2}^*(\log\mathcal C))=m+1$; thus the untwisted Euler characteristic
does not have a negative cubic leading term. This naturally raises the corresponding one-jet vanishing question,
which is completely resolved by the following theorem.

\begin{pro}[One-jet vanishing for three conics]
\label{pro:one-jet-vanishing}
For a generic simple-normal-crossing configuration of three smooth conics
$\mathcal C=\mathcal C_1+\mathcal C_2+\mathcal C_3$ in $\mathbb P^2$, the
logarithmic cotangent bundle
\[
 E=T_{\mathbb P^2}^*(\log\mathcal C)
\]
is nef.  Moreover, for every pair of integers $m,t\geqslant1$,
\begin{equation}\label{eq:one-jet-total-vanishing}
 H^0\!\left(
 \mathbb P^2,S^mT_{\mathbb P^2}^*(\log\mathcal C)
 \otimes\mathcal O_{\mathbb P^2}(-t)
 \right)=0.
\end{equation}
\end{pro}

The proof is given in Appendix~\ref{appendix:one-jet-vanishing}.  It is
instructive to compare the theorem with known borderline cases.  For certain
special configurations of six lines, negatively twisted logarithmic
symmetric differentials do exist; see \cite[p.~133]{Hirzebruch} and
\cite[Theorem~1.4]{VZ2000}.  For a generic sextic curve, they do not exist
\cite[Lemma~1.4.1]{Elgoul2003}.  Proposition~\ref{pro:one-jet-vanishing} shows
that three generic conics satisfy the stronger vanishing for every negative
twist.  Thus the one-jet route is genuinely unavailable in this low-degree
setting.  Invariant logarithmic $2$-jets are therefore not merely a technical
convenience: they are the first possible jet order for the negatively twisted
jet-differential strategy.  This necessity gives the two-jet construction in
the present paper additional geometric significance.
This is a precise instance of the principle that ``\textit{the lower the
degree, the harder the hyperbolicity}'' \cite[p.~3]{Duval2021}.  Our main
result is as follows.

\begin{thm}[Second Main Theorem for Three Conics]
    \label{main theorem 3 conics}
For a generic configuration of three conics $\{\mathcal{C}_i\}_{i=1}^3$ in $\mathbb{P}^2$ and a nonconstant algebraically nondegenerate entire holomorphic curve $f \colon \mathbb{C} \to \mathbb{P}^2$ whose image is not contained in any $\mathcal{C}_i$, the following inequality holds:
\begin{equation}
\label{eq:smt-three-conics}
T_f(r) \leqslant 5 \sum_{i=1}^3 N_f^{[1]}(r, \mathcal{C}_i) + o\big(T_f(r)\big) \quad\parallel.
\end{equation}
\end{thm}

%Moreover, the term ``generic'' admits an explicit characterization: the three conics must form a simple normal crossing divisor and satisfy the vanishing condition~\eqref{key vanishings to define generic} specified below. In Section~\ref{section 7}, we develop an effective method which, in principle, can check whether condition~\eqref{key vanishings to define generic} holds for any given explicit equations of three simple normal crossing conics; the underlying algorithm, however, would require extensive computational resources in general. For brevity and to illustrate the method concretely, we therefore carry out the detailed computations only for some tractable Fermat-type case, see Appendix~\ref{app:maple}.

\begin{rmk}
\label{remark 1.7}
The constant $c=5$ in this theorem is smaller than the constant
${\approx}11.74$ appearing in
Theorem~\ref{smt for three components, case large degrees} for the case of
two conics and one cubic.  This may seem counterintuitive, since the smaller
total degree $2+2+2<2+2+3$ would normally suggest a larger constant.
Proposition~\ref{pro:one-jet-vanishing} shows that negatively twisted $1$-jets
do not exist for three generic conics; the improvement therefore comes from
the more refined invariant $2$-jet method.  A detailed $2$-jet analysis is
likewise expected to sharpen the constants for the two-conics-and-one-cubic
case and for the other configurations covered by
Theorem~\ref{smt for three components, case large degrees}.
\end{rmk}

\begin{rmk}
An unexpected computational byproduct appears for the three Fermat-type conics $\mathcal{C}$ defined by the quadratic polynomials
\[
     2 Z_0^2 + Z_1^2 + Z_2^2,\qquad
     Z_0^2 + 2Z_1^2 + Z_2^2, \qquad
     Z_0^2 + Z_1^2 + 2Z_2^2.
\]
For this explicit configuration, the finite calculation detects a
$9$-dimensional space of negatively twisted invariant logarithmic $2$-jets
of weighted degree $5$ and twist $-1$; see
Remark~\ref{m=5 for three Fermat-type conics}.  Together with the preceding
vanishings in weighted degrees $1,2,3,4$, this forces the space to contain
two algebraically independent members and hence yields a Second Main Theorem
with constant $5$ for this Fermat-type configuration.  This byproduct is not
used in the proof of the Main Theorem, whose constant is instead governed by
the finite set of vanishings in Key Vanishing Lemma~II.
\end{rmk}

\begin{rmk}
\label{rmk:algebraically-degenerate-curves}
For a generic configuration of three conics, an effective Second Main
Theorem can also be obtained without the algebraic nondegeneracy
assumption by combining the results of
Song-Yan Xie and Shengyuan Zhao~\cite{XieZhao2026}
on Poincar\'e-type degree bounds for invariant algebraic curves of
foliations arising from jet differentials
with a more refined treatment of algebraically degenerate entire curves.
The constant in the estimate without the algebraic nondegeneracy
assumption, however, may be larger than~$5$.
\end{rmk}

We now outline the main ideas underlying the proof of Theorem~\ref{main theorem 3 conics} for a finite constant $c < \infty$, before refining the argument to achieve the sharper goal $c = 5$. 
Let \(f\) be algebraically nondegenerate, as in Theorem~\ref{main theorem 3 conics}.

The core idea is to construct two independent negatively twisted invariant
logarithmic $2$-jet differentials $\omega_1$ and $\omega_2$:

\begin{itemize}
    \item If $f^*\omega_1\not\equiv0$ or $f^*\omega_2\not\equiv0$, then
    Theorem~\ref{smt-form-logarithmic-diff-jet} directly gives the desired
    Second Main Theorem type estimate.
    
    \item If both pullbacks vanish, then their independence forces the
    lifting $f_{[1]}(\mathbb C)\subset X_1$ to be algebraically degenerate.
    Hence McQuillan's Theorem~\ref{McQuillan's result} applies.
\end{itemize}

In either case, we obtain the desired Second Main Theorem type estimate.

\medskip
In the numerical range where differentiation preserves a negative twist,
Siu's slanted vector fields realize this idea directly: the natural action of
a slanted vector field on the universal logarithmic jet space sends
$\omega_1$ to a genuine global invariant logarithmic $2$-jet differential
$\omega_2$.  When this direct construction is not numerically available, we
retain Miyaoka's zero-locus argument (see also Lu--Yau), in the
jet-differential form used by Demailly--El Goul
~\cite{Miyaoka,Lu-Yau,Demailly-Elgoul2000}, as a logically necessary
fallback. The conceptual core of our method is the use of two independent negatively twisted invariant logarithmic \(2\)-jet differentials. The zero-locus argument is only a fallback for the remaining numerical cases; sufficiently many additional Key Vanishing Lemmas would rule out those cases and thereby make the fallback unnecessary.

The central difficulty is to control the common base locus of the resulting invariant logarithmic \(2\)-jet differentials. Constructions with sufficiently controlled base locus are available in high degree, as in Siu--Yeung~\cite{siu_yeung1996}, but become far more difficult in low degree because of the algebraic complexity of \(2\)-jet differentials and the associated elimination problem; compare \cite{Merker2008}. In the case of three generic conics, this difficulty cannot be bypassed by reverting to the simpler first-order argument: Proposition~\ref{pro:one-jet-vanishing} shows that no negatively twisted logarithmic \(1\)-jet differentials exist.

A notable illustration of these difficulties arises even in the $1$-jet case. Debarre~\cite{Debarre2005} conjectured that for a generic choice of high-degree hypersurfaces $H_1, \dots, H_c \subset \mathbb{P}^N$ with $c \geqslant N/2$, the complete intersection $X := H_1 \cap \cdots \cap H_c$ should have an ample cotangent bundle $T_X^*$, which is equivalent to the existence of sufficiently many negatively twisted $1$-jet differentials on $X$ with empty common base locus. The first proof of Debarre's ampleness conjecture was established by Xie~\cite{Xie2018}, using explicit $1$-jet differentials obtained by Brotbek~\cite{Brotbek2016Symmetricdifferentialforms}  for Fermat-type hypersurfaces. A key step in Xie's approach is a core lemma~\cite[p.~969]{Xie2018} that overcomes substantial difficulties in the elimination process. Shortly thereafter, Brotbek and Darondeau~\cite{BD2018}  found an alternative  proof, which instead applies Nakamaye's theorem~\cite{Nakamaye} to tackle the elimination.

\subsection{Roadmap for Theorem~\ref{main theorem 3 conics}: Towards an Effective Finite Constant $c<\infty$}

We now outline the proof of our main result, focusing first on establishing an explicit finite constant $c$ in the inequality~\eqref{eq:smt-three-conics}. The argument proceeds in several steps.

\medskip

First, for generic  three conics $\mathcal{C} = \mathcal{C}_1 + \mathcal{C}_2 + \mathcal{C}_3$, standard Riemann--Roch computations (detailed in Section~\ref{Section 3}) guarantee the existence of an irreducible, negatively twisted invariant logarithmic $2$-jet differential
\begin{equation}
    \label{omega1, m, t}
\omega_1 \in H^0\big(\mathbb{P}^2, E_{2,m}T^*_{\mathbb{P}^2}(\log \mathcal{C}) \otimes \mathcal{O}_{\mathbb{P}^2}(-t)\big),
\end{equation}
with weighted degree $m \geqslant 1$, vanishing order $t \geqslant 1$. 

If $f^*\omega_1 \not\equiv 0$, then Theorem~\ref{smt-form-logarithmic-diff-jet} immediately yields the desired Second Main Theorem with a finite constant $c$. We may therefore assume $f^*\omega_1 \equiv 0$ in what follows.

\smallskip
The subsequent strategy depends on the values of $(m, t)$:

\smallskip

\smallskip\noindent \textbf{Case 1: $t > 7$}. We apply Siu's strategy~\cite{Siu2004} of slanted vector fields, as realized in dimension two by P\u{a}un~\cite{mihaipaun2008} (compact case) and Rousseau~\cite{Rousseau2009} (logarithmic case).  On the universal logarithmic jet space, a pole-order-$7$ slanted vector field acts naturally on the weighted homogeneous polynomial associated with a holomorphic family extending $\omega_1$ and produces a genuine global invariant logarithmic $2$-jet differential
\[
 \omega_2\in H^0\!\left(\mathbb P^2,
 E_{2,m}T_{\mathbb P^2}^*(\log\mathcal C)\otimes
 \mathcal O_{\mathbb P^2}(7-t)\right).
\]
The field can be chosen so that $\omega_1$ and $\omega_2$ are independent modulo $\Gamma_2$.  Theorem~\ref{smt-form-logarithmic-diff-jet} applied to $\omega_2$ gives the coefficient $m/(t-7)$; if both pullbacks vanish, their independence forces $f_{[1]}$ to be algebraically degenerate.  See Section~\ref{section: Siu's fields}.

\medskip

\noindent\textbf{Case 2:\, $t < m$.}
Following Miyaoka's zero-locus method (see also Lu--Yau), in the form developed by Demailly--El Goul~\cite{Miyaoka,Lu-Yau,Demailly-Elgoul2000} (see Section~\ref{sect: refined Demailly-EG}), Riemann--Roch computations give a nonzero negatively twisted section on the residual zero divisor of $\omega_1$ in the second level $X_2$ of the Demailly--Semple tower; see Proposition~\ref{prop: bigness on Z}. We then distinguish two subcases:
    \begin{itemize}
        \item If the canonical evaluation of the restricted section along $f_{[2]}$ is nonzero, we derive a Second Main Theorem-type estimate from the intrinsic tautological evaluation and the partition-of-unity argument of~\cite[Theorem~3.1]{Huynh-Vu-Xie-2017}; see Theorem~\ref{pro:SMT on Z}.
        \smallskip
        
        \item If the canonical evaluation is identically zero, Lemma~\ref{lem:restricted-zero-evaluation} shows that the lifting $f_{[1]}(\mathbb{C})$ is algebraically degenerate in the first level $X_1$ of the logarithmic Demailly-Semple tower for $(\mathbb{P}^2, \mathcal{C})$, and hence McQuillan's Theorem~\ref{McQuillan's result} applies.
    \end{itemize}

\medskip
\noindent
\textbf{Case 3:  The Remaining Possibilities.}
For the remaining pairs $(m, t)$ not covered by the previous two cases, we establish the following Key Vanishing Lemma I to rule out the existence of any nonzero $\omega_1$ in these settings. Note that vanishing for $(m, t)$ implies vanishing for all $(m, t')$ with $t' > t$. Thus we complete the proof of a Second Main Theorem for three generic conics with a finite constant $c$. Moreover, careful numerical computation shows that we can take $c = 19$; see Subsection~\ref{computing c}.

\subsection{New Vanishing Results and Proof Strategy}

\begin{namedthm*}{Key Vanishing Lemma I}\label{thm:vanishing1}
For a generic configuration of three conics $\mathcal{C} = \cup_{i=1}^3 \mathcal{C}_i$ in $\mathbb{P}^2$, one has:
\[
H^0\big(\mathbb{P}^2, E_{2, m}T_{\mathbb{P}^2}^*(\log \mathcal{C}) \otimes \mathcal{O}_{\mathbb{P}^2}(-t)\big) = 0
\]
for all $(m,t)$ in the set
$
\{ (3,3), (4,4), (5,5), (6,6), (7,7)\}.
$
\end{namedthm*}

While such vanishing results have long been recognized as fundamental in complex hyperbolicity theory, a proof has remained out of reach for lack of a viable method. The decisive innovation is a new exponent bound that converts the a priori infinite extension problem into finitely many exact linear-algebra calculations.

The project began with the question of whether three conics admit an effective Second Main Theorem and with an attempt to use slanted vector fields of pole order $3$. After substantial work, that route did not yield a complete argument; moreover, lowering the pole order would not by itself remove the Key Vanishing obstruction. The present proof therefore uses the classical pole-order-$7$ package and resolves the remaining finite obstruction through the new exponent bound.

We later observed, with the benefit of the pole-lowering ideas developed in
the subsequent papers~\cite{CuiHouLiuXie2026,HouWangXie2026}, that a
pole-order-$4$ version can also be made to close.  We do not pursue that
variant here: the pole-order-$7$ argument keeps the exposition focused on
the genuinely new ingredient of the present paper, namely the exponent
bound and the resulting finite Key Vanishing problem.  Lowering the pole
order changes the finite list, but does not eliminate the need to prove its
members.

\smallskip

Each case in Key Vanishing Lemma~I follows from the corresponding stronger
case in Key Vanishing Lemma~II by twist monotonicity.  The substantive finite
verification is therefore the proof of Key Vanishing Lemma~II, presented in
Subsection~\ref{subsec:all-layer-cpp}.  It has three steps:

\smallskip
\begin{itemize}
    \item \textbf{All-layer finite-dimensional reduction:} The exponent bound is established uniformly in the Wronskian layer \(q\). It restricts each coefficient to a finite-dimensional space of homogeneous polynomials, including the coefficients in all layers that can occur when \(m\) is large.
    
    \smallskip
    
    \item \textbf{Exact C++ verification:} For explicit Fermat-type conics,
    a C++ program constructs the complete coupled transition matrix over
    $\mathbb F_5$ and proves full column rank by exact sparse elimination.
    
    \smallskip
    
    \item \textbf{Geometric extension:} Upper semicontinuity globalizes the
    verified vanishing from this explicit configuration to a nonempty
    Zariski-open set of triples of conics.
\end{itemize}

%The approach of Demailly and El Goul~\cite[Proposition 5.4]{Demailly-Elgoul2000} established the vanishing of global $2$-jet differentials of weighted degrees $3$, $4$, and $5$ for generic smooth surfaces of degree $d \geqslant 6$ in $\mathbb{P}^3$. This result was subsequently extended by El Goul~\cite[Lemma 1.4.4]{Elgoul2003} to the logarithmic case of the complement of a smooth curve of degree $d \geqslant 7$ in $\mathbb{P}^2$. Both proofs rely fundamentally on Pa\u{u}n's ``proportionality lemma'' and Nadel's construction of meromorphic connections~\cite{Nadel}.  The differences in the underlying setting, however, render these approaches insufficient for proving our Key Vanishing Lemma I.

%Another possible approach is to adapt the argument of Pa\u{u}n \cite[pp.~886-889]{mihaipaun2008} and Rousseau \cite[Proposition 17]{Rousseau2009} by employing Siu's slanted vector fields with a twist by $\mathcal{O}_{\mathbb{P}^2}(3)$, aiming to show that a certain associated ordinary differential equation admits no nontrivial solution. If successful, this method would significantly reduce the number of vanishing results required (only $(m, t)=(3, 3)$ would be needed). However, in the case of three conics, we were unable to overcome the substantial algebraic complexity that arises in the analysis of the corresponding ODE; see also Remark~\ref{how to use O(3)}.

\begin{rmk}
An alternative approach to establishing such vanishing results is through cohomological diagram chasing based on the filtration structure~\eqref{filtration of E_2,m}. This method enables systematic construction of long exact sequences in cohomology, thereby facilitating stepwise computation of 
 $H^0\big(\mathbb{P}^2, E_{2,m}T^*_{\mathbb{P}^2}(\log\mathcal{C}) \otimes \mathcal{O}(-t)\big)$. When $m=t$ is not divisible by $3$, the vanishing reduces to negatively twisted logarithmic $1$-jet differentials. This observation is due to Pengchao Wang; see Appendix~\ref{appendix:results of Pengchao Wang}. Proposition~\ref{pro:one-jet-vanishing} supplies all the required symmetric-power vanishings. Consequently, Proposition~\ref{pro:vanishing of E_2,3q+r} gives, in particular, the cases $(m,t)=(4,4),(5,5),(7,7)$ in Key Vanishing Lemma~I.

However, in more general cases—such as when $m > t$ or when $3$ divides $m = t$—this diagram-chasing approach becomes substantially more complicated and its feasibility remains uncertain. Determining the essential information required for the diagram chasing presents a significant challenge. In particular, computing $H^i$ for $i = 0, 1, 2$ for each term in the long exact sequence requires either analyzing the connecting homomorphisms from $H^i$ to $H^{i+1}$, or computing other $H^{i}$ groups that typically do not vanish, thereby introducing considerable additional complexity.
\end{rmk}

\begin{rmk}
The significance of the present method extends well beyond the numerical
constant in the Second Main Theorem.  To the best of our knowledge, this is
the first effective Second Main Theorem in higher-dimensional Nevanlinna
theory whose decisive input is a pair of algebraically independent,
negatively twisted invariant logarithmic $2$-jet differentials, rather than
a Cartan-type Wronskian construction.  Here the distinction concerns the
proof mechanism: the local algebra of invariant $2$-jet differentials still contains its
intrinsic Wronskian generator, but the Second Main Theorem is not produced by
the classical Cartan--Wronskian method.  This introduces a new
algebraic-geometric mechanism for effectiveness.  Equally importantly, the
exponent bound introduced here
reduces the extension problem to finite exact algebra and became a
foundational input in the subsequent
papers~\cite{HHMX2026,CuiHouLiuXie2026}.  In those works, the compact
threshold case $d=15$ and the logarithmic threshold case $d=11$ are reduced
to explicit finite rank-verification problems.  The corresponding computer
calculations are now in progress; no further conceptual ingredient is
required, so reaching the two threshold degrees predicted by the $2$-jet
method is now a matter of computation time.
\end{rmk}

\subsection{Proof Blueprint of Theorem~\ref{main theorem 3 conics} with  $c=5$}\label{subsec:intro more vanishing results}

We first show that in the existence result for $\omega_1$ in~\eqref{omega1, m, t}, the ratio $t/m$ can be made close to $\frac{4-\sqrt{10}}{3}$; see Proposition~\ref{estimate of 2-jet threshold}. If $f^*\omega_1 \not\equiv 0$, then Theorem~\ref{smt-form-logarithmic-diff-jet} yields the estimate~\eqref{eq:smt-three-conics} with constant $c$ around $\frac{3}{4-\sqrt{10}} = 3.58\ldots$. This value represents the theoretical limit of our current approach. Our goal of achieving $c=5$ is strategically set near this limit.

\smallskip
In {\bf Case 1} ($t > 7$), the original differential $\omega_1$ gives the coefficient $m/t$, while the global differential $\omega_2$ produced by a pole-order-$7$ slanted vector field gives the coefficient $m/(t-7)$. Requiring both to be at most $5$ is equivalent to $t \geqslant m/5 + 7$.

\smallskip
In {\bf Case 2} ($t < m$), we refine Demailly–El Goul's approach~\cite{Demailly-Elgoul2000} to estimate the ratio of the vanishing order of $\omega_2$ to its weighted degree. Through the analysis in Section~\ref{sect: refined Demailly-EG}, we show that the restricted-section argument gives any constant $c>3/\tau_1(\frac{t}{m})$, where
{\footnotesize
\[
\tau_1(\tfrac{t}{m}) = \frac{3}{2}\left(4 - \frac{t}{m}\right) - \frac{\sqrt{3}}{2} \sqrt{32 - 8\frac{t}{m} + 3\left(\frac{t}{m}\right)^2}.
\]}
A direct computation shows that $3/\tau_1(\frac{t}{m})<5$ is equivalent to $t/m<43/85$.

\smallskip
Although most positive integer pairs $(m,t)$ satisfy either $t/m<43/85$ or $t \geqslant m/5 + 7$, exceptional cases remain where neither condition holds. To address these, we need the following strengthened vanishing results:

\begin{namedthm*}{Key Vanishing Lemma II}\label{thm:vanishing}
For a generic configuration of three conics $\mathcal{C} = \cup_{i=1}^3 \mathcal{C}_i$ in $\mathbb{P}^2$, one has:
\begin{equation}
    \label{key vanishings to define generic}
H^0\big(\mathbb{P}^2, E_{2, m}T_{\mathbb{P}^2}^*(\log \mathcal{C}) \otimes \mathcal{O}_{\mathbb{P}^2}(-t)\big) = 0
\end{equation}
for all
\[
\begin{aligned}
(m,t)\in\big\{& (3,2),(4,3),(5,3),(6,4),(7,4),(8,5),(9,5),\\
                 & (10,6),(11,6),(12,7),(13,7),(14,8),(15,8),(16,9),\\
                 & (17,9),(18,10),(19,10),(21,11)\big\}.
\end{aligned}
\]
\end{namedthm*}

The proof in Subsection~\ref{subsec:all-layer-cpp} first establishes the
uniform all-layer finite shape.  It then specializes to Fermat-type conics
and proves, by exact full-coupled sparse elimination over $\mathbb F_5$, that
the resulting integral matrix has trivial kernel.

%\medskip
%\subsection{A  Challenge} To avoid the use of any Key Vanishing Lemma and directly obtain, for three generic conics $\mathcal{C}$, a constant \[c \approx \frac{1}{\overline{\theta}_2(\mathbb{P}^2, \mathcal{C})} < 3.58\ldots\quad\text{(see Proposition~\ref{estimate of 2-jet threshold})}\]in the Second Main Theorem, it suffices to verify the following conjecture.

%\begin{conj}\label{key question}\rm
%Let $\mathcal{C} \subset \mathbb{P}^2$ be a configuration of three (generic) conics with simple normal crossings. For any $\epsilon > 0$, there exists a nonzero irreducible negatively twisted invariant logarithmic $2$-jet differential
%\[
%\omega \in H^0\!\left( \mathbb{P}^2,\; E_{2,m}T^*_{\mathbb{P}^2}(\log\mathcal{C}) \otimes \mathcal{O}_{\mathbb{P}^2}(-t) \right)
%\]
%with $\overline{\theta}_2(\mathbb{P}^2, \mathcal{C}) - \epsilon < \frac{t}{m} < \overline{\theta}_2(\mathbb{P}^2, \mathcal{C})$. Here the irreducibility of $\omega$ means that its zero locus $\{\omega = 0\} \subset X_2$ is irreducible and reduced modulo $\Gamma_2$.
%\end{conj}

%We expect that such a property should extend to more general settings, in particular to any “sufficiently varying” family.

\medskip

\subsection{Structure of the article}

Section~\ref{sect: preliminaries} provides necessary background.  
Section~\ref{section: The case of three curves with large degrees} presents detailed Riemann–Roch computations that form the core of the proof of Theorem~\ref{smt for three components, case large degrees}.  
Section~\ref{Section 3} advances to more sophisticated Riemann–Roch computations, specifically designed to establish the existence of negatively twisted invariant logarithmic $2$-jet differentials for three generic conics.  
Section~\ref{section: Siu's fields} details the application of Siu's slanted vector fields in \textbf{Case 1}.  
Section~\ref{sect: refined Demailly-EG} then introduces a refinement of the Demailly–El Goul argument, adapted to handle \textbf{Case 2}.  
The {main innovation} of the article is contained in Section~\ref{section 7}, where we overcome fundamental obstructions in our low-degree setting by establishing Key Vanishing Lemmas.  
Appendix~\ref{appendix:miracle-factorizations} records the explicit quotient
polynomials in the four factorization identities,
Appendix~\ref{appendix:one-jet-vanishing} proves the one-jet vanishing
theorem, and Appendix~\ref{appendix:results of Pengchao Wang} presents
Pengchao Wang's cohomological observation.

\medskip

\section{\bf Preliminaries}
\label{sect: preliminaries}

%In this section, we briefly recall some fundamental concepts and results concerning jet differentials, the Demailly–Semple tower, and a Bogomolov-type vanishing theorem. For a more comprehensive treatment, we refer the reader to \cite{Bloch1926, Green-Griffiths1980, Noguchi1986, Siu-Yeung1996, siu_yeung1997, Demailly1997, Dethloff-Lu2001, Elgoul2003, Siu2004, mihaipaun2008, Rousseau2009}.

\subsection{Logarithmic Jet Bundles and Invariant Differential Forms}

Let $X$ be an $n$-dimensional complex manifold. The $k$-th jet bundle $J_k(X)$ is defined by considering holomorphic germs $f \colon (\mathbb{C}, 0) \to (X, x)$, where two germs are considered equivalent if their Taylor expansions coincide up to order $k$ in some local coordinate system. The equivalence class of a germ $f$, denoted by $j_k(f)$, is said to be regular if $f'(0) \ne 0$.

The total space
$
J_k(X) = \cup_{x \in X} J_k(X)_x
$
carries a natural topology and, together with the natural projection $\pi_k \colon J_k(X) \to X$, forms a holomorphic fiber bundle. Note that $J_1(X)$ coincides with the tangent bundle $T_X$.

Now let $D \subset X$ be a normal crossing divisor. The logarithmic cotangent bundle $T_X^*(\log D)$ is a locally free sheaf generated by meromorphic $1$-forms with at most logarithmic poles along $D$. Precisely, for $x \in X \setminus D$, the stalk $T_{X,x}^*(\log D)$ coincides with the stalk $T_{X,x}^*$ of the usual cotangent sheaf $T_X^*$ at $x$. For $x \in D$, take a neighborhood $U$ of $x$ and irreducible holomorphic functions $\eta_1, \dots , \eta_{\ell}$ on $U$ such that $D \cap U = \{ \eta_1 \cdots \eta_{\ell} = 0 \}$; then the stalk is given by $T_{X,x}^*(\log D) = \sum_{j=1}^{\ell} \mathcal{O}_{X,x} \frac{\mathrm{d} \eta_j}{\eta_j} + T_{X,x}^*$. Its dual, the logarithmic tangent bundle $T_X(-\log D)$, consists of holomorphic vector fields tangent to each component of $D$ (cf.~\cite{Noguchi1986, Elgoul2003}).

Consider an open subset $U \subset X$, a section $\omega \in H^0(U, T_X^*(\log D))$. For a holomorphic curve germ
$f:(\mathbb C,0)\to U\setminus D$, write
\[
f^*\omega=A_f(t)\,dt.
\]
Since $\omega$ is holomorphic on $U\setminus D$, the quantities
$A_f^{(j)}(0)$, for $0\leq j\leq k-1$, depend only on the
ordinary $k$-jet $j_k(f)$. Thus they define a holomorphic map
\begin{equation}\label{eq:log-jet-evaluation}
\widetilde{\omega}:
J_k(U\setminus D)\longrightarrow\mathbb C^k,
\qquad
j_k(f)\longmapsto
\bigl(A_f(0),A_f'(0),\ldots,A_f^{(k-1)}(0)\bigr).
\end{equation}
In particular, for every logarithmic one-form $\omega$ on $U$,
the map in \eqref{eq:log-jet-evaluation} has a unique
holomorphic extension
\begin{equation}\label{trivialization-jet}
    \widetilde{\omega}:
    J_k(X,-\log D)|_U\longrightarrow\mathbb C^k.
\end{equation}
Here the extension is taken on the logarithmic jet bundle,
rather than on the ordinary jet bundle.

A holomorphic section $s:U\to J_kX|_U$ is called a
\emph{logarithmic $k$-jet field} if, for every open subset
$V\subset U$ and every
$\omega\in H^0(V,T_{X}^{*}(\log D))$, the map
\[
\widetilde{\omega}\circ s|_{V\setminus D}
:V\setminus D\longrightarrow\mathbb C^k
\]
extends holomorphically to $V$.
Such an extension, when it exists, is unique.
The resulting sheaf is the sheaf of holomorphic sections
of the logarithmic $k$-jet bundle along $D$, denoted by $J_k(X, -\log D)$ (cf.~\cite{Noguchi1986, Elgoul2003, Rousseau2009}).

Furthermore, given a local basis $\omega_1, \dots, \omega_n$ of $T_X^*(\log D)$ on $U$, we obtain a trivialization:
\begin{equation*}
J_k(X, -\log D)|_{U} \to U \times (\mathbb{C}^k)^n, \quad \sigma \mapsto \big( \pi_k \circ \sigma;\; \tilde\omega_1 \circ \sigma, \dots, \tilde\omega_n \circ \sigma \big).
\end{equation*}

A logarithmic jet differential of order $k$ and weighted degree $m$ at a point $x \in X$ is a polynomial $Q(f^{(1)}, \dots, f^{(k)})$ on the fiber of $J_k(X, -\log D)$ over $x$, satisfying the weighted homogeneity condition:
\begin{equation*}
    Q \big( \lambda f', \lambda^2 f'', \dots, \lambda^k f^{(k)} \big) = \lambda^m Q \big( f', f'', \dots, f^{(k)} \big),
\end{equation*}
for all $\lambda \in \mathbb{C}^*$ and all $(f', f'', \dots, f^{(k)}) \in J_k(X, -\log D)_x$. 

The vector space $E_{k,m}^{GG}T_X^*(\log D)_x$ consists of all such jet differential polynomials defined at the point $x \in X$. Then the total space
\begin{equation*}
    E_{k,m}^{GG}T_X^*(\log D) := \cup_{x \in X} E_{k,m}^{GG}T_X^*(\log D)_x,
\end{equation*}
acquires the structure of a holomorphic vector bundle over $X$ by Faà di Bruno’s formula (cf. e.g.~\cite{Constantine1996, Merker2015}), referred to as the logarithmic Green–Griffiths bundle of jet differentials of order $k$ and weighted degree $m$.

Let $\mathbb{G}_k$ denote the group of germs of $k$-jets of biholomorphisms of $(\mathbb{C}, 0)$:
\begin{equation}\label{reparameterization group G_k}
\mathbb{G}_k \coloneqq \Big{\{} \varphi \colon t \mapsto \sum_{j=1}^k c_j t^j \;\Big|\; c_1 \in \mathbb{C}^*,\; c_j \in \mathbb{C} \text{ for } j \geqslant 2 \Big{\}},
\end{equation}
where the group operation is defined by composition modulo $t^{k+1}$. This group acts naturally on $J_k(X)$ via the right action:$\varphi \cdot j_k(f) \coloneqq j_k(f \circ \varphi)$.

A logarithmic jet differential operator $P \in \Gamma\big(U, E_{k,m}^{GG}T_X^*(\log D)\big)$ over an open set $U \subset X$ is said to be invariant under the reparameterization group $\mathbb{G}_k$ if for every $\varphi \in \mathbb{G}_k$ the following identity holds:
\begin{equation}
\label{define invariant jet differentials}
P \big( (f \circ \varphi)', (f \circ \varphi)'', \dots, (f \circ \varphi)^{(k)} \big) = (\varphi')^m \cdot P \big( f', f'', \dots, f^{(k)} \big).
\end{equation}
The subbundle $E_{k,m} T_X^*(\log D) \subset E_{k,m}^{GG}T_X^*(\log D)$ is then defined as the holomorphic vector bundle whose fibers consist of all such invariant logarithmic jet differential polynomials.

This notion was initially introduced by Demailly \cite{Demailly1997} in the compact case, and subsequently generalized to the logarithmic case by Dethloff-Lu \cite{Dethloff-Lu2001}. %As noted in \cite[p.~4]{Demailly1997}, The subbundle $E_{k,m}$ turns out to have better positivity properties than $E_{k,m}^{GG}$. At least in the case of surfaces of general type, the conditions involved for the existence of sections in $E_{k,m}$ are better than those for $E_{k,m}^{GG}$.

Throughout this article, we only concern the situation  $k=1$ or $2$. In the case that $\dim_{\mathbb{C}} X = 2$, the bundle of invariant logarithmic jet differentials of order $2$ admits a filtration structure  with the grading~\cite{Demailly1997}:
\begin{equation}\label{filtration of E_2,m}
    \mathrm{Gr}^{\bullet} E_{2,m} T_X^*(\log D) = \bigoplus_{0 \leqslant j \leqslant \lfloor m/3 \rfloor} S^{m-3j} T_X^*(\log D) \otimes \overline{K}_X^{j},
\end{equation}
where  $\overline{K}_X $ is  the logarithmic canonical  line bundle $K_X \otimes \mathcal{O}_X(D)$. This filtration corresponds to the expression of an invariant logarithmic $2$-jet differential $P(f', f'')$ outside $D$ as:
\begin{equation*}
    P(f', f'') = \sum_{0 \leqslant j \leqslant \lfloor m/3 \rfloor} \;
    \sum_{\substack{\alpha_1 + \alpha_2 = m - 3j \\ \alpha_1, \alpha_2\geqslant 0}} 
    R_{\alpha_1, \alpha_2, j}(f_1, f_2) \cdot 
    (f_1')^{\alpha_1} (f_2')^{\alpha_2} (f_1' f_2'' - f_2' f_1'')^j,
\end{equation*}
where $f=(f_1, f_2)$ represents the holomorphic germs in local coordinates, and where $R_{\bullet}$ are holomorphic functions. 

\subsection{Logarithmic Demailly-Semple Tower}\label{subsec:log DS tower}
The main reference of this subsection is~\cite{Demailly1997, Dethloff-Lu2001}. 

For simplicity, consider the log-direct manifold $(X_0,D_0,V_0)$ associated with three generic conics $\mathcal{C}_1$, $\mathcal{C}_2$, $\mathcal{C}_3$ in the complex projective plane $\mathbb{P}^2$, where $X_0=X:=\mathbb{P}^2$, $D_0:=\mathcal{C}=\sum_{i=1}^3\mathcal{C}_i$ is a simple normal crossing divisors on $X$ whose irreducible components $\mathcal{C}_i$ are smooth conics and $V_0:=T_X(-\log D_0)$. Define $X_1:=\mathbb{P}(V_0)$ with the natural projection $\pi_{1,0}:X_1\rightarrow X_0$. Set $D_1:=\pi_{1,0}^*D_0$. The subbundle $V_1\subset T_{X_1}(-\log D_1)$ is defined fiberwise by:
\begin{equation}
    \label{define V_1}
V_{1,(x,[v])}:=\{\xi\in T_{X_1,(x,[v])}(-\log D_1)\, |\, \left(\pi_{1,0}\right)_{*}\xi\in \mathbb{C} \cdot v\},
\end{equation}
with $\mathbb{C} \cdot v \subset V_{0,x} = T_{X_{0},x}(-\log D_0)$. Thus we obtain the first level $(X_1, D_1, V_1)$ of the logarithmic Demailly-Semple tower for $(\mathbb{P}^2, \mathcal{C})$.

Proceeding likewise, we can construct the second level of logarithmic Demailly–Semple tower:
\[
(X_2, D_2, V_2)\rightarrow (X_1,D_1,V_1)\rightarrow(X_0,D_0,V_0),
\] 
together with the projections
$\pi_{2,1} \colon X_2\rightarrow X_1$ and
$\pi_{2,0} \colon X_2\rightarrow X_0$. 

Let $\mathcal{O}_{X_{1}}(-1)$ denote the tautological line bundle of $X_{1} = \mathbb{P}(V_0)$ such that $\mathcal{O}_{X_{1}}(-1)_{(x,[v])} = \mathbb{C} \cdot v$. By definition, there is a short exact sequence
\begin{equation}\label{SES:V_1}
    \xymatrix{
    0 \ar[r] & T_{X_{1} / X_{0}} \ar[r] & V_{1} \ar[r]^(0.35){\left(\pi_{1,0}\right)_{*}} & \mathcal{O}_{X_{1}}(-1) \ar[r] & 0 
    },
\end{equation}
where $T_{X_{1} / X_{0}}$ is the relative tangent bundle of the fibration $\pi_{1,0} \colon X_1 \rightarrow X_0$. Moreover, there is a relative version of the Euler sequence of the tangent bundle of the fiber $\mathbb{P}(V_{0,x})$:
\begin{equation}\label{ses:relative Euler sequence}
    \xymatrix{
    0 \ar[r] & \mathcal{O}_{X_1} \ar[r] & \pi_{1,0}^{*} V_{0} \otimes \mathcal{O}_{X_1}(1) \ar[r] & T_{X_{1} / X_{0}} \ar[r] & 0
    }.
\end{equation}

Furthermore, there exists a canonical injection $\mathcal{O}_{X_{2}}(-1) \hookrightarrow \pi_{2,1}^{*} V_{1}$. Combining this with the above sequence yields a natural morphism of line bundles:
\begin{equation*}
    \xymatrix@=1.5cm{
    \mathcal{O}_{X_{2}}(-1) \ar@{^{(}->}[r] & \pi_{2,1}^{*} V_{1} \ar[r]^(0.4){ \pi_{2,1}^{*} \left(\pi_{1,0}\right)_{*}} & \pi_{2,1}^{*} \mathcal{O}_{X_{1}}(-1)
    }, 
\end{equation*}
whose zero divisor is the hyperplane subbundle:
\begin{equation}\label{equ:def of Gamma_2}
    \Gamma_2 \coloneqq \mathbb{P}(T_{X_{1} / X_{0}}) \subset \mathbb{P}(V_1) = X_2.
\end{equation}
Consequently, we can obtain the isomorphism: 
\begin{equation}
\label{gamma_2=O(-1, 1)}
\mathcal{O}_{X_2}(\Gamma_2)
\cong
%%\mathcal{O}_{X_{2}}(-1,1),
\pi_{2,1}^{*} \mathcal{O}_{X_{1}}(-1)\otimes  \mathcal{O}_{X_{2}}(1).
\end{equation}

\medskip

Let $f \colon \mathbb{D}_R \rightarrow X_0$ be a nonconstant holomorphic map from the disc $\mathbb{D}_R=\{z\in\mathbb{C}:|z|<R\}$ of radius $R\in \mathbb{R}_+\cup\{\infty\}$ satisfying $f(\mathbb{D}_R) \not\subset D_0$. To define its liftings to the first level $X_1$ of the logarithmic Demailly–Semple tower, we first consider the natural morphism
\[
\rho:T_{X_0}(-\log D_0)\longrightarrow T_{X_0}
\]
is an isomorphism over $X_0 \setminus D_0$. 
On $\mathbb{D}_R \setminus f^{-1}(D_0)$, define
\[
v_f=\rho^{-1}(f').
\]
With respect to a local parameter $t$ on $\mathbb{D}_R$ and adapted
coordinates satisfying $D=\{z_1\cdots z_q=0\}$, the local expression of $v_f$ is 
\[
v_f=
\left(
\frac{(z_1\circ f)'}{z_1\circ f},
\ldots,
\frac{(z_q\circ f)'}{z_q\circ f},
(z_{q+1}\circ f)',\ldots,(z_2\circ f)'
\right).
\]
Then $v_f$ extends meromorphically across the points of $f^{-1}(D)$.
Locally at any point of $\mathbb{D}_R$, write
\[
v_f=t^\nu b(t),
\qquad
b(t)\ \text{holomorphic},\qquad b(0)\ne0,
\]
by removing the common lowest-order factor.
The regularized first lift is defined by
\[
f_{[1]}(t)=(f(t); [b(t)])\in X_1.
\]
where $[b(t)]$ agrees with $[v_f]$ wherever the latter is directly
defined and is independent of the choices involved.
With this convention, the first lifting $f_{[1]}$ of $f$ is a well-defined holomorphic map from $\mathbb{D}_R$ to $X_1$. Repeating the process for $f_{[1]}$, the second lifting $f_{[2]}$ of $f$ is also well-defined and holomorphic from $\mathbb{D}_R$ to $X_2$.

Moreover, from the definition of $\Gamma_2$ and the identity $\pi_{1,0} \circ f_{[1]} = f$, at every point \(t\) such that \(f(t)\notin D_0\) and \(f'(t)\neq0\), one has \(f_{[2]}(t)\notin\Gamma_2\).

We next explain the corresponding construction of
sections. Let $m\geq0$ and
\[
P\in H^0\bigl(U,E_{2,m}T_X^*(\log D)\bigr).
\]
Write
\[
\alpha_2:
J_2(X,-\log D)^{\mathrm{reg}}
\longrightarrow X_2^{\mathrm{reg}}
\]
for the quotient map and
\[
\widetilde\alpha_2:
J_2(X,-\log D)^{\mathrm{reg}}
\longrightarrow
\mathcal O_{X_2}(-1)|_{X_2^{\mathrm{reg}}}
\]
for its tautological lift, whose values are nonzero.
Over $X\setminus D$, the latter sends $j_2(f)$ to
$f_{[1]}'(0)$; over $D$, it is defined by the
holomorphic logarithmic jet construction.

On $X_2^{\mathrm{reg}}\cap\pi_{2,0}^{-1}(U)$, define
$\sigma_P$ by
\[
\left\langle
\sigma_P(\alpha_2(\xi)),
\widetilde\alpha_2(\xi)^{\otimes m}
\right\rangle
=P(\xi).
\]
The transformation laws
\[
P(\xi\circ\varphi)=\varphi'(0)^mP(\xi),
\qquad
\widetilde\alpha_2(\xi\circ\varphi)
=\varphi'(0)\widetilde\alpha_2(\xi)
\]
make this definition independent of the representative.
The local holomorphic sections of the logarithmic
quotient map show that $\sigma_P$ is holomorphic,
including at regular tower points over $D$.

Over $\pi_{2,0}^{-1}(U\setminus D)$, the ordinary
Demailly--Semple correspondence gives a holomorphic
section of $\mathcal O_{X_2}(m)$, including along
$\Gamma_2$ in this open set. It agrees with the preceding
construction on their common domain. Consequently,
$\sigma_P$ is defined and holomorphic away from
\(
D_2\cap\Gamma_2.
\)
Since $X$ is a surface, $\Gamma_2\to X_1$ is a section
of the projective bundle $X_2\to X_1$, and
$D_2\cap\Gamma_2$ corresponds to the divisor
$D_1\subset X_1$. Hence this intersection has
codimension two in $X_2$.
As $X_2$ is smooth, $\sigma_P$ extends uniquely across
this set as a holomorphic section of
$\mathcal O_{X_2}(m)$.

More generally, the following fundamental isomorphism holds:
\begin{namedthm*}{Direct image formula}[cf.~{\cite[Theorem 6.8]{Demailly1997}}, {\cite[Proposition 3.9]{Dethloff-Lu2001}}]
    The direct image sheaf
    \begin{equation}\label{direct image formula}
        \left( \pi_{2,0} \right)_{*} \mathcal{O}_{X_2} (m) \cong \mathcal{O} \left( E_{2,m}T_{\mathbb{P}^2}^*(\log \mathcal{C}) \right)
    \end{equation}
    is isomorphic to the sheaf of invariant logarithmic $2$-jet differentials of weighted degree $m$.
\end{namedthm*}

\subsection{A Bogomolov-type Vanishing Theorem}

In the subsequent Sections~\ref{section: The case of three curves with large degrees} and \ref{Section 3}, the existence of negatively twisted logarithmic $1$-jet and $2$-jet differentials will be established via  Riemann--Roch calculations. This reduces to verifying the vanishing of the corresponding second cohomology groups:
\[ H^2\big(\mathbb{P}^2, S^{m}T^*_{\mathbb{P}^2}(\log \mathcal{C}) \otimes \mathcal{O}_{\mathbb{P}^2}(-t)\big) \quad \text{and} \quad H^2\big(\mathbb{P}^2, E_{2,m}T^*_{\mathbb{P}^2}(\log \mathcal{C}) \otimes \mathcal{O}_{\mathbb{P}^2}(-t)\big),\]
respectively. A key tool in proving these vanishing results is the following Bogomolov-type vanishing theorem.

\begin{thm}[cf.~{\cite[Theorem 4.4]{Chen2025}}]\label{thm:Bogomolov vanishing theorem}
    Let $(X,D)$ be a log pair where $X$ is a smooth projective surface and $D$ is a reduced simple normal crossing divisor. 
    Assume that $(X,D)$ is of general type (i.e., $K_{X} + D$ is big) and $K_{X} + D$ is nef. Then
    \begin{equation*}
        H^0 \big( X , S^p T_{X}(- \log D) \otimes \mathcal{O}_X(q(K_{X} + D)) \big) = 0,
    \end{equation*}
    for all integers $p, q$ satisfying $p > 2q$.
\end{thm}

In the compact case, Bogomolov \cite{Bogomolov1978} established a seminal vanishing theorem under the assumption that tangent bundle $T_{X}$ is semi-stable. Later, Enoki \cite{Enoki1987} showed that $T_{X}$ is $K_{X}$-semi-stable for minimal K\"ahler space $X$ with big $K_{X}$, and Tsuji \cite{Tsuji1988} proved that $T_{X}$ is $K_{X}$-semi-stable for smooth minimal algebraic variety (i.e., $K_{X}$ is nef). Combining these results, the Bogomolov-type vanishing theorem holds for projective variety $X$ with big and nef $K_{X}$; see~\cite[Theorem 14.1]{Demailly1997}.

%For the logarithmic case, the proof of the theorem \ref{thm:Bogomolov vanishing theorem} also ties to the semi-stability of $T_{X}(- \log D)$. Under the assumptions that $K_{X} + D$ is ample and $T_{X}(- \log D)$ is semi-stable, a brief proof is provided in the Appendix. In our three conics case, $K_{\mathbb{P}^2}(\mathcal{C})$ is ample, so $T_{\mathbb{P}^2}(-\log\mathcal{C})$ is semi-stable with respect to $K_{\mathbb{P}^2}(\mathcal{C})$ by Tsuji's result \cite[Theorem 5.1]{Tsuji1988} (see also \cite{Guenancia2016}).

\medskip

\section{\bf Proof of Theorem~\ref{smt for three components, case large degrees}}
\label{section: The case of three curves with large degrees}
Let $\mathcal{C} = \sum_{i=1}^3 \mathcal{C}_i$ be a simple normal crossing divisor of total degree $d = d_1 + d_2 + d_3$, where $d_i = \deg \mathcal{C}_i$ satisfy  
$
d_1 \geqslant d_2 \geqslant d_3 \geqslant 2$ and $d_1 \geqslant 3.
$  
To derive an effective Second Main Theorem, we seek  nontrivial global sections of the twisted sheaf  
\[
S^m T_{\mathbb{P}^2}^*(\log \mathcal{C}) \otimes \mathcal{O}_{\mathbb{P}^2}(1)^{-(d-3)\delta m} 
\cong 
S^m T_{\mathbb{P}^2}^*(\log \mathcal{C}) \otimes \overline{K}_{\mathbb{P}^2}^{\, -\delta m},
\]  
where $\delta > 0$ is a small rational number to be chosen, and $m$ is taken sufficiently divisible.

\begin{defi}[cf.~{\cite[p.~517]{Demailly-Elgoul2000}}, {\cite[Definition 1.3.1]{Elgoul2003}}]
    \label{def:one-jet-log-threshold}
    \normalfont
    The \emph{$1$-jet log threshold} of $(\mathbb{P}^2,\mathcal{C})$ is defined by
    \[
    \overline{\theta}_1(\mathbb{P}^2,\mathcal{C})
    \coloneqq
    \sup\left\{
    \frac{t}{m} \,\middle|\, m, t > 0,\;
    H^0\big(\mathbb{P}^2, S^m T_{\mathbb{P}^2}^*(\log \mathcal{C}) \otimes \mathcal{O}_{\mathbb{P}^2}(-t)\big) \neq 0
    \right\}.
    \]
    If no such pair $(m, t)$ exists, the threshold is defined to be zero.
\end{defi}

For a generic configuration of three conics,
Proposition~\ref{pro:one-jet-vanishing} gives
$\overline{\theta}_1(\mathbb P^2,\mathcal C)=0$.

\begin{pro}[Proof of~\eqref{key 1-jet differential existence}]
\label{existence of symmetric differentials 1 form, three components case}
Let $\mathcal{C}=\cup_{i=1}^3\mathcal{C}_i$ be a union of three smooth curves $\mathcal{C}_i$ of degrees $d_i$ with simple normal crossings. Assume that $d_1\geqslant d_2\geqslant d_3\geqslant 2$ and $d_1\geqslant 3$. Then one has the estimate
\[
0 < \frac{1}{\overline{\theta}_1(\mathbb{P}^2,\mathcal{C})} \leqslant \frac{1}{(d-3)\delta_1},
\]
where{\footnotesize
\begin{equation}
	\label{definition of delta 1, three components case with large degree}
	\delta_1 = \frac{1}{2}-\dfrac{\sqrt{9(d-1)^2-12(d_1d_2+d_2d_3+d_3d_1)}}{6(d-3)}>0.
\end{equation}}
\end{pro}

\begin{proof}
We follow the standard Riemann-Roch argument (cf.~\cite{Demailly-Elgoul2000, Elgoul2003}).

	By Serre duality, the dual of
	\begin{equation}\label{use Serre duality_1 jet}
	H^2\big(\mathbb{P}^2,S^mT_{\mathbb{P}^2}^*(\log \mathcal{C})
	\otimes \overline{K}_{\mathbb{P}^2}^{-\delta m}\big)
	\end{equation}
	is
	\[
	H^0\big(\mathbb{P}^2,K_{\mathbb{P}^2}\otimes
	S^mT_{\mathbb{P}^2}(-\log\mathcal{C})\otimes
	\overline{K}_{\mathbb{P}^2}^{\delta m}\big).
	\]
	Take $m$ sufficiently divisible that $\delta m$ is an integer.  Multiplication
	by a nonzero section of $\mathcal O_{\mathbb P^2}(3)$ injects the last sheaf
	into
	\[
	S^mT_{\mathbb{P}^2}(-\log\mathcal C)\otimes
	\overline K_{\mathbb P^2}^{\delta m}.
	\]
	For $0<\delta<1/2$ one has $m>2\delta m$.  The
	Bogomolov-type Vanishing Theorem~\ref{thm:Bogomolov vanishing theorem}
	therefore annihilates the target, and hence also the Serre-dual group.  We
	obtain
	\[
	H^2\big(\mathbb{P}^2,S^mT_{\mathbb{P}^2}^*(\log \mathcal{C}) \otimes \overline{K}_{\mathbb{P}^2}^{-\delta m}\big)=0
	\]
	which yields
	\[
	\dim H^0\Big(\mathbb{P}^2, S^mT_{\mathbb{P}^2}^*(\log \mathcal{C}) \otimes \overline{K}_{\mathbb{P}^2}^{-\delta m}\Big) 
	\geqslant \chi\Big(S^mT_{\mathbb{P}^2}^*(\log \mathcal{C}) \otimes \overline{K}_{\mathbb{P}^2}^{-\delta m}\Big),
	\]
	for large $m$. A standard Riemann-Roch computation gives
	\begin{align*}
		\chi\Big(S^mT_{\mathbb{P}^2}^*(\log \mathcal{C}) \otimes \overline{K}_{\mathbb{P}^2}^{-\delta m}\Big)
		&=
		\dfrac{1}{6}\big((m-\delta m)^3-(-\delta m)^3\big)\bar{c}_1^2
		-
		\dfrac{1}{6}m^3\bar{c}_2+O(m^2)\\
		&=\dfrac{1}{2}m^3\bigg((d-3)^2\delta^2-(d-3)^2\delta+\frac{1}{3}(d_1d_2+d_2d_3+d_3d_1)-d+2\bigg)+O(m^2).
	\end{align*}
	The leading coefficient $\varphi(\delta)=(d-3)^2\delta^2-(d-3)^2\delta+\frac{1}{3}(d_1d_2+d_2d_3+d_3d_1)-d+2$ has two positive roots given by
	{\footnotesize
	\[
	\delta_1 = \frac{1}{2}-\dfrac{\sqrt{9(d-1)^2-12(d_1d_2+d_2d_3+d_3d_1)}}{6(d-3)}, \quad  \quad \delta_2 =  \frac{1}{2}+\dfrac{\sqrt{9(d-1)^2-12(d_1d_2+d_2d_3+d_3d_1)}}{6(d-3)}>0.
	\]} We can check the
	 positivity of $\delta_1$  directly. Hence $\varphi(\delta) > 0$ for all $\delta \in (0,\delta_1)$. Thus we conclude the proof.
\end{proof}

\medskip

\section{\bf Twisted Invariant Logarithmic $2$-Jet Differentials and Threshold}
 \label{Section 3}

To obtain a constant $c$ as small as possible in the Second Main Theorem for three conics $\mathcal{C} \subset  \mathbb{P}^2$, 
by Theorem~\eqref{smt-form-logarithmic-diff-jet}, we seek  negatively twisted invariant logarithmic $2$-jet differentials 
with relatively high vanishing order $t$ compared to the degree $m$. This motivates the study of the following threshold:

\begin{defi}[cf.~{\cite[p.~517]{Demailly-Elgoul2000}}, {\cite[Definition 1.3.1]{Elgoul2003}}]
    \normalfont
    The \emph{$2$-jet log threshold} of $(\mathbb{P}^2,\mathcal{C})$ is defined by
    \[
    \overline{\theta}_2(\mathbb{P}^2, \mathcal{C})
    \coloneqq
    \sup\left\{
    \frac{t}{m} \,\middle|\, m, t > 0,\;
    H^0\big(X_2, \mathcal{O}_{X_2}(m) \otimes \pi_{2,0}^* \mathcal{O}_X(-t)\big) \neq 0
    \right\}.
    \]
    If no such pair $(m, t)$ exists, the threshold is defined to be zero.
\end{defi}

If $\overline{\theta}_2 > 0$, then for any sufficiently small $\epsilon > 0$, there exists a negatively twisted logarithmic $2$-jet differential 
\[
\omega \in H^0\big(\mathbb{P}^2, E_{2,m}T_{\mathbb{P}^2}^*(\log \mathcal{C}) \otimes \mathcal{O}_{\mathbb{P}^2}(-t)\big)
\]
with ratio $\overline{\theta}_2 - \epsilon < \frac{t}{m} \leqslant \overline{\theta}_2$.
Furthermore, if $f^* \omega \not\equiv 0$ also holds, then Theorem~\ref{smt-form-logarithmic-diff-jet} yields the Second Main Theorem with constant
\[
\frac{1}{\overline{\theta}_2} \leqslant c = \frac{m}{t} < \frac{1}{\overline{\theta}_2 - \epsilon}.
\]

The significance of the $2$-jet log  threshold stems from its role as the fundamental barrier that dictates the possibility of constructing jet differentials with a prescribed ratio of vanishing orders to weighted degrees.

The following result establishes an explicit upper bound for $\frac{1}{\overline{\theta}_2}$ in the case of three smooth conics with simple normal crossings:
 
\begin{pro}\label{estimate of 2-jet threshold}
Let $\mathcal{C}=\cup_{i=1}^3\mathcal{C}_i$ be a union of three smooth conics $\mathcal{C}_i$ with simple normal crossings. Then one has the estimate
\[
0 < \frac{1}{\overline{\theta}_2(\mathbb{P}^2, \mathcal{C})} \leqslant \frac{3}{4-\sqrt{10}}=3.58\ldots.
\]
\end{pro}

The proof employs a Riemann–Roch calculation, which is standard  (cf. \cite{Demailly-Elgoul2000}).

\begin{proof}
We seek nonzero invariant logarithmic $2$-jet differential forms with high vanishing order, which correspond to global sections of the twisted vector bundle
\[
E_{2,m}T_{\mathbb{P}^2}^*(\log \mathcal{C})\otimes\mathcal{O}_{\mathbb{P}^2}(1)^{-3\delta m} 
\cong 
E_{2,m}T_{\mathbb{P}^2}^*(\log \mathcal{C}) \otimes \overline{K}_{\mathbb{P}^2}^{-\delta m},
\]
where the logarithmic canonical line bundle $$\overline{K}_{\mathbb{P}^2}=K_{\mathbb{P}^2}\otimes\mathcal{O}(\mathcal{C})\cong \mathcal{O}_{\mathbb{P}^2}(-3)\otimes \mathcal{O}_{\mathbb{P}^2}(6)\cong\mathcal{O}_{\mathbb{P}^2}(3),$$ with $\delta>0$ a rational number to be chosen and $m$ sufficiently divisible.

The standard approach begins with Serre duality:
\begin{align}
\label{use Serre duality}
H^2\big(\mathbb{P}^2,E_{2,m}T_{\mathbb{P}^2}^*(\log \mathcal{C}) \otimes \overline{K}_{\mathbb{P}^2}^{-\delta m}\big)^\vee = H^0\big(\mathbb{P}^2,K_{\mathbb{P}^2}\otimes \big(E_{2,m}T_{\mathbb{P}^2}^*(\log\mathcal{C})\big)^\vee\otimes \overline{K}_{\mathbb{P}^2}^{\delta m}\big).
\end{align}
The bundle on the right admits the filtration induced by the dual of
\eqref{filtration of E_2,m}, with graded pieces
\[
S^{m-3j}T_{\mathbb{P}^2}(-\log\mathcal{C})\otimes \overline{K}_{\mathbb{P}^2}^{-j}\otimes K_{\mathbb{P}^2}\otimes\overline{K}_{\mathbb{P}^2}^{\delta m} 
= S^{m-3j}T_{\mathbb{P}^2}(-\log\mathcal{C})\otimes \overline{K}_{\mathbb{P}^2}^{\delta m-j-1}
\]
for $0\leqslant j\leqslant \frac{m}{3}$. 

For any rational number $\delta\in (0,\frac{1}{3})$ and sufficiently divisible $m\gg 1$, the inequality
\[
m - 3j > 2(\delta m - j - 1)
\]
automatically holds for all $0 \leqslant j \leqslant \frac{m}{3}$. By applying the Bogomolov-type Vanishing Theorem~\ref{thm:Bogomolov vanishing theorem}, we deduce:
\[
H^0\Big(\mathbb{P}^2, S^{m-3j}T_{\mathbb{P}^2}(-\log\mathcal{C}) \otimes \overline{K}_{\mathbb{P}^2}^{\delta m - j - 1}\Big) = 0 \quad \,\,\,\text{for all }\quad 0 \leqslant j \leqslant \tfrac{m}{3}.
\]

This establishes the vanishing of the right-hand side in~\eqref{use Serre duality}, which consequently implies the vanishing of the left-hand side as well. We thus obtain the lower bound
\[
\dim H^0\Big(\mathbb{P}^2, E_{2,m}T_{\mathbb{P}^2}^*(\log \mathcal{C}) \otimes \overline{K}_{\mathbb{P}^2}^{-\delta m}\Big) 
\geqslant \chi\Big(E_{2,m}T_{\mathbb{P}^2}^*(\log \mathcal{C}) \otimes \overline{K}_{\mathbb{P}^2}^{-\delta m}\Big),
\]
since the Euler characteristic decomposes as
\[
\chi(\cdots) = \dim H^0(\cdots) - \dim H^1(\cdots) + \dim H^2(\cdots) = \dim H^0(\cdots) - \dim H^1(\cdots).
\]

To estimate the asymptotic growth of the Euler characteristic, we employ the Schur bundle filtration (see~\eqref{filtration of E_2,m}):
\begin{align}\label{equ:filtration of 2-jet}
\mathrm{Gr}^{\bullet}
E_{2,m}T_{\mathbb{P}^2}^*(\log \mathcal{C}) 
&= \bigoplus_{a+3b=m}\Gamma^{(a+b,\,b)}T_{\mathbb{P}^2}^*(\log \mathcal{C}) 
\qquad \text{(cf.~\cite[p.~67]{Demailly1997}, \cite[p.~474]{Elgoul2003})} \notag \\
&= \bigoplus_{0\leqslant b\leqslant\frac{m}{3}}\Gamma^{(m-2b,\, b)}T_{\mathbb{P}^2}^*(\log \mathcal{C}).
\end{align}
Together with the Pieri formula~(cf.~\cite[p.~455]{FultonHarris1991}, \cite[p.~66]{Demailly1997}):
\[
\Gamma^{(\ell_1,\, \ell_2)}T_{\mathbb{P}^2}^*(\log \mathcal{C})\otimes \overline{K}_{\mathbb{P}^2}^{-\delta m} 
= \Gamma^{(\ell_1-\delta m,\, \ell_2-\delta m)}T_{\mathbb{P}^2}^*(\log \mathcal{C}),
\]
we derive the following computation (cf.~\cite[p.~77]{Merker2015}):
\begin{align*}
&\chi\big(E_{2,m}T_{\mathbb{P}^2}^*(\log \mathcal{C}) \otimes \overline{K}_{\mathbb{P}^2}^{-\delta m}\big) \\
=\,& \sum_{b=0}^{\frac{m}{3}} \chi\big(\Gamma^{(m-2b-\delta m,\, b-\delta m)}T_{\mathbb{P}^2}^*(\log \mathcal{C})\big) \\
=\, &\bigg(\frac{1}{6} \int_{0}^{\frac{m}{3}} \big[(m-2b-\delta m)^3-(b-\delta m)^3\big] \mathrm{d} b\bigg) \bar{c}_1^2 \,-\, \bigg(\frac{1}{6} \int_0^{\frac{m}{3}} (m-3b)^3 \mathrm{d} b\bigg) \bar{c}_2 + O(m^3) \\
=\,& \frac{m^2(m+3)}{648}(54\delta^2 m-48\delta m+13m-18\delta+15)\bar{c}_1^2 
 - \frac{m^4}{72}\bar{c}_2 + O(m^3) \\
=\,& \frac{m^4}{648}\big((54\delta^2-48\delta+13)\bar{c}_1^2 - 9\bar{c}_2\big) + O(m^3).
\end{align*}
Here the logarithmic Chern numbers are determined by the formula $\bar{c}_1^2=(d-3)^2=9$ and $\bar{c}_2=d^2 - 3d + 3 - \sum_{i< j}d_id_j=9$, where $d_1=d_2=d_3=2$ and $d=d_1+d_2+d_3=6$.

The quadratic leading coefficient (after removing the common factor  $\frac{9}{648}$)
\[
\varphi(\delta) := (54\delta^2 - 48\delta + 13) - 9 = 54\delta^2 - 48\delta + 4
\]
has two positive roots given by
$
\delta_1 = \frac{4-\sqrt{10}}{9}$ and $ \delta_2 = \frac{4+\sqrt{10}}{9}$.

Hence, $\varphi(\delta) > 0$ for all $\delta \in (0,\delta_1)$. It follows that $\overline{\theta}_2 \geqslant 3\delta_1=0.279\ldots$, which  completes the proof of Proposition~\ref{estimate of 2-jet threshold}.
\end{proof}
 
Siu's strategy of slanted vector fields requires families of negatively twisted jet differentials; see~\cite[pp.~558--559]{Siu2004}. In Proposition~\ref{prop:step1}, we will use algebraic cohomology and base change on the relative logarithmic tower to extend the section defining the residual divisor.

Now we apply Demailly--El Goul's argument (cf.~\cite[p.~532]{Demailly-Elgoul2000}, \cite[Lemma 1.3.2]{Elgoul2003}) to show the existence of ``irreducible'' negatively twisted invariant logarithmic $2$-jet differentials.
 
\begin{pro}\label{prop:step1}
Let $S=(\mathbb P^5)^3$ be the parameter space of ordered triples of
plane conics, and let $S^{\mathrm{snc}}\subset S$ be the nonempty
Zariski-open subset parametrizing triples of smooth conics whose union
$\mathcal C_a=\mathcal C_{1,a}+\mathcal C_{2,a}+\mathcal C_{3,a}$
is a simple normal crossing divisor.
There exist positive integers $m,t$, integers $b_1,b_2$ with
\[
 b_1\geqslant2b_2>0,\qquad b_1+b_2=m,\qquad \frac{t}{m}>\frac15,
\]
and a nonempty Zariski-open subset $U^{m,t}\subset S^{\mathrm{snc}}$
on which there is an algebraic, hence holomorphic, family of nonzero
invariant logarithmic $2$-jet differentials
\[
 \omega_{1,a}\in H^0\!\left(
 \mathbb P^2,E_{2,m}T_{\mathbb P^2}^*(\log\mathcal C_a)
 \otimes\mathcal O_{\mathbb P^2}(-t)\right),\qquad a\in U^{m,t}.
\]
Viewed as a section on the second logarithmic Demailly--Semple stage
$X_{2,a}$ via~\eqref{direct image formula}, each $\omega_{1,a}$ has
zero divisor
\[
 \operatorname{div}(\omega_{1,a})=Z_a+b_1\Gamma_{2,a},
\]
where $Z_a$ is irreducible and reduced and $Z_a\neq\Gamma_{2,a}$.
In particular,
\[
 Z_a\sim b_1u_{1,a}+b_2u_{2,a}-t\pi_{2,0}^*h,
\]
where $u_{1,a}$ and $u_{2,a}$ denote the classes of
$\pi_{2,1}^*\mathcal O_{X_{1,a}}(1)$ and $\mathcal O_{X_{2,a}}(1)$,
respectively, and $h$ denotes the class of $\mathcal O_{\mathbb P^2}(1)$.
\end{pro}

\begin{proof}
We first construct the logarithmic Demailly--Semple tower relatively
over $S^{\mathrm{snc}}$.  Let
$\mathscr C\subset\mathbb P^2\times S^{\mathrm{snc}}$ be the universal
relative simple normal crossing divisor (see \eqref{family of curves in P^2}).  The relative logarithmic
tangent bundle
\[
 T_{\mathbb P^2\times S^{\mathrm{snc}}/S^{\mathrm{snc}}}(-\log\mathscr C)
\]
is locally free of rank two.  Its relative tower consists of projective
bundles
\[
 \mathscr X_2\xrightarrow{\varPi_{2,1}}\mathscr X_1
 \xrightarrow{\varPi_{1,0}}\mathbb P^2\times S^{\mathrm{snc}},
\]
with a smooth projective structure morphism
$q\colon\mathscr X_2\to S^{\mathrm{snc}}$ and a relative vertical divisor
$\varGamma_2$.  Their fibres are $X_{2,a}$, $X_{1,a}$ and
$\Gamma_{2,a}$, respectively.  We write
$\varPi_{2,0}=\varPi_{1,0}\circ\varPi_{2,1}$ and suppress the projection to
$\mathbb P^2$ when pulling back line bundles from the base surface.
For $(c_1,c_2,s)\in\mathbb Z^3$, set
\[
 \mathscr L_{c_1,c_2,s}:=
 \varPi_{2,1}^*\mathcal O_{\mathscr X_1}(c_1)
 \otimes\mathcal O_{\mathscr X_2}(c_2)
 \otimes\varPi_{2,0}^*\mathcal O_{\mathbb P^2}(-s).
\]
Generic cohomology and base change for this projective flat family
(cf.~\cite[Chapter~III, Section~12]{HartshorneGTM52}) gives a dense
Zariski-open subset $S_{c_1,c_2,s}\subset S^{\mathrm{snc}}$ on which
$q_*\mathscr L_{c_1,c_2,s}$ is locally free and the natural maps
\[
 (q_*\mathscr L_{c_1,c_2,s})\otimes\mathbb C(a)
 \longrightarrow
 H^0(X_{2,a},\mathscr L_{c_1,c_2,s}|_{X_{2,a}})
\]
are isomorphisms.
Since there are only countably many integer triples and
$S^{\mathrm{snc}}$ is irreducible over the uncountable field $\mathbb C$,
we may choose a complex point
\[
 a_0\in\bigcap_{(c_1,c_2,s)\in\mathbb Z^3}S_{c_1,c_2,s}
\]
which also lies in the generic locus of
Proposition~\ref{pro:one-jet-vanishing}.

Proposition~\ref{estimate of 2-jet threshold} gives
\[
 \overline\theta_2(\mathbb P^2,\mathcal C_{a_0})
 \geqslant\frac{4-\sqrt{10}}3>\frac15.
\]
By the definition of the threshold and~\eqref{direct image formula},
there exist positive integers $M,T$ with $T/M>1/5$ and a nonzero section
\[
 \omega_0\in H^0\!\left(X_{2,a_0},
 \mathcal O_{X_{2,a_0}}(M)
 \otimes\pi_{2,0}^*\mathcal O_{\mathbb P^2}(-T)\right).
\]
Write $\operatorname{div}(\omega_0)=\sum_jp_jZ_j$, where the $Z_j$ are
distinct irreducible reduced divisors and $p_j\geqslant1$.
With $u_i=u_{i,a_0}$, the projective-bundle formula for the Picard group gives
\[
 \operatorname{Pic}(X_{2,a_0})
 =\mathbb Z\pi_{2,0}^*h\oplus\mathbb Zu_1\oplus\mathbb Zu_2.
\]
Thus we can write
\[
 Z_j\sim b_{1,j}u_1+b_{2,j}u_2-t_j\pi_{2,0}^*h,
 \qquad b_{1,j},b_{2,j},t_j\in\mathbb Z.
\]
The divisor-class argument of
\cite[Lemma~3.3(d)]{Demailly-Elgoul2000}, applied to the logarithmic tower
as in~\cite[Lemma~1.3.2]{Elgoul2003}, gives the following three cases:
\begin{itemize}
 \item $(b_{1,j},b_{2,j})=(0,0)$ and $t_j<0$;
 in this case $Z_j$ is pulled back from a curve in $\mathbb P^2$;
 \item $Z_j=\Gamma_{2,a_0}$, with
 $(b_{1,j},b_{2,j},t_j)=(-1,1,0)$;
 \item $b_{1,j}\geqslant2b_{2,j}\geqslant0$ and
 $m_j:=b_{1,j}+b_{2,j}>0$.
\end{itemize}
In the third case, the relation
\[
\pi_{2,1}^*\mathcal O_{X_{1,a_0}}(b_1)
\otimes
\mathcal O_{X_{2,a_0}}(b_2)
\cong
\mathcal O_{X_{2,a_0}}(b_1+b_2)
\otimes
\mathcal O_{X_{2,a_0}}(-b_1\Gamma_{2,a_0}), \quad \text{see \eqref{gamma_2=O(-1, 1)}}
\]
gives
\[
Z_j
\sim
b_{1,j}u_1+b_{2,j}u_2-t_j\pi_{2,0}^*h
\sim
m_j u_2-b_{1,j}\Gamma_{2,a_0}
-t_j\pi_{2,0}^*h.
\]

Let \(I\) denote the set of indices corresponding to the third case.
Since
\(
\sum_j p_j[Z_j]
=
M u_2-T\pi_{2,0}^*h,
\)
comparison of the \(u_1\)- and \(u_2\)-coefficients, yields
\[
M=\sum_{j\in I}p_jm_j.
\]
Indeed, divisors of the first type have zero relative coefficients, while
for \(Z_j=\Gamma_{2,a_0}\) the two relative coefficients sum to
\(-1+1=0\).  Since \(M>0\), it follows that \(I\neq\varnothing\).

Comparing the coefficients of \(\pi_{2,0}^*h\), we obtain
\[
 T=\sum_jp_jt_j\leqslant\sum_{j\in I}p_jt_j,
\]
since the first two cases contribute only nonpositive $t_j$.
Consequently,
\[
 \max_{j\in I}\frac{t_j}{m_j}
 \geqslant\frac{\sum_{j\in I}p_jt_j}{\sum_{j\in I}p_jm_j}
 \geqslant\frac TM>\frac15.
\]
Choose a component $Z=Z_j$ attaining this maximum and put
$b_1=b_{1,j}$, $b_2=b_{2,j}$, $m=b_1+b_2$ and $t=t_j$.
Then $m,t>0$, $b_1\geqslant2b_2\geqslant0$ and $t/m>1/5$.
The canonical section of $\mathcal O_{X_{2,a_0}}(Z)$, after choosing
an isomorphism with $\mathscr L_{b_1,b_2,t}|_{X_{2,a_0}}$, gives
\[
 \widetilde\omega_{a_0}\in
 H^0(X_{2,a_0},\mathscr L_{b_1,b_2,t}|_{X_{2,a_0}}),
 \qquad \operatorname{div}(\widetilde\omega_{a_0})=Z.
\]
If $b_2=0$, then $b_1=m$.  Pushing this nonzero section down the two
projective bundles and applying the projection formula would give
\[
 0\neq H^0\!\left(\mathbb P^2,
 S^mT_{\mathbb P^2}^*(\log\mathcal C_{a_0})
 \otimes\mathcal O_{\mathbb P^2}(-t)\right),
\]
contrary to Proposition~\ref{pro:one-jet-vanishing}.
Hence $b_1\geqslant2b_2>0$.

We now extend the section defining $Z$.
Set $\mathscr L=\mathscr L_{b_1,b_2,t}$.
By the choice of $a_0$, the sheaf $q_*\mathscr L$ is locally free
near $a_0$ and commutes with base change there.
After choosing a local algebraic frame, the vector
$\widetilde\omega_{a_0}$ extends to a section of $q_*\mathscr L$
on a Zariski-open neighbourhood $U$ of $a_0$, hence to
\[
 \widetilde\omega\in H^0(q^{-1}(U),\mathscr L).
\]
Since 
\(\widetilde\omega|_{X_{2,a_0}}=\widetilde\omega_{a_0}\neq0\), the locus of points \(a\in U\) for which
\(\widetilde\omega|_{X_{2,a}}\) vanishes identically is Zariski closed in \(U\).
After replacing \(U\) by a smaller Zariski-open neighbourhood of
\(a_0\), we may assume that
\[
\widetilde\omega|_{X_{2,a}}\neq0
\qquad\text{for every }a\in U.
\]
Since each fibre \(X_{2,a}\) is smooth and geometrically integral,
these restrictions are nonzerodivisors.  The fibrewise criterion for relative effective
Cartier divisors therefore shows that
\[
 \mathscr Z:=V(\widetilde\omega)\subset q^{-1}(U)
\]
is an effective Cartier divisor flat over $U$.  It is also proper and of finite
presentation over $U$.  Its fibre at $a_0$ is the irreducible reduced
divisor $Z$, which is geometrically integral over $\mathbb C$.
The openness of the geometrically integral fibre locus for proper flat
morphisms of finite presentation
\cite[Theorem~12.2.1(x)]{EGAIV3} now allows us to shrink $U$ around
$a_0$ so that every $Z_a:=\mathscr Z_a$ is geometrically integral.
In particular, for every complex parameter $a\in U$, the divisor
$Z_a$ is irreducible and reduced, with
\[
 Z_a\sim b_1u_{1,a}+b_2u_{2,a}-t\pi_{2,0}^*h.
\]
It cannot equal $\Gamma_{2,a}$: their classes in the projective-bundle
Picard group differ, since $\Gamma_{2,a}$ has class
$-u_{1,a}+u_{2,a}$ and $t>0$.

Finally, the relative form of~\eqref{gamma_2=O(-1, 1)} is
\[
 \mathcal O_{\mathscr X_2}(\varGamma_2)
 \simeq\varPi_{2,1}^*\mathcal O_{\mathscr X_1}(-1)
 \otimes\mathcal O_{\mathscr X_2}(1).
\]
Let $s_{\varGamma_2}$ be its canonical section and set
$\omega_1=s_{\varGamma_2}^{\,b_1}\widetilde\omega$ on $q^{-1}(U)$.
Since $m=b_1+b_2$, this is an algebraic section of
$\mathcal O_{\mathscr X_2}(m)\otimes
\pi_{2,0}^*\mathcal O_{\mathbb P^2}(-t)$ and satisfies
\[
 \operatorname{div}(\omega_{1,a})=Z_a+b_1\Gamma_{2,a}
 \qquad(a\in U).
\]
The relative direct-image correspondence underlying
\eqref{direct image formula} identifies its restrictions with an
algebraic, hence holomorphic, family
\begin{equation}\label{define omega_1}
 0\neq\omega_{1,a}\in H^0\!\left(
 \mathbb P^2,E_{2,m}T_{\mathbb P^2}^*(\log\mathcal C_a)
 \otimes\mathcal O_{\mathbb P^2}(-t)\right),
 \qquad \frac tm>\frac15.
\end{equation}
Taking $U^{m,t}=U$ proves the proposition.
\end{proof}

\section{\bf Generating Independent Twisted  Jet Differentials via Slanted Vector Fields}\label{section: Siu's fields}

To construct an additional family of independent negatively twisted jet differentials, we employ Siu's method of slanted vector fields~\cite{Siu2004}. This approach generalizes to $k$-jet spaces the original technique developed for $1$-jet spaces by Clemens~\cite{Clemens1986}, Ein \cite{Ein1988}, and Voisin~\cite{Voisin1996} in their study of rational curves.

The method was first implemented by P\u{a}un~\cite{mihaipaun2008} for compact surfaces in $\mathbb{P}^{3}$. Rousseau~\cite{Rousseau2009} subsequently established the logarithmic theory for pairs $(\mathbb{P}^{2},D)$ with curve configurations $D$, while in~\cite{Rousseau2007cpt,Rousseau2007log} he further extended both P\u{a}un's setting and the logarithmic case to threefolds. For hypersurfaces $H\subset\mathbb{P}^n$ of arbitrary dimension, Merker~\cite{Merker2009} established the complete generalization, thereby realizing Siu's original vision concerning generic global generation. The logarithmic counterpart for pairs $(\mathbb{P}^n,H)$ was subsequently obtained by Darondeau~\cite{Darondeau_vectorfield}.

This section adapts the standard technique of Siu's slanted vector fields to the case of three generic conics in $\mathbb{P}^2$. The content is entirely classical and no originality is claimed.

\subsection{Logarithmic Slanted Vector Fields for Three Conics in $\mathbb{P}^2$}\label{subsec:groundwork of slanted vector fields}

The framework presented in this subsection specializes Rousseau's logarithmic setting in $\mathbb{P}^2$ to the configurations of three conics. We emphasize that this subsection contains no original contributions and is included solely for completeness and the reader's convenience. 

Let
\begin{equation}\label{family of curves in P^2}
    \mathcal{X} = \cup_{\lambda=1}^3 \mathcal{X}_{\lambda} \subset \mathbb{P}^2 \times {\textstyle\prod}_{\lambda=1}^3 \mathbb{P}^{N_{d_{\lambda}}}
\end{equation}
denote the universal family of algebraic curves in \(\mathbb{P}^2\) consisting of three components of degrees \(d_1\), \(d_2\), and \(d_3\), respectively. Here, each \(\mathcal{X}_{\lambda}\) is the universal curve of degree \(d_{\lambda}\), defined by the equation $\sum_{|\alpha| = d_{\lambda}} A_{\alpha}^{\lambda} Z^{\alpha}=0$, where \([A^{\lambda}] = [A_{\alpha}^{\lambda}]_{|\alpha|=d_{\lambda}} \in \mathbb{P}^{N_{d_{\lambda}}}\) are the homogeneous coefficients and \([Z] = [Z_0 : Z_1 : Z_2] \in \mathbb{P}^2\) are the projective coordinates. For a multi-index \(\alpha = (\alpha_0, \alpha_1, \alpha_2) \in \mathbb{N}^3\) with length \(|\alpha| := \alpha_0 + \alpha_1 + \alpha_2\), we employ the notation
$
Z^{\alpha} := Z_0^{\alpha_0} Z_1^{\alpha_1} Z_2^{\alpha_2}$. 
The dimension of the parameter space for each component is given by
$
N_{d_{\lambda}} = \binom{d_{\lambda} + 2}{2} - 1$.

For any integer $k \geqslant 1$, the \emph{vertical logarithmic $k$-jet space} is then defined as 
\[
\overline{J_k^v}\left(\mathbb{P}^2 \times {\textstyle\prod}_{\lambda=1}^3 \mathbb{P}^{N_{d_\lambda}}\right)
:= \ker\left(\mathrm{pr}_2\right)_* 
\subset J_k\left(\mathbb{P}^2 \times {\textstyle\prod}_{\lambda=1}^3 \mathbb{P}^{N_{d_\lambda}}, -\log \mathcal{X} \right),
\]
where $\mathrm{pr}_2 \colon \mathbb{P}^2 \times {\textstyle\prod}_{\lambda=1}^3 \mathbb{P}^{N_{d_\lambda}} \longrightarrow {\textstyle\prod}_{\lambda=1}^3 \mathbb{P}^{N_{d_\lambda}}$ denotes the natural projection onto the parameter space and $J_k \big(\mathbb{P}^2 \times {\textstyle\prod}_{\lambda=1}^3 \mathbb{P}^{N_{d_\lambda}}, -\log \mathcal{X} \big)$ is the logarithmic $k$-jet bundle along the universal family $\mathcal{X}$.

We define the universal intersection family
\[
\mathcal{Y} = \cap_{\lambda=1}^3 \mathcal{Y}_{\lambda} \subset \mathbb{P}^5 \times {\textstyle\prod}_{\lambda=1}^3 \mathbb{P}^{N_{d_{\lambda}} + 1},
\]
where for each \(1 \leqslant \lambda \leqslant 3\), the hypersurface \(\mathcal{Y}_{\lambda}\) is the universal family of degree \(d_\lambda\) hypersurfaces in \(\mathbb{P}^5\), given explicitly by $\mathcal{Y}_{\lambda} := \left\{ 
A_{\mathbf{0}}^{\lambda} Z_{2+\lambda}^{d_{\lambda}} + \sum_{|\alpha| = d_{\lambda}} A_{\alpha}^{\lambda} Z^{\alpha} = 0 
\right\} \subset \mathbb{P}^5 \times \mathbb{P}^{N_{d_{\lambda}} + 1}$. Here, the space \(\mathbb{P}^5\) has homogeneous coordinates $[Z : Z_3 : Z_4 : Z_5] = [Z_0 : \cdots : Z_5]$, and the parameter space \(\mathbb{P}^{N_{d_{\lambda}} + 1}\) is equipped with homogeneous coordinates $[A^{\lambda} : A_{\mathbf{0}}^{\lambda}] = [\ldots : A_{\alpha}^{\lambda} : \ldots : A_{\mathbf{0}}^{\lambda}]_{|\alpha| = d_{\lambda}}$.

For notational convenience, we denote the exceptional locus $\{Z = \mathbf{0}\} \cup \cup_{\lambda=1}^{3} \{A^{\lambda} = \mathbf{0}\}$, where $\{Z = \mathbf{0}\}$ denotes the subset of $\mathbb{P}^5$ defined by $Z_0 = Z_1 = Z_2 = 0$ and $\{A^{\lambda} = \mathbf{0}\}$ refers to the locus in $\mathbb{P}^{N_{d_{\lambda}} + 1}$ where all coordinates vanish except $A_{\mathbf{0}}^{\lambda}$. There exists a natural forgetful map
\[
\pi \colon \mathbb{P}^{5} \times {\textstyle\prod}_{\lambda=1}^3 \mathbb{P}^{N_{d_{\lambda}} + 1} \setminus \left(\{Z = \mathbf{0}\} \cup \cup_{\lambda=1}^{3} \{A^{\lambda} = \mathbf{0}\}\right) \longrightarrow \mathbb{P}^2 \times {\textstyle\prod}_{\lambda=1}^3 \mathbb{P}^{N_{d_{\lambda}}}
\]
that forgets the auxiliary coordinates. Explicitly, $\pi$ is given by
\[
\left([Z:Z_3:Z_4:Z_5], [A^1:A_{\mathbf{0}}^1], [A^2:A_{\mathbf{0}}^2], [A^3:A_{\mathbf{0}}^3]\right) \longmapsto \left([Z], [A^1], [A^2], [A^3]\right).
\]

We observe the containment relation
\begin{equation}
    \label{containment relation}
\mathcal{Y} \cap \left( \{ Z = \mathbf{0} \} \cup \cup_{\lambda=1}^{3} \{ A^{\lambda} = \mathbf{0} \} \right) 
\subset \mathcal{Y} \cap \left( \cup_{\lambda=1}^{3} \left( \{ A_{\mathbf{0}}^{\lambda} = 0 \} \cup \{ A^{\lambda} = \mathbf{0} \} \right) \right),
\end{equation}
from the inclusion $\mathcal{Y} \cap \{ Z = \mathbf{0} \} \subset \mathcal{Y} \cap \cup_{\lambda=1}^{3} \{ A_{\mathbf{0}}^{\lambda} = 0 \}$, which is deduced by the definition of $\mathcal{Y}$.

Now we consider $\mathcal{Y}^* \coloneqq \mathcal{Y} \setminus \left( \cup_{\lambda=1}^{3} \left( \{ A_{\mathbf{0}}^{\lambda} = 0 \} \cup \{ A^{\lambda} = \mathbf{0} \} \right) \right)$ on which the restricted forgetful map \[\pi|_{\mathcal{Y}^*} \colon \mathcal{Y}^* \rightarrow \mathbb{P}^2 \times {\textstyle\prod}_{\lambda=1}^3 \mathbb{P}^{N_{d_\lambda}}\] is well-defined according to the
containment relation~\eqref{containment relation}. Furthermore, we obtain the identification of the preimage $\left( \pi|_{\mathcal{Y}^*} \right)^{-1}(\mathcal{X}) = \{ Z_3 Z_4 Z_5 = 0 \} =: \mathcal{D}$.
This shows that \(\pi\) induces a logarithmic morphism from $\left( \mathcal{Y}^*, \mathcal{D} \right)$ to $\big( \mathbb{P}^2 \times {\textstyle\prod}_{\lambda=1}^3 \mathbb{P}^{N_{d_\lambda}}, \mathcal{X} \big)$, which consequently gives rise to a dominant map between jet spaces:
\begin{equation}\label{pi_[2]}
    \pi_{[2]} \colon \overline{J_2^v}(\mathcal{Y}^*) \longrightarrow \overline{J_2^v}\left( \mathbb{P}^2 \times {\textstyle\prod}_{\lambda=1}^3 \mathbb{P}^{N_{d_\lambda}} \right).
\end{equation}
Here, $\overline{J_2^v}(\mathcal{Y}^*)$ represents the space of vertical logarithmic 2-jets with respect to the normal crossing divisor $\mathcal{D} = \{ Z_3 Z_4 Z_5 = 0 \}$. This construction allows us to perform all local computations on the more tractable space $\mathcal{Y}^*$ through pullback via $\pi$, where we have access to explicit coordinate expressions.

On the affine chart $\{Z_0 \neq 0\} \times {\textstyle\prod}_{\lambda=1}^3 \{A_{\mathbf{0}}^{\lambda} \neq 0\}$, we employ the affine coordinates $z_i := Z_i/Z_0$ for $i = 1,\dots,5$ and $a_{\alpha}^{\lambda} := A_{\alpha_0 \alpha_1 \alpha_2}^{\lambda}/A_{\mathbf{0}}^{\lambda}$ with $\alpha_0 = d_{\lambda} - \alpha_1 - \alpha_2$.  To streamline notation, for each jet of order $j = 1,2$ we define:
\begin{equation*}\label{notation xi of 2-jets}
    \xi_i^{(j)} := \begin{cases}
        z_i^{(j)} & \text{for } i = 1,2, \\
        (\log z_{2+\lambda})^{(j)} & \text{for } i = 3,4,5 \ (\lambda = i-2).
    \end{cases}
\end{equation*} 
Then local coordinates for $\overline{J_2}\big( \{ Z_0 \neq 0 \} \times {\textstyle\prod}_{\lambda=1}^3 \{ A_{\mathbf{0}}^{\lambda} \neq 0 \} \big) \subset J_2\big( \mathbb{P}^5 \times {\textstyle\prod}_{\lambda=1}^3 \mathbb{P}^{N_{d_\lambda} + 1} , -\log \mathcal{D} \big)$ are given by the coordinate tuple $\left( z_i , a_{\alpha}^{\lambda} , \xi_{i_1}' , \xi_{i_2}'' \right)$.

\medskip

Following the strategy of Siu \cite{Siu2004}
and the realizations of Păun~\cite{mihaipaun2008} and Rousseau~\cite{Rousseau2009}, we consider the following family of Siu's slanted vector fields with  pole order $7$:

\begin{thm}\label{thm:globally generated vector field of pole order 7}
Let $\Sigma_0$ be the zero locus of the $2 \times 2$ Wronskian $\det ( \xi_i^{(j)} )_{1 \leqslant i,j \leqslant 2}$:
\begin{equation*}
    \Sigma_0 = \left\{ \left( z_i , a_{\alpha}^{\lambda} , \xi_{i_1}' , \xi_{i_2}'' \right) \,\middle\vert\, \det \left( \xi_i^{(j)} \right)_{1 \leqslant i,j \leqslant 2} = 0 \right\} \subset \overline{J_2^v} \left( \mathcal{Y} \cap \left( \{ Z_0 \neq 0 \} \times {\textstyle\prod}_{\lambda=1}^3 \{ A_{\mathbf{0}}^\lambda \neq 0 \} \right) \right),
\end{equation*}
and $\Sigma$ be the closure of $\Sigma_0$ in $\overline{J_2^v}(\mathcal{Y}^*)$. Then for some fixed $q_1 , q_2 , q_3 \in \mathbb{Z}_{>0}$, the twisted tangent bundle
\begin{equation}\label{the twisted tangent bundle}
    T_{ \overline{J_2^v}\left( \mathbb{P}^2 \times {\textstyle\prod}_{\lambda=1}^3 \mathbb{P}^{N_{d_\lambda}} \right) } \otimes p^* \mathrm{pr}_1^* \mathcal{O}_{\mathbb{P}^2}(7) \otimes p^* \mathrm{pr}_2^* \mathcal{O}_{\prod_{\lambda=1}^3 \mathbb{P}^{N_{d_\lambda}}}(q_1,q_2,q_3), 
\end{equation}
is generated by its global sections at every point of $\overline{J_2^v} ( \mathbb{P}^2 \times {\textstyle\prod}_{\lambda=1}^3 \mathbb{P}^{N_{d_\lambda}} ) \setminus \left( \pi_{[2]} \left( \Sigma \right) \cup p^{-1} \left( \mathcal{X} \right) \right)$, where $\pi_{[2]}$ is introduced in \eqref{pi_[2]}, $p$ is the natural projection 
\begin{equation*}
    p 
    \colon 
    \overline{J_2^v} \left( \mathbb{P}^2 \times {\textstyle\prod}_{\lambda=1}^3 \mathbb{P}^{N_{d_\lambda}} \right) 
    \longrightarrow 
    \mathbb{P}^2 \times {\textstyle\prod}_{\lambda=1}^3 \mathbb{P}^{N_{d_\lambda}},
\end{equation*}
and $\mathrm{pr}_1$ (resp. $\mathrm{pr}_2$) is the projection from $\mathbb{P}^2 \times {\textstyle\prod}_{\lambda=1}^3 \mathbb{P}^{N_{d_\lambda}}$ to $\mathbb{P}^2$ (resp. ${\textstyle\prod}_{\lambda=1}^3 \mathbb{P}^{N_{d_\lambda}}$). \qed
\end{thm}

\begin{rem}\rm
By the construction and proof of
\cite[Corollary~8 and Proposition~10]{Rousseau2009}, together with
\cite[Section~4.3]{Merker2009}, the vector fields in
Theorem~\ref{thm:globally generated vector field of pole order 7}
are $\mathbb G_2$-invariant (see~\eqref{reparameterization group G_k}). Consequently, their jet-direction coefficients
carry exactly the weights compensating those lost under differentiation, and
hence their action preserves the weighted degree in the jet variables.

Moreover, as in \cite[Remark~10]{Rousseau2009}, if the image of the $2$-jet
lift $j_2(f)$ is contained in $\pi_{[2]}(\Sigma)$, then the image of $f$ is
contained in a line and hence is algebraically degenerate. Therefore, for an
algebraically nondegenerate entire curve $f$, its $2$-jet lift is not entirely
contained in $\pi_{[2]}(\Sigma)$, and it suffices to work on the nonempty open
set where the lifted jet lies outside $\pi_{[2]}(\Sigma)$.
\end{rem}

\medskip

\subsection{Producing an Independent $2$-Jet via Siu's Slanted Vector Fields of Pole Order~$7$}

\begin{pro}[The second global jet differential]
\label{prop:global-slanted-second-jet}
Let $P \coloneqq \{\omega_{1,a}\}_{a\in U}$ be the holomorphic family supplied by
Proposition~\ref{prop:step1}, with weighted degree $m$ and negative base
twist $t\geqslant8$.  For a generic parameter $a$, one can choose a
pole-order-$7$ slanted vector field $V$ such that
\[
 \omega_{2,a}:=\left.\mathcal L_VP\right|_a
 \in H^0\!\left(\mathbb P^2,
 E_{2,m}T_{\mathbb P^2}^*(\log\mathcal C_a)\otimes
 \mathcal O_{\mathbb P^2}(7-t)\right)
\]
is nonzero and the common zero locus of $\omega_{1,a}$ and $\omega_{2,a}$
in $X_{2,a}$ has codimension at least $2$ modulo $\Gamma_{2,a}$.
Consequently, if $f^*\omega_{1,a}\equiv0$, then either
\[
 T_f(r)\leqslant\frac{m}{t-7}
 \sum_{i=1}^3N_f^{[1]}(r,\mathcal C_i)
 +O\bigl(\log T_f(r)+\log r\bigr)\quad\parallel,
\]
or $f_{[1]}(\mathbb C)$ is algebraically degenerate in $X_{1,a}$.
\end{pro}

This is Siu's global differentiation mechanism~\cite{Siu2004}, in the
logarithmic form developed by Rousseau~\cite{Rousseau2009}; see also the
explicit invariant formulation in
\cite[pp.~171--173]{DiverioMerkerRousseau2010}. On the universal vertical jet space, $P$ is locally expressed as an
equivariant polynomial function of the jet variables, and a
slanted vector field acts on it naturally by derivation.  Equivariance makes
the result independent of the homogeneous chart and preserves the weighted
degree and reparameterization invariance.

\begin{proof}
Regard the family $\{\omega_{1,a}\}_{a\in U}$ as the equivariant holomorphic
map
\[
 P:\overline{J_2^v}\!\left(\mathbb P^2\times
 {\textstyle\prod}_{\lambda=1}^3\mathbb P^{N_{d_\lambda}}\right)
 \big|_{\mathrm{pr}_2^{-1}(U)}
 \longrightarrow
 p^*\mathrm{pr}_1^*\mathcal O_{\mathbb P^2}(-t).
\]
If $V$ is a global invariant section of the twisted tangent bundle
\eqref{the twisted tangent bundle}, its natural action on equivariant
functions gives
\[
 \mathcal L_VP:\overline{J_2^v}\!\left(\mathbb P^2\times
 {\textstyle\prod}_{\lambda=1}^3\mathbb P^{N_{d_\lambda}}\right)
 \big|_{\mathrm{pr}_2^{-1}(U)}
 \longrightarrow
 p^*\mathrm{pr}_1^*\mathcal O_{\mathbb P^2}(7-t),
\]
where the parameter-space twist becomes a constant one-dimensional factor
after restriction to a fixed fibre.  Since $V$ is invariant, $\mathcal L_VP$
is again an invariant polynomial of weighted degree $m$ in the jet
variables.  Its restriction to the fibre over $a$ is therefore the asserted
global differential $\omega_{2,a}$.

By Proposition~\ref{prop:step1}, the residual zero divisor of
$\omega_{1,a}$ is irreducible and reduced modulo $\Gamma_{2,a}$.  The global
generation in Theorem~\ref{thm:globally generated vector field of pole order 7}
allows us to choose $V$ non-tangent to that divisor at a general smooth
point.  Hence $\mathcal L_VP$ does not vanish identically there, and the
common zero locus has codimension at least $2$ modulo $\Gamma_{2,a}$.

If $f^*\omega_{2,a}\not\equiv0$, apply
Theorem~\ref{smt-form-logarithmic-diff-jet} to $\omega_{2,a}$, whose negative
twist is $t-7$.  If both pullbacks vanish, then $f_{[2]}(\mathbb C)$ lies in
their common zero locus.  Its image in $X_{1,a}$ has dimension at most $2$,
strictly smaller than $\dim X_{1,a}=3$; hence $f_{[1]}$ is algebraically
degenerate.
\end{proof}

%\begin{rmk}\label{how to use O(3)}\rm 
%We hypothesize that an alternative strategy, which deepens the technical machinery of Paŭn \cite[pp.~886-889]{mihaipaun2008} and Rousseau \cite[Proposition 17]{Rousseau2009} via a refined application of Siu's slanted vector fields with twists by $\mathcal{O}_{\mathbb{P}^2}(3)$, $\mathcal{O}_{\mathbb{P}^2}(4)$, $\mathcal{O}_{\mathbb{P}^2}(5)$, or $\mathcal{O}_{\mathbb{P}^2}(6)$, could demonstrate the nonexistence of nontrivial solutions for a certain associated ODE. Achieving this would substantially lower the number of required vanishing results. Nevertheless, in the setting of three conics, the substantial algebraic complexity encountered in the ODE analysis precluded further progress.
%\end{rmk}

\section{\bf Refinement of Demailly--El Goul's Arguments}\label{sect: refined Demailly-EG}

\subsection{More Riemann-Roch Computations}
\label{subsection:4.1}

In Proposition~\ref{prop:step1}, we show the existence of a nonzero negatively twisted invariant logarithmic $2$-jet differential
\[
\omega_1 \in H^0\big(X_2, \mathcal{O}_{X_2}(m) \otimes \pi_{2,0}^* \mathcal{O}_X(-t)\big),
\]
for some positive integers $m, t \geqslant 1 $ with $t/m\geqslant 1/5$, whose zero divisor on the second level of the Demailly-Semple tower $X_2$ is of the form $Z + b_1 \Gamma_2$, where $Z\subset X_2$ is irreducible and reduced. Moreover, we have the linear equivalence
\begin{equation}\label{equivalent class of Z}
    Z \sim b_1 u_1 + b_2 u_2 - t \pi_{2,0}^* h, \quad (b_1, b_2) \in \mathbb{Z}^2, \quad b_1 \geqslant 2b_2 > 0, \quad b_1 + b_2 = m, \quad t \in \mathbb{Z}_{+},
\end{equation}
where $u_1 = \pi_{2,1}^* \mathcal{O}_{X_1}(1)$, $u_2 = \mathcal{O}_{X_2}(1)$, and $h = \mathcal{O}_{\mathbb{P}^2}(1)$. 

Hereafter, we use the standard notation:
\begin{equation*}
    \mathcal{O}_{X_{2}}(b_1 , b_2) \coloneqq \pi_{2,1}^{*} \mathcal{O}_{X_{1}}(b_1) \otimes \mathcal{O}_{X_{2}}(b_2), \quad b_1, b_2\in \mathbb{Z}.
\end{equation*}

Using the same argument as in Demailly--El Goul~\cite[Proposition 3.4]{Demailly-Elgoul2000} --- Miyaoka's method~\cite{Miyaoka}; see also Lu--Yau~\cite{Lu-Yau} --- one can show that $\mathcal{O}_{X_2}(2,1)|_{Z}$ is big, provided that $t < m$. Now, our goal  is to determine an effective small rational number $\tau > 0$ such that the $\mathbb{Q}$-line bundle
$
\big( \mathcal{O}_{X_2}(2,1) \otimes \pi_{2,0}^* \mathcal{O}_{\mathbb{P}^2}(-\tau) \big)|_{Z}
$
remains big. This will be achieved through Riemann-Roch computations and the following proposition observed by Demailly~\cite[p.~73]{Demailly1997}.

\begin{pro}\label{pro:the higher direct images vanishing}
Let $a_1,a_2\in\mathbb Z$.
\begin{enumerate}
    \item If $a_2=-1$, then
    \[
    R^q(\pi_{2,0})_*\mathcal O_{X_2}(a_1,-1)=0
    \qquad\text{for every }q\geqslant0.
    \]

    \item If $a_2\geqslant0$, then
    \[
    R^q(\pi_{2,0})_*\mathcal O_{X_2}(a_1,a_2)=0
    \qquad\text{for every }q\geqslant2,
    \]
    and
    \[
    R^1(\pi_{2,0})_*\mathcal O_{X_2}(a_1,a_2)=0
    \quad\Longleftrightarrow\quad
    a_1-2a_2\geqslant-1.
    \]
\end{enumerate}
In particular, for every $p\geqslant0$,
\[
R^q(\pi_{2,0})_*\mathcal O_{X_2}(2p,p)=0
\qquad\text{for every }q\geqslant1.
\]
\end{pro}

\begin{proof}
Set
\[
\pi_{2,1}\colon X_2=\mathbb P(V_1)\longrightarrow X_1,
\qquad
\pi_{1,0}\colon X_1\longrightarrow X_0.
\]
Recall that $V_1$ is the rank-two directed bundle defined by
\eqref{SES:V_1}, and that $V_1^*$ denotes its dual. We use the
convention that $\mathbb P(V_1)$ parametrizes lines in the fibers of
$V_1$, so that $\mathcal O_{X_2}(-1)\subset\pi_{2,1}^*V_1$ is the
tautological line subbundle.

Suppose first that $a_2=-1$. The restriction of
$\mathcal O_{X_2}(a_1,-1)$ to every fiber of $\pi_{2,1}$ is
$\mathcal O_{\mathbb P^1}(-1)$, whose cohomology vanishes in every
degree. Cohomology and base change therefore gives
\[
R^q(\pi_{2,1})_*\mathcal O_{X_2}(a_1,-1)=0
\qquad(q\geqslant0),
\]
and the first assertion follows from $\pi_{2,0}=\pi_{1,0}\circ\pi_{2,1}$.

Assume now that $a_2\geqslant0$. By the projective bundle formula and
the projection formula,
\begin{equation}\label{eq:rho-pushforward-corrected}
R^q(\pi_{2,1})_*\mathcal O_{X_2}(a_1,a_2)=0
\quad(q\geqslant1),
\qquad
(\pi_{2,1})_*\mathcal O_{X_2}(a_1,a_2)
=
\mathcal O_{X_1}(a_1)\otimes S^{a_2}V_1^*.
\end{equation}
Hence the Leray spectral sequence yields
\[
R^q(\pi_{2,0})_*\mathcal O_{X_2}(a_1,a_2)
\simeq
R^q(\pi_{1,0})_*
\bigl(\mathcal O_{X_1}(a_1)\otimes S^{a_2}V_1^*\bigr).
\]

Fix $x\in X_0$ and set
\[
F_x:=\pi_{1,0}^{-1}(x)\simeq\mathbb P^1.
\]
Restricting \eqref{SES:V_1} to $F_x$ gives
\[
0\longrightarrow\mathcal O_{\mathbb P^1}(2)
\longrightarrow V_1|_{F_x}
\longrightarrow\mathcal O_{\mathbb P^1}(-1)
\longrightarrow0.
\]
This sequence splits, since
\[
\operatorname{Ext}^1
\bigl(\mathcal O_{\mathbb P^1}(-1),
      \mathcal O_{\mathbb P^1}(2)\bigr)
=
H^1\bigl(\mathbb P^1,\mathcal O_{\mathbb P^1}(3)\bigr)
=0.
\]
Thus
\[
V_1^*|_{F_x}
\simeq
\mathcal O_{\mathbb P^1}(-2)
\oplus
\mathcal O_{\mathbb P^1}(1),
\]
and, since
$\mathcal O_{X_1}(1)|_{F_x}\simeq\mathcal O_{\mathbb P^1}(1)$,
\begin{equation}\label{eq:fiber-splitting-corrected}
\left(
\mathcal O_{X_1}(a_1)\otimes S^{a_2}V_1^*
\right)\big|_{F_x}
\simeq
\bigoplus_{j=0}^{a_2}
\mathcal O_{\mathbb P^1}(a_1+a_2-3j).
\end{equation}

The smallest degree in \eqref{eq:fiber-splitting-corrected} is
$a_1-2a_2$. Since
\[
H^1\bigl(\mathbb P^1,\mathcal O_{\mathbb P^1}(d)\bigr)=0
\quad\Longleftrightarrow\quad
d\geqslant-1,
\]
the fiberwise first cohomology vanishes precisely when
$a_1-2a_2\geqslant-1$. The splitting type is independent of $x$, so
Grauert's theorem gives
\[
R^1(\pi_{2,0})_*\mathcal O_{X_2}(a_1,a_2)=0
\quad\Longleftrightarrow\quad
a_1-2a_2\geqslant-1.
\]

Finally, $\pi_{1,0}$ has one-dimensional fibers; hence
\[
R^q(\pi_{1,0})_*
\bigl(\mathcal O_{X_1}(a_1)\otimes S^{a_2}V_1^*\bigr)=0
\qquad(q\geqslant2).
\]
Together with \eqref{eq:rho-pushforward-corrected}, this proves the
remaining assertion. Taking $(a_1,a_2)=(2p,p)$ gives
$a_1-2a_2=0$, and therefore all higher direct images vanish.
\end{proof}

\begin{pro}\label{prop: bigness on Z}
If $0 < \frac{t}{m} < 1$ and $0 \leqslant \tau < \tau_1(\frac{t}{m}), \tau\in \mathbb{Q}$, where
 {\footnotesize
 \[
 \tau_1 \left( \frac{t}{m} \right) = \frac{3}{2} \left(4 - \frac{t}{m}\right) - \frac{\sqrt{3}}{2} \sqrt{32 - 8 \frac{t}{m} + 3 \left( \frac{t}{m} \right)^2},
 \]}
 then the restricted $\mathbb{Q}$-line bundle
$
\big( \mathcal{O}_{X_2}(2,1) \otimes \pi_{2,0}^* \mathcal{O}_{\mathbb{P}^2}(-\tau) \big)|_{Z}
$
is big on $Z$.
\end{pro}

\begin{proof}
We begin by recalling some key intersection formulas on the Demailly–Semple tower. By~\cite[p.~478]{Elgoul2003}, we have the identities:
\[
u_1^2 + \bar{c}_1 u_1 + \bar{c}_2 = 0 \quad \text{and} \quad u_2^2 + \bar{c}_1(V_1) u_2 + \bar{c}_2(V_1) = 0,
\]
where $\bar{c}_1(V_1) = \bar{c}_1 + u_1$ (see~\eqref{define V_1}) and $\bar{c}_2(V_1) = \bar{c}_2 - u_1^2 = 2\bar{c}_2 + \bar{c}_1 u_1$.

On the second level $X_2$ of the Demailly-Semple tower, the following identities hold:
\[
u_{1}^{4} = 0, \quad u_{1}^{3} u_{2} = \bar{c}_{1}^{2} - \bar{c}_{2}, \quad u_{1}^{2} u_{2}^{2} = \bar{c}_{2}, \quad u_{1} u_{2}^{3} = \bar{c}_{1}^{2} - 3\bar{c}_{2}, \quad u_{2}^{4} = 5\bar{c}_{2} - \bar{c}_{1}^{2}.
\]
Furthermore, for any pullback $F \in \pi_{2,0}^* \pic (\mathbb{P}^2)$, we have:
\[
u_{1}^{3} \cdot F = 0, \quad u_{1}^{2} u_{2} \cdot F = -\bar{c}_{1} \cdot F, \quad u_{1} u_{2}^{2} \cdot F = 0, \quad u_{2}^{3} \cdot F = 0.
\]
A crucial part of our approach relies on the following additional identities:
\[
u_{1}^2 \cdot F_1 \cdot F_2 = 0, \quad u_{1} u_{2} \cdot F_1 \cdot F_2 = F_1 \cdot F_2, \quad u_{2}^2 \cdot F_1 \cdot F_2 = - F_1 \cdot F_2,
\qquad
\forall\, F_1, F_2 \in \pi_{2,0}^* \pic (\mathbb{P}^2).
\]

Then by direct computation, we can obtain
\begin{equation}\label{equ:self intersection on Z}
    \left( 2 u_1 + u_2 - \tau \pi_{2,0}^* h \right)^3 \cdot Z = m\left( 13\overline{c}_{1}^{2} -9\overline{c}_{2}\right) + 12t(\overline{c}_{1} \cdot h) +12m(\overline{c}_{1} \cdot h) \tau +3m\left( h^{2}\right) \tau^{2} + 9\left( h^{2}\right) t \tau .
\end{equation}

In our setting of three generic conics, we have $\overline{c}_{1} = -3h$ and $\overline{c}_{2} = 9 h^2$. Substituting these into~\eqref{equ:self intersection on Z}, together with $h^2 = 1$, yields:
\begin{equation}\label{equ:self intersection on Z for 3 conics}
    \left( 2 u_1 + u_2 - \tau \pi_{2,0}^* h \right)^3 \cdot Z = 3 \left( m\tau ^{2} - 3(4m-t) \tau +12(m-t) \right).
\end{equation}
This quadratic polynomial has two real roots:
{\footnotesize
\begin{equation}
\label{what is tau}
    \tau_1 \left( \frac{t}{m} \right) = \frac{3}{2} \left(4 - \frac{t}{m}\right) - \frac{\sqrt{3}}{2} \sqrt{32 - 8 \frac{t}{m} + 3 \left( \frac{t}{m} \right)^2} 
    \,
    <\, \tau_2 \left( \frac{t}{m} \right) = \frac{3}{2} \left(4 - \frac{t}{m}\right) + \frac{\sqrt{3}}{2} \sqrt{32 - 8 \frac{t}{m} + 3 \left( \frac{t}{m} \right)^2},
\end{equation}
}
because $32 - 8 \frac{t}{m} + 3 ( \frac{t}{m} )^2 > 0$.  Moreover, $\tau_1 ( \frac{t}{m} ) > 0$ if and only if $\frac{t}{m} < 1$.

%For $0 < \frac{t}{m} < 1$, both functions $\tau_1(\frac{t}{m})$, $\tau_2(\frac{t}{m})$ are monotonically decreasing with respect to $\frac{t}{m}$, with $0 < \tau_1(\frac{t}{m}) < 6-2 \sqrt{6}$ and $9 < \tau_2(\frac{t}{m}) < 6+2 \sqrt{6}$. Furthermore, if $0 < \tau < \tau_1(\frac{t}{m})$ or $\tau > \tau_2(\frac{t}{m})$, then
%\begin{equation*}
    %\left( 2 u_1 + u_2 - \tau \pi_{2,0}^* h \right)^3 \cdot Z = 3 m \left( \tau - \tau_1\left(\frac{t}{m}\right) \right) \left( \tau - \tau_2\left(\frac{t}{m}\right) \right) > 0.
%\end{equation*}

Thus when $0<\frac{t}{m}<1$   and $0 \leqslant \tau < \tau_1(\frac{t}{m})$, we have $( 2 u_1 + u_2 - \tau \pi_{2,0}^* h )^3 \cdot Z > 0$.   Riemann-Roch Theorem then gives, for sufficiently divisible integer $p\gg 1$, that:
\begin{align*}
    &h^0 \big( Z, (\mathcal{O}_{X_2}(2,1) \otimes \pi_{2,0}^* \mathcal{O}_{\mathbb{P}^2}(-\tau))^{p} |_{Z} \big) + h^2 \big( Z, (\mathcal{O}_{X_2}(2,1) \otimes \pi_{2,0}^* \mathcal{O}_{\mathbb{P}^2}(-\tau))^{p} |_{Z} \big) \\
    \geqslant& 
    \chi \big( Z, (\mathcal{O}_{X_2}(2,1) \otimes \pi_{2,0}^* \mathcal{O}_{\mathbb{P}^2}(-\tau))^{p} |_{Z} \big) 
    =
    \frac{\big( 2 u_1 + u_2 - \tau \pi_{2,0}^* h \big)^3 \cdot Z}{3!} p^3 + O(p^2)\gg 1.
\end{align*}

Now we prove that
\begin{equation}\label{eq:h2-on-z-corrected}
H^2\!\left(
Z,
\left(\mathcal O_{X_2}(2,1)\otimes
\pi_{2,0}^*\mathcal O_{\mathbb P^2}(-\tau)\right)^p\big|_Z
\right)=0
\end{equation}
for all sufficiently divisible $p\gg1$.  Put
\[
\mathcal L_p=
\mathcal O_{X_2}(2p,p)\otimes
\pi_{2,0}^*\mathcal O_{\mathbb P^2}(-\tau p).
\]
Tensoring
$0\to\mathcal O_{X_2}(-Z)\to\mathcal O_{X_2}\to\mathcal O_Z\to0$
with $\mathcal L_p$ gives
\begin{align}\label{eq:h2-z-bound-corrected}
 h^2(Z,\mathcal L_p|_Z)
 \leqslant h^2(X_2,\mathcal L_p)
 +h^3(X_2,\mathcal L_p\otimes\mathcal O_{X_2}(-Z)).
\end{align}

For the first term, Proposition~\ref{pro:the higher direct images vanishing}
applies because $2p-2p=0\geqslant-1$.  Hence
\[
R^q(\pi_{2,0})_*\mathcal O_{X_2}(2p,p)=0\qquad(q\geqslant1)
\]
and
\[
(\pi_{2,0})_*\mathcal O_{X_2}(2p,p)
\simeq E_{2,3p}T_{\mathbb P^2}^*(\log\mathcal C).
\]
Indeed, in filtration layer $q$ a weighted-degree $3p$ invariant monomial
has $3p-3q$ first-derivative factors and $q$ Wronskian factors.  On $X_2$
its class is
\[
 \mathcal O_{X_2}(3p-q,q)
 =\mathcal O_{X_2}(3p)\bigl(-(3p-q)\Gamma_2\bigr).
\]
Since $0\leqslant q\leqslant p$, every section of
$\mathcal O_{X_2}(3p)$ coming from an invariant differential vanishes along
$\Gamma_2$ at least to order $2p$.  Hence
\[
 (\pi_{2,0})_*\mathcal O_{X_2}(2p,p)
 \simeq(\pi_{2,0})_*\mathcal O_{X_2}(3p)
 \simeq E_{2,3p}T_{\mathbb P^2}^*(\log\mathcal C),
\]
which justifies the displayed zeroth-direct-image identification.
By Leray and the projection formula,
\[
H^2(X_2,\mathcal L_p)
\simeq
H^2\!\left(
\mathbb P^2,
E_{2,3p}T_{\mathbb P^2}^*(\log\mathcal C)
\otimes\mathcal O_{\mathbb P^2}(-\tau p)
\right).
\]
The argument of Proposition~\ref{estimate of 2-jet threshold}, with
jet degree $3p$ and $\delta=\tau/9$, gives the vanishing of this group.
Indeed,
\[
0\leqslant\tau<\tau_1(t/m)\leqslant6-2\sqrt6<3,
\]
so that $0\leqslant\delta<1/3$.  Therefore
\begin{equation}\label{eq:first-ambient-vanishing-corrected}
H^2(X_2,\mathcal L_p)=0.
\end{equation}

It remains to treat the second term in
\eqref{eq:h2-z-bound-corrected}.  By~\eqref{equivalent class of Z}, set
\[
\mathcal M_p:=\mathcal L_p\otimes\mathcal O_{X_2}(-Z)
=\mathcal O_{X_2}(2p-b_1,p-b_2)
\otimes\pi_{2,0}^*\mathcal O_{\mathbb P^2}(t-\tau p).
\]
Write
\[
r:=b_1-2b_2\geqslant0.
\]
Choose $p$ so large that $p-b_2\geqslant0$ and, when $r\geqslant2$,
\[
p-b_2\geqslant\left\lfloor\frac{r-2}{3}\right\rfloor.
\]
Since $p-b_2\geqslant0$,
\[
R^q(\pi_{2,1})_*\mathcal M_p=0\quad(q\geqslant1),
\]
whereas $(\pi_{2,1})_*\mathcal M_p$ has a filtration induced by
\[
0\longrightarrow\mathcal O_{X_1}(1)
\longrightarrow V_1^*
\longrightarrow\pi_{1,0}^*\overline K_{\mathbb P^2}
       \otimes\mathcal O_{X_1}(-2)
\longrightarrow0.
\]
After reindexing by $j=p-b_2-k$, the corresponding graded terms are
\begin{equation}\label{eq:bad-piece-on-x1}
\mathcal Q_{p,k}=
\pi_{1,0}^*\!\left(
\overline K_{\mathbb P^2}^{\,p-b_2-k}
\otimes\mathcal O_{\mathbb P^2}(t-\tau p)
\right)
\otimes\mathcal O_{X_1}(-r+3k),
\qquad 0\leqslant k\leqslant p-b_2.
\end{equation}
Set
\[
A_{p,k}:=
\overline K_{\mathbb P^2}^{\,p-b_2-k}
\otimes\mathcal O_{\mathbb P^2}(t-\tau p),
\qquad
d_k:=-r+3k,
\]
so that
\[
\mathcal Q_{p,k}
=
\pi_{1,0}^*A_{p,k}\otimes\mathcal O_{X_1}(d_k).
\]
If $d_k\geqslant-1$, then the projective bundle formula and the
projection formula give
\[
R^1(\pi_{1,0})_*\mathcal Q_{p,k}
\simeq
A_{p,k}\otimes
R^1(\pi_{1,0})_*\mathcal O_{X_1}(d_k)
=0.
\]
Since $\pi_{1,0}$ has one-dimensional fibers and
$\dim\mathbb P^2=2$, the Leray spectral sequence yields
\[
H^3(X_1,\mathcal Q_{p,k})
\simeq
H^2\bigl(\mathbb P^2,R^1(\pi_{1,0})_*\mathcal Q_{p,k}\bigr)
=0.
\]
Consequently, only the indices satisfying $d_k\leqslant-2$ may
contribute; equivalently,
\[
0\leqslant k\leqslant
\left\lfloor\frac{r-2}{3}\right\rfloor.
\]

The only possible nonvanishing terms are therefore indexed by the fixed
finite set
\begin{equation}\label{eq:finite-bad-index-set}
0\leqslant k\leqslant
\left\lfloor\frac{r-2}{3}\right\rfloor.
\end{equation}
For such $k$, relative Serre duality for
$\pi_{1,0}\colon X_1=\mathbb P(T_{\mathbb P^2}(-\log\mathcal C))\to\mathbb P^2$
gives
\begin{equation}\label{eq:r1-bad-piece-corrected}
R^1(\pi_{1,0})_*\mathcal Q_{p,k}
\simeq 
S^{r-3k-2}T_{\mathbb P^2}(-\log\mathcal C)
\otimes
\mathcal O_{\mathbb P^2}
\bigl((3-\tau)p+t-3(b_2+k+1)\bigr),
\end{equation}
because $\overline K_{\mathbb P^2}\simeq\mathcal O_{\mathbb P^2}(3)$.
For each $k$ in~\eqref{eq:finite-bad-index-set}, the symmetric power in
\eqref{eq:r1-bad-piece-corrected} is independent of $p$, while the twist
tends to $+\infty$ since $\tau<3$.  Serre vanishing therefore gives
\[
H^2\bigl(\mathbb P^2,R^1(\pi_{1,0})_*\mathcal Q_{p,k}\bigr)=0
\qquad(p\gg1).
\]
The Leray spectral sequence for $\pi_{1,0}$ now yields
$H^3(X_1,\mathcal Q_{p,k})=0$ for every graded term of the filtration.
Induction along the filtration gives
\[
H^3(X_1,(\pi_{2,1})_*\mathcal M_p)=0.
\]
Using $R^q(\pi_{2,1})_*\mathcal M_p=0$ once more, we obtain
\begin{equation}\label{eq:second-ambient-vanishing-corrected}
H^3(X_2,\mathcal M_p)=0.
\end{equation}
Combining~\eqref{eq:h2-z-bound-corrected},
\eqref{eq:first-ambient-vanishing-corrected}, and
\eqref{eq:second-ambient-vanishing-corrected} proves
\eqref{eq:h2-on-z-corrected}.  Thus Riemann--Roch gives
\[
h^0(Z,\mathcal L_p|_Z)
\geqslant
\frac{(2u_1+u_2-\tau\pi_{2,0}^*h)^3\cdot Z}{3!}\,p^3
+O(p^2),
\]
and the positive leading coefficient proves that the restricted
$\mathbb Q$-line bundle is big.
\end{proof}

\subsection{A Refined Estimate}
Here is a key step in proving our main theorem.

\begin{lem}[Tautological evaluation of a restricted section]
\label{lem:tautological-evaluation-on-Z}
Let $Z\subset X_2$ be an effective Cartier divisor, let
$p,\widetilde t\in\mathbb Z_{>0}$, and let
\[
 \sigma\in H^0\!\left(
 Z,\mathcal O_{X_2}(2p,p)\otimes
 \pi_{2,0}^*\mathcal O_{\mathbb P^2}(-\widetilde t)\big|_Z
 \right).
\]
Suppose that $g\colon\mathbb C\to\mathbb P^2$ is nonconstant,
$g(\mathbb C)\not\subset\mathcal C$, and that its logarithmic Semple lift
extends across its isolated indeterminacy points and satisfies
$g_{[2]}(\mathbb C)\subset Z$.  On
$\mathbb C\setminus g^{-1}(\mathcal C)$, regard
$dg_{[1]}\colon T_{\mathbb C}\to g_{[1]}^*V_1$ as a differential in the
directed structure.  Wherever it is nonzero, its image is the tautological
line $g_{[2]}^*\mathcal O_{X_2}(-1)$.  It therefore defines, by meromorphic
continuation across the remaining discrete set, the intrinsic tautological
derivative
\[
 \mathcal Dg_{[1]}\in
 H^0_{\mathrm{mer}}\!\left(
 \mathbb C,K_{\mathbb C}\otimes
 g_{[2]}^*\mathcal O_{X_2}(-1)
 \right).
\]
Here \(H^0_{\mathrm{mer}}(Y,L)\) denotes the vector space of global
meromorphic sections of a line bundle \(L\) on a complex manifold \(Y\).
The canonical inclusion
\[
 \iota_p\colon
 \mathcal O_{X_2}(2p,p)
 =\mathcal O_{X_2}(3p)(-2p\Gamma_2)
 \hookrightarrow\mathcal O_{X_2}(3p)
\]
defines an intrinsic meromorphic section
\begin{equation}\label{eq:canonical-restricted-evaluation}
 s_{\sigma,g}:=
 \left\langle
 g_{[2]}^*(\iota_p\sigma),
 (\mathcal Dg_{[1]})^{\otimes3p}
 \right\rangle
 \in H^0_{\mathrm{mer}}\!\left(
 \mathbb C,K_{\mathbb C}^{\otimes3p}\otimes
 g^*\mathcal O_{\mathbb P^2}(-\widetilde t)
 \right),
\end{equation}
This construction is independent of every
choice of an ambient local lifting of $\sigma$.  If $s_{\sigma,g}\not\equiv0$,
then its pole divisor satisfies
\begin{equation}\label{eq:canonical-restricted-pole-bound}
 \operatorname{div}_{\infty}(s_{\sigma,g})
 \leqslant3p\sum_{i=1}^3(g^*\mathcal C_i)_{\mathrm{red}},
\end{equation}
and, outside a set of finite Lebesgue measure,
\begin{equation}\label{eq:canonical-restricted-proximity}
 \int_0^{2\pi}\log^+
 \|s_{\sigma,g}(re^{\sqrt{-1}\theta})\|\,
 \frac{d\theta}{2\pi}
 =O\!\left(\log T_g(r)+\log r\right).
\end{equation}
\end{lem}

\begin{proof}
Since the map $g_{[2]}$ factors through $Z$, one may
pull back $\iota_p\sigma$ directly from $Z$, without extending $\sigma$ to a
neighborhood.  Pairing the two line bundles gives
\eqref{eq:canonical-restricted-evaluation}.  If one nevertheless represents
$\sigma$ by ambient local liftings, two such liftings differ by a section in
the ideal of $Z$; their pullbacks by $g_{[2]}$ therefore agree.

Equip the line bundles on $X_2$ and $\mathbb P^2$ with smooth Hermitian
metrics.  Since $Z$ is projective, the norm of the holomorphic section
$\iota_p\sigma$ is bounded on $Z$.  Consequently, for a uniform constant
$C>0$,
\[
 \|s_{\sigma,g}\|
 \leqslant C\,\|\mathcal Dg_{[1]}\|^{3p}.
\]
This estimate uses only the intrinsic section on $Z$ and does not extend it
to a saturated inverse image in the Semple tower.

We next analyze the tautological derivative.  Near a point of
$g^{-1}(\mathcal C)$, choose logarithmic coordinates with
$\mathcal C=\{x_1\cdots x_q=0\}$, $q\leqslant2$.  The derivative \(\mathrm{d}\,g_{[1]}\) has both base and fiber components. Its base components are those of \(g'\); in the logarithmic boundary directions, they are $(x_i\circ g)'/(x_i\circ g)$. 
Since the extended first lift $g_{[1]}$ is holomorphic, its fiberwise projective-coordinate components are holomorphic. Thus $\mathcal Dg_{[1]}$ has at most a simple pole along each
$(g^*\mathcal C_i)_{\mathrm{red}}$. At an isolated zero of $g'$
outside $g^{-1}(\mathcal C)$, the regular Semple lift extends and the
evaluation has only a removable singularity or a zero.  If the critical point lies over $\mathcal C$, a logarithmic component such as $(x_i\circ g)'/(x_i\circ g)$ may still have a simple pole. It follows that \( \mathcal Dg_{[1]} \) has at most a simple pole along each $(g^*\mathcal C_i)_{\mathrm{red}}$. Taking its $3p$-th tensor power therefore gives \eqref{eq:canonical-restricted-pole-bound}.

Finally choose a finite logarithmic coordinate cover and, above it, a finite
cover by relatively compact Semple charts.  A partition of unity is used
only to estimate the norm of the already-defined global section.  On each
chart, the base components of $\mathcal Dg_{[1]}$ are logarithmic derivatives of
local coordinate functions.  A projective Semple coordinate $u$ is a rational
function on $X_1$ and, in a logarithmic frame, a ratio of two first
logarithmic derivatives.  Hence
\[
 T_{u\circ g_{[1]}}(r)=O\bigl(T_g(r)+\log r\bigr).
\]
For a partition function $\chi$ supported in a relatively compact Semple
chart, choose $a\notin\overline{u(\operatorname{supp}\chi)}$.  On this
support, $u-a$ is bounded above and away from zero, and
\[
 (u\circ g_{[1]})'
 =(u\circ g_{[1]}-a)
 \frac{(u\circ g_{[1]})'}{u\circ g_{[1]}-a}.
\]
The logarithmic-derivative lemma applied to $u\circ g_{[1]}-a$ controls the
vertical component by
\[
 O\bigl(\log^+T_g(r)+\log r\bigr).
\]
Summing this estimate and
the analogous horizontal estimates over the finite cover gives
\eqref{eq:canonical-restricted-proximity}; this is the partition-of-unity
argument of \cite[Thm.~3.1, (3.3)--(3.6)]{Huynh-Vu-Xie-2017}
with the preceding Semple-coordinate adaptation made explicit.
\end{proof}

\begin{lem}[Tautological evaluation of a restricted $\mathcal O_{X_2}(m)$-section]
\label{lem:tautological-evaluation-general}
Let $Z\subset X_2$ be a projective subvariety and let
\[
 \sigma\in H^0\!\left(
 Z,\mathcal O_{X_2}(m)\otimes
 \pi_{2,0}^*\mathcal O_{\mathbb P^2}(-s)|_Z
 \right),\qquad m,s\geqslant1.
\]
If $g_{[2]}(\mathbb C)\subset Z$, the tautological pairing defines an
intrinsic meromorphic section
\[
 \left\langle g_{[2]}^*\sigma,
 (\mathcal Dg_{[1]})^{\otimes m}\right\rangle
 \in H^0_{\mathrm{mer}}\!\left(
 \mathbb C,K_{\mathbb C}^{\otimes m}\otimes
 g^*\mathcal O_{\mathbb P^2}(-s)
 \right).
\]
If it is not identically zero, then
\[
 T_g(r)\leqslant\frac{m}{s}\sum_{i=1}^3
 N_g^{[1]}(r,\mathcal C_i)
 +O\bigl(\log T_g(r)+\log r\bigr)\quad\parallel.
\]
\end{lem}

\begin{proof}
The proof is the same intrinsic pairing and logarithmic/Semple-coordinate
estimate as in Lemma~\ref{lem:tautological-evaluation-on-Z}, with
$\mathcal O_{X_2}(m)$ paired directly against
$\mathcal O_{X_2}(-m)$.  Thus the reduced boundary pole divisor has
coefficient at most $m$; Poincar\'e--Lelong and Jensen's formula give the
stated inequality.
\end{proof}

\begin{thm}\label{pro:SMT on Z}
Let $g \colon \mathbb{C} \rightarrow \mathbb{P}^2$ be an algebraically nondegenerate entire curve whose lift $g_{[2]} \colon \mathbb{C} \rightarrow X_2$  satisfies $g_{[2]}(\mathbb{C}) \subset Z$. Suppose there exists a section
\begin{equation}
    \label{sigma of thm 5.3}
\sigma \in H^0\Big(Z, \big( \mathcal{O}_{X_2}(2,1)^{\tilde{m}} \otimes \pi_{2,0}^*\mathcal{O}_X(-\tilde{t}) \big)|_Z \Big) \cong H^0\Big(Z, \big( \mathcal{O}_{X_2}(3) \otimes \mathcal{O}(- 2 \Gamma_2) \big)^{\tilde{m}} \otimes \pi_{2,0}^*\mathcal{O}_X(-\tilde{t}) \big|_Z \Big),
\end{equation}
such that the canonical evaluation
$s_{\sigma,g}$ of Lemma~\ref{lem:tautological-evaluation-on-Z}, with
$p=\tilde m$, is not identically zero.
Then the following Second Main Theorem type estimate holds:
\[
T_{g}(r) \leqslant \frac{3\tilde{m}}{\tilde{t}} \sum_{i=1}^{3} N_{g}^{[1]}(r,\mathcal{C}_i) + o\big(T_g(r)\big) \quad\parallel.
\]
\end{thm}

The fixed double vertical factor is retained through the canonical inclusion
in Lemma~\ref{lem:tautological-evaluation-on-Z}.  The lemma, rather than an
identification of an arbitrary local lifting with a jet polynomial, supplies
the global meromorphic evaluation and the coefficient
$3\tilde m/\tilde t$.

\begin{proof}
%%On the Zariski-open set $U_{0 , C}$, the logarithmic $1$-forms $\left(\log a \right)'$ and $\left(\log b \right)'$ generate the logarithmic cotangent bundle $T_{\mathbb{P}^2}^{*} \left( \log \mathcal{C} \right)$. Thus we can choose the dual basis $\mathbf{e}_a$ and $\mathbf{e}_b$ of the logarithmic tangent bundle $T_{\mathbb{P}^2} \left( - \log \mathcal{C} \right)$. 

%%Let $f \colon (\mathbb{C},0) \ni t \mapsto f(t) \in U_{0 , C} \subset \mathbb{P}^2$ be a germ of holomorphic curve, and $f = (f_1,f_2)$ be the local representation of $f$ over the open set $U_{0 , C}$. Then we can write $f' = f_1' \partial / \partial x_1 + f_2' \partial / \partial x_2 = \left(\log a(f) \right)' \mathbf{e}_a + \left(\log b(f) \right)' \mathbf{e}_b$. 

%%The construction of Demailly-Semple tower yields $f'(t) \in \mathcal{O}_{X_1}(-1)_{f_{[1]}(t)}$ and $f_{[1]}'(t) \in \mathcal{O}_{X_2}(-1)_{f_{[2]}(t)}$. Moreover, by the definition of $\Gamma_{2}$ and $\pi_{1,0} \circ f_{[1]} = f$, we can conclude that $f_{[2]}(t) \in \Gamma_{2}$ if and only if ${\pi_{1,0}}_{*} f_{[1]}'(t) = f'(t) = 0$. 

Equip $\mathcal{O}_{X}(-1)$ with the standard norm, whose Chern curvature is
$-\omega_{\mathrm{FS}}$. By~\eqref{sigma of thm 5.3}, the canonical
evaluation $s_{\sigma,g}$ is a nonzero meromorphic section of
$K_{\mathbb C}^{\otimes3\tilde m}\otimes
g^*\mathcal O_{\mathbb P^2}(-\tilde t)$.  Give $K_{\mathbb C}$ the flat
metric defined by $dz$.  Poincar\'e--Lelong and Jensen's formula give
\[
 \tilde t\,T_g(r)
 \leqslant N\!\left(r,\operatorname{div}_{\infty}(s_{\sigma,g})\right)
 +\int_0^{2\pi}\log^+
 \|s_{\sigma,g}(re^{\sqrt{-1}\theta})\|\frac{d\theta}{2\pi}+O(1).
\]
Applying \eqref{eq:canonical-restricted-pole-bound} and
\eqref{eq:canonical-restricted-proximity}, with $p=\tilde m$, yields
\[
 \tilde t\,T_g(r)
\leqslant3\tilde m\sum_{i=1}^3N_g^{[1]}(r,\mathcal C_i)
+O\!\left(\log T_g(r)+\log r\right)\quad\parallel,
\]
which is the asserted estimate.
\end{proof}

\begin{lem}[The zero-evaluation branch]
\label{lem:restricted-zero-evaluation}
In the setting of Theorem~\ref{pro:SMT on Z}, assume that $Z$ is irreducible
and that $\sigma\neq0$.  If the canonical evaluation
$s_{\sigma,g}$ is identically zero, then
\[
 g_{[2]}(\mathbb C)\subset Z(\sigma)\subsetneq Z,
\]
and $g_{[1]}(\mathbb C)$ is contained in a proper algebraic subset of $X_1$.
\end{lem}

\begin{proof}
The tautological derivative $\mathcal Dg_{[1]}$ is not identically zero,
since otherwise $g_{[1]}$, and hence $g$, would be constant.  On the dense
open set where this derivative is nonzero and $g_{[2]}$ is outside
$\Gamma_2$, the canonical inclusion $\iota_p$ is nonvanishing.  The equality
$s_{\sigma,g}=0$ therefore implies
$\sigma(g_{[2]}(z))=0$.  Analytic continuation gives the first inclusion.
Since $\sigma$ is a nonzero section on the irreducible threefold $Z$,
$Z(\sigma)$ is proper and has dimension at most $2$.  Its image under the
proper morphism $\pi_{2,1}$ is consequently a proper algebraic subset of the
threefold $X_1$ containing $g_{[1]}(\mathbb C)$.
\end{proof}

\subsection{An Effective Constant $c$ in the Second Main Theorem for Three Conics}\label{computing c}

Suppose there exists a nonzero negatively twisted invariant logarithmic $2$-jet differential  
\[
\omega_1 \in H^0\bigl(\mathbb{P}^2,\, E_{2,m}T^*_{\mathbb{P}^2}(\log \mathcal{C}) \otimes \mathcal{O}_{\mathbb{P}^2}(-t)\bigr)
\]  
satisfying  
$
0 < t<m
$ 
and such that the divisor $Z$ of $\omega_1$ on the second level $X_2$ of the logarithmic Demailly-Semple tower for $(\mathbb{P}^2, \mathcal{C})$ is irreducible and reduced modulo $\Gamma_2$. Then, by the argument in Section~\ref{subsection:4.1}, for any rational number $\tau \in (0, \tau_1(t/m))$, we obtain a nontrivial global section  
\[
\sigma \in H^0\left( Z,\, \left( \mathcal{O}_{X_2}(2,1) \otimes \pi_{2,0}^* \mathcal{O}_{\mathbb{P}^2}(-\tau) \right)^{\otimes \ell} \right)
\]  
for sufficiently divisible $\ell \gg 1$ (depending on $\tau$).

Furthermore, if $f^* \omega_1 \equiv 0$ and the canonical evaluation of
$\sigma$ along $f_{[2]}$ is nonzero, then
Theorem~\ref{pro:SMT on Z} yields the estimate~\eqref{eq:smt-three-conics}
with constant (see~\eqref{what is tau})
$
    c = \frac{3}{\tau} > \frac{3}{\tau_1(t/m)}.
$  
If the canonical evaluation is zero, Lemma~\ref{lem:restricted-zero-evaluation}
shows that $f_{[1]}$ is algebraically degenerate, and
Theorem~\ref{McQuillan's result} gives the stronger estimate
$3T_f\leqslant\sum_iN_f^{[1]}+o(T_f)$.  Thus, in either branch and by the
arbitrariness of $\tau$, the estimate~\eqref{eq:smt-three-conics} holds for
any constant $c > 3/\tau_1(\frac{t}{m})$.

Direct computation shows that  
$
 \inf_{0<x<1} \frac{3}{\tau_1(x)} 
 =
 \frac{3+\sqrt{6}}{2}.
$
Thus 
$
 \frac{3}{\tau_1(t/m)} < c
$  
is equivalent to $
c > \frac{3+\sqrt{6}}{2}$ and  $0 < \frac{t}{m} < \frac{4c^2-12c+3}{c(4c-3)}.
$

Meanwhile, in {\bf Case 1}, Proposition~\ref{prop:global-slanted-second-jet}
gives the coefficient $m/(t-7)$.  For this coefficient to be at most $c$,
the positive integers $m$ and $t$ associated to the first family $\omega_1$
must satisfy
$
t \geqslant \frac{m}{c} + 7.
$

\smallskip
Taking $c = 19$, we obtain the desired estimate~\eqref{eq:smt-three-conics} for positive integer pairs $(m,t)$ satisfying  
$
0 < t < \frac{1219}{1387}m$ or $t \geqslant \frac{m}{19} + 7.$ 
Note that $8\cdot\frac{1219}{1387}>7$, whereas $7\cdot\frac{1219}{1387}<7$. To exclude the exceptional cases of $(m,t)$, it therefore suffices to prove Key Vanishing Lemma~I. In fact, the exact real boundary is $c>(75+\sqrt{5241})/8$, so $19$ is the smallest integral constant supplied by this argument.

\smallskip
Taking $c = 5$, we obtain the desired estimate~\eqref{eq:smt-three-conics} for positive integer pairs $(m,t)$ satisfying  
$
0 < t < \frac{43}{85}m $ or $ t \geqslant \frac{m}{5} + 7$. 
For $m\geqslant23$ these two regions overlap, since
$43m/85\geqslant m/5+7$. For $3\leqslant m\leqslant22$, direct integer
enumeration shows that, after using monotonicity in the negative twist, the
only minimal uncovered pairs are precisely the eighteen pairs in Key
Vanishing Lemma~II. (For $m=20$ and $m=22$ there is no gap.)

\begin{proof}[Proof of Theorem~\ref{main theorem 3 conics}]
Work on the common nonempty Zariski-open locus on which
Proposition~\ref{prop:step1}, the slanted-vector-field construction, and Key
Vanishing Lemma~II hold.  Let $f$ be algebraically nondegenerate.  The
threshold and irreducible-component argument give a section $\omega_1$ with
$m<5t$ and residual irreducible divisor $Z$ modulo $\Gamma_2$.  If
$f^*\omega_1\not\equiv0$, Theorem~\ref{smt-form-logarithmic-diff-jet} gives
the result because $m/t<5$.

Assume $f^*\omega_1\equiv0$.  A nonconstant lifted curve is not contained in
$\Gamma_2$, so analytic continuation gives $f_{[2]}(\mathbb C)\subset Z$.
If $t\geqslant m/5+7$, Proposition~\ref{prop:global-slanted-second-jet}
either gives the constant-$5$ estimate or makes $f_{[1]}$ algebraically
degenerate, in which case Theorem~\ref{McQuillan's result} is stronger.  If
$t/m<43/85$, then $\tau_1(t/m)>3/5$; choose
$3/5<\tau<\tau_1(t/m)$.  The restricted-section construction and
Theorem~\ref{pro:SMT on Z} give coefficient $3/\tau<5$ when the canonical
evaluation is nonzero, while Lemma~\ref{lem:restricted-zero-evaluation} and
Theorem~\ref{McQuillan's result} handle the zero-evaluation branch.

It remains only to see that one of these alternatives occurs.  For $m=1,2$
this follows from Proposition~\ref{pro:one-jet-vanishing}.  For $m\geqslant23$
the two numerical regions overlap.  For $3\leqslant m\leqslant22$, the
eighteen minimal uncovered pairs are exactly those excluded by Key Vanishing
Lemma~II; twist monotonicity excludes every larger twist with the same $m$.
Finally, algebraic nondegeneracy implies
$O(\log T_f(r)+\log r)=o(T_f(r))$ outside the usual exceptional set.  This
proves~\eqref{eq:smt-three-conics}.
\end{proof}

\section{\bf Proof of the Key Vanishing Lemmas}\label{section 7}

\subsection{Basic Calculations}
\label{subsection 7.1}
Let $\mathcal{C} = \mathcal{C}_1 + \mathcal{C}_2 + \mathcal{C}_3$ be a divisor of three smooth conics in $\mathbb{P}^2$, defined by homogeneous quadratic polynomials:
\begin{align*}
A(Z) &= A_{200}Z_0^2 + A_{020}Z_1^2 + A_{002}Z_2^2 + A_{110}Z_0Z_1 + A_{101}Z_0Z_2 + A_{011}Z_1Z_2, \\
B(Z) &= B_{200}Z_0^2 + B_{020}Z_1^2 + B_{002}Z_2^2 + B_{110}Z_0Z_1 + B_{101}Z_0Z_2 + B_{011}Z_1Z_2, \\
C(Z) &= C_{200}Z_0^2 + C_{020}Z_1^2 + C_{002}Z_2^2 + C_{110}Z_0Z_1 + C_{101}Z_0Z_2 + C_{011}Z_1Z_2,
\end{align*}
where $Z = (Z_0,Z_1,Z_2)$ and all coefficients $A_{\bullet}, B_{\bullet}, C_{\bullet} \in \mathbb{C}$ are generic.

On the affine chart $U_0 = \{Z_0 \neq 0\}$ with coordinates $(x_1,x_2) = (Z_1/Z_0,Z_2/Z_0)$, the dehomogenized polynomials are:
\begin{align*}
a(x_1,x_2) &= A_{200} + A_{020}x_1^2 + A_{002}x_2^2 + A_{110}x_1 + A_{101}x_2 + A_{011}x_1x_2, \\
b(x_1,x_2) &= B_{200} + B_{020}x_1^2 + B_{002}x_2^2 + B_{110}x_1 + B_{101}x_2 + B_{011}x_1x_2, \\
c(x_1,x_2) &= C_{200} + C_{020}x_1^2 + C_{002}x_2^2 + C_{110}x_1 + C_{101}x_2 + C_{011}x_1x_2.
\end{align*}

Their partial derivatives are given by:
\begin{align*}
a_1 &:= \frac{\partial a}{\partial x_1} = 2A_{020}x_1 + A_{110} + A_{011}x_2, \quad
a_2 := \frac{\partial a}{\partial x_2} = 2A_{002}x_2 + A_{101} + A_{011}x_1, \\
b_1 &:= \frac{\partial b}{\partial x_1} = 2B_{020}x_1 + B_{110} + B_{011}x_2, \quad
b_2 := \frac{\partial b}{\partial x_2} = 2B_{002}x_2 + B_{101} + B_{011}x_1, \\
c_1 &:= \frac{\partial c}{\partial x_1} = 2C_{020}x_1 + C_{110} + C_{011}x_2, \quad
c_2 := \frac{\partial c}{\partial x_2} = 2C_{002}x_2 + C_{101} + C_{011}x_1.
\end{align*}

The Jacobian determinant is the homogeneous cubic polynomial:
\[
J(Z) := 
\begin{vmatrix}
\frac{\partial}{\partial Z_0} A & \frac{\partial}{\partial Z_0} B & \frac{\partial}{\partial Z_0} C \\
\frac{\partial}{\partial Z_1} A & \frac{\partial}{\partial Z_1} B & \frac{\partial}{\partial Z_1} C \\
\frac{\partial}{\partial Z_2} A & \frac{\partial}{\partial Z_2} B & \frac{\partial}{\partial Z_2} C
\end{vmatrix}.
\]
Its vanishing locus $\{J=0\}$ is the ramification divisor of the endomorphism:
\begin{equation}
    \label{what is sigma}
\sigma \colon \mathbb{P}^2 \rightarrow \mathbb{P}^2 , \quad [Z_0:Z_1:Z_2] \mapsto [A:B:C].
\end{equation}
The utility of such Jacobian curves has been well-documented in the literature, cf. e.g.~\cite{Francois-Duval-2001,Corvaja-Zannier-2008,Guo-Sun-Wang-2021,Guo-Sun-Wang-2022,CHSX2025}.

Under generic coefficient assumptions, the Jacobian curve $\{J = 0\}$ satisfies the following transversality properties (cf., e.g., \cite{Guo-Sun-Wang-2022, CHSX2025}):
\begin{itemize}
    \item[(\textbf{TP})] Transverse intersection with each  conic $\mathcal{C}_i$ for $i = 1, 2, 3$.
    %\item[(\textbf{TP2})] Transverse intersection with each coordinate hyperplane $\{Z_i = 0\}$ for $i = 0, 1, 2$.
\end{itemize}
For the Key Vanishing Lemma, we will work  with parameters satisfying this transversality conditions.

\medskip

We begin with the following two results, which follow from direct computation on the affine charts $\{Z_i \neq 0\}\subset \mathbb{P}^2$, $i=0,1,2$.

\begin{pro}\label{prop:log-generating}
Let $\mathcal{H} :=  \{Z_0=0\} + \{Z_1=0\} + \{Z_2=0\}$ 
be a divisor consisting of the three coordinate hyperplanes in $\mathbb{P}^2$. Then
the logarithmic cotangent bundle $T^*_{\mathbb{P}^2}(\log \mathcal{H})$ 
is globally generated by the sections
$(\log\frac{Z_0}{Z_2})'$ and $(\log\frac{Z_1}{Z_2})'$. \qed
\end{pro}

\begin{pro}\label{E23 for hyperplanes}
The logarithmic $2$-jet differential bundle $E_{2,m}T^*_{\mathbb{P}^2}(\log \mathcal{H})$ of weighted degree $m$ is globally generated by the following sections:
\begin{equation}\label{basis about Z0/Z2, Z1/Z2}
    \left(\left(\log\frac{Z_0}{Z_2}\right)'\right)^k\left(\left(\log\frac{Z_1}{Z_2}\right)'\right)^{m-3j-k}
    \begin{vmatrix}
        \left(\log\frac{Z_0}{Z_2}\right)' & \left(\log\frac{Z_1}{Z_2}\right)' \\
        \left(\log\frac{Z_0}{Z_2}\right)'' & \left(\log\frac{Z_1}{Z_2}\right)''
    \end{vmatrix}^j
    \quad (0 \leqslant j \leqslant \lfloor m/3 \rfloor , \; 0 \leqslant k \leqslant m-3j).
\end{equation} \qed
\end{pro}

The endomorphism $\sigma$ defined in~\eqref{what is sigma} is in fact a morphism 
$
\sigma : (\mathbb{P}^2, \mathcal{C})  \to
(\mathbb{P}^2, \mathcal{H})
$
between two log pairs. The functoriality  immediately yields the relations between
logarithmic
differential $1$-forms
\[
\Big(\log\frac{A}{C}\Big)' = \sigma^* \Big(\log\frac{Z_0}{Z_2}\Big)'
\quad\text{and}\quad
\Big(\log\frac{B}{C}\Big)' = \sigma^* \Big(\log\frac{Z_1}{Z_2}\Big)'.
\]
Combining this with Proposition~\ref{prop:log-generating}, we obtain:

\begin{pro}\label{prop:log-generating-proper-transform}
On the Zariski-open set $\{J \neq 0\} \subset \mathbb{P}^2$, the logarithmic 1-forms
$(\log\frac{A}{C})'$ and $(\log\frac{B}{C})'$ 
generate the sheaf $T^*_{\mathbb{P}^2}(\log \mathcal{C})$.
\end{pro}

We complement this conceptual proof with explicit computations. While the functorial approach is elegant, the detailed coordinate computations reveal specific algebraic structures that will play an essential role in our later  analysis.

\begin{proof}
The sections $\left(\log\frac{A}{C}\right)'$ and $\left(\log\frac{B}{C}\right)'$ are clearly holomorphic on $\{J \neq 0\}$. We prove generation by working on affine charts.

On $U_0 = \{Z_0 \neq 0\}$ with coordinates $(x_1,x_2)$, consider the transformation matrix between differentials:
\begin{equation*}
    \label{equ:2 times 2 transition formula}
\begin{pmatrix}
\left(\log\frac{a}{c}\right)' \\ \left(\log\frac{b}{c}\right)'
\end{pmatrix} 
= M \begin{pmatrix}
x_1' \\ x_2'
\end{pmatrix}, \quad
M = \begin{pmatrix}
\frac{\partial}{\partial x_1}\log\frac{a}{c} & \frac{\partial}{\partial x_2}\log\frac{a}{c} \\
\frac{\partial}{\partial x_1}\log\frac{b}{c} & \frac{\partial}{\partial x_2}\log\frac{b}{c}
\end{pmatrix}.
\end{equation*}
A direct computation gives the determinant of $M$:
\begin{equation}
\label{what is D}
\det M = \frac{D}{abc}, \quad \text{where} \quad D = \begin{vmatrix}
a & b & c \\ 
a_1 & b_1 & c_1 \\ 
a_2 & b_2 & c_2
\end{vmatrix} = \frac{J}{2Z_0^3}.
\end{equation}

Thus the matrix $M$ is invertible on the Zariski open set $U_0 \cap \{J\neq 0\} \cap \{ABC \neq 0\}$, with its inverse given explicitly by:
\[
M^{-1} = \frac{1}{D}
\begin{pmatrix}
a(c b_2 - b c_2) & b(a c_2 - c a_2) \\
a(b c_1 - c b_1) & b(c a_1 - a c_1)
\end{pmatrix}.
\]

Hence inverting the relation~\eqref{what is D} yields the explicit coordinate transformation:
\begin{equation}
\label{equ:inverse of 2 times 2 transition formula}
\begin{pmatrix}
x_1' \\ x_2'
\end{pmatrix} 
= \frac{1}{D}
\begin{pmatrix}
a(c b_2 - b c_2) & b(a c_2 - c a_2) \\
a(b c_1 - c b_1) & b(c a_1 - a c_1)
\end{pmatrix}
\begin{pmatrix}
\left(\log\frac{a}{c}\right)' \\ \left(\log\frac{b}{c}\right)'
\end{pmatrix}.
\end{equation}
This establishes that the logarithmic differential forms $\left(\log\frac{A}{C}\right)'$ and $\left(\log\frac{B}{C}\right)'$ generate the sheaf of logarithmic cotangent bundle $T^*_{\mathbb{P}^2}(\log\mathcal{C})$ over the Zariski-open set $\{J\neq 0\}\cap \{ABC\neq 0\}$, as the transformation matrix consists of regular functions  on this domain.

Furthermore, the differential $(\log c)'$ can be expressed in terms of our generators:
\begin{align*}
    (\log c)' &= \frac{c'}{c} = \begin{pmatrix}
        \frac{c_1}{c} & \frac{c_2}{c}
    \end{pmatrix}
    \begin{pmatrix}
        x_1' \\ x_2'
    \end{pmatrix} \\
    \text{[use~\eqref{equ:inverse of 2 times 2 transition formula}]}\quad
    &= \frac{1}{D}\begin{pmatrix}
        a \begin{vmatrix}
            c_1 & b_1 \\
            c_2 & b_2
        \end{vmatrix} & 
        b \begin{vmatrix}
            a_1 & c_1 \\
            a_2 & c_2
        \end{vmatrix}
    \end{pmatrix}
    \begin{pmatrix}
        \left(\log\frac{a}{c}\right)' \\ \left(\log\frac{b}{c}\right)'
    \end{pmatrix}.
\end{align*}
Thus the sections $\left(\log\frac{A}{C}\right)'$ and $\left(\log\frac{B}{C}\right)'$ continue to generate the logarithmic cotangent bundle $T^*_{\mathbb{P}^2}(\log\mathcal{C})$ in a neighborhood of $\{J\neq 0\}\cap\{C=0\}$ inside $\{J\neq 0\}$, as the transformation matrix remains regular when approaching this locus.
The analogous argument holds on charts $U_1=\{Z_1 \neq 0\}$ and $U_2=\{Z_2 \neq 0\}$, completing the proof.
\end{proof}

Using Proposition~\ref{E23 for hyperplanes} and considering the morphism $\sigma : (\mathbb{P}^2, \mathcal{C}) \to (\mathbb{P}^2, \mathcal{H})$, by functoriality, we obtain:

\begin{pro}\label{prop4}
The bundle $E_{2,m}T^*_{\mathbb{P}^2}(\log\mathcal{C})$ is freely generated over $\{J \neq 0\}$ by the logarithmic jet differentials:
\begin{equation}\label{basis about A/C, B/C}
    \left(\left(\log\frac{A}{C}\right)'\right)^k\left(\left(\log\frac{B}{C}\right)'\right)^{m-3j-k}
    \begin{vmatrix}
        \left(\log\frac{A}{C}\right)' & \left(\log\frac{B}{C}\right)' \\
        \left(\log\frac{A}{C}\right)'' & \left(\log\frac{B}{C}\right)''
    \end{vmatrix}^j
    \quad (0 \leqslant j \leqslant \lfloor m/3 \rfloor , \; 0 \leqslant k \leqslant m-3j).
    \qed
\end{equation}
\end{pro}

Let $P \in \Gamma \big( \{J \neq 0\} , E_{2,m}T^*{\mathbb{P}^2}(\log\mathcal{C}) \otimes \mathcal{O}(-t) \big)$. To check and verify the global extendability of $P$, it suffices to demonstrate that $P$ admits a holomorphic extension across the Jacobian divisor $\{J = 0\}$. Note that  the intersection $\{J = 0\} \cap \{Z_0ABC = 0\}$ constitutes a Zariski-closed subset of codimension $2$ in $\mathbb{P}^2$. By Hartog's extension theorem, it is thus sufficient to extend $P$ to $\{J = 0\} \cap \{Z_0ABC \neq 0\}$. This extension problem necessitates an analysis of the transition formula for $E_{2,m}T^*_{\mathbb{P}^2}(\log\mathcal{C}) \otimes \mathcal{O}(-t)$ from $\{J \neq 0\}$ to $\{Z_0ABC \neq 0\}$.

Our carefully chosen basis \eqref{basis about A/C, B/C} induces a highly
homogeneous representation of the logarithmic jet differentials and keeps
the affine transition problem finite.

\medskip
The remaining argument treats all Wronskian layers simultaneously.  This is
essential for the strengthened list, whose largest target has $m=21$; a
small-$m$ pole bound does not determine the higher layers that occur there.

\subsection{All-layer finite shape and exact C++ verification}
\label{subsec:all-layer-cpp}

In what follows, we will abbreviate the first-order jet transfer formulas from standard jets $(x_1', x_2')$ to logarithmic jets $((\log\frac{a}{c})', (\log\frac{b}{c})')$ as:
\[
\begin{aligned}
\Big(\log\frac{a}{c}\Big)' &= \alpha\,x_1' + \beta\,x_2', \\
\Big(\log\frac{b}{c}\Big)' &= \gamma\,x_1' + \delta\,x_2',
\end{aligned}
\qquad\qquad
\begin{aligned}
\alpha &:= \frac{a_{1}c - a c_{1}}{ac}, &
\beta &:= \frac{a_{2}c - a c_{2}}{ac}, \\
\gamma &:= \frac{b_{1}c - b c_{1}}{bc}, &
\delta &:= \frac{b_{2}c - b c_{2}}{bc}.
\end{aligned}
\]

Because logarithmic jet differentials are by definition assumed to exist and to be holomorphic outside the divisor $\mathcal{C} = \{ABC = 0\}$, the presence of $a$, $b$, $c$ (within the affine chart $\C^2 \subset \P^2$) in denominators does not produce any singularity.

The only possible singularity---which clarifies why Proposition \ref{prop4} holds only on the complement $\P^2 \setminus \{J = 0\}$---comes from inverting the transfer formulas from logarithmic jets back to standard jets:
\begin{equation}\label{transfer formulas from logarithmic jets to standard jets}
\begin{aligned}
 x_1'&=
\frac{\delta}{\alpha\delta - \beta\gamma}\,\Big(\log\frac{a}{c}\Big)' 
- \frac{\beta}{\alpha\delta - \beta\gamma}\,\Big(\log\frac{b}{c}\Big)' , \\
x_2'&= 
-\frac{\gamma}{\alpha\delta - \beta\gamma}\,\Big(\log\frac{a}{c}\Big)' 
+ \frac{\alpha}{\alpha\delta - \beta\gamma}\,\Big(\log\frac{b}{c}\Big)' ,
\end{aligned}
\end{equation}
and this requires the nonvanishing of the $2\times 2$ determinant (see \eqref{what is D}):
\begin{equation}
    \label{key point about D}
\alpha\delta - \beta\gamma=
\det\begin{pmatrix}
\alpha & \beta \\
\gamma & \delta
\end{pmatrix}
= \det M
= \frac{D}{abc},
\quad \text{where} \quad 
D = \begin{vmatrix}
a & b & c \\ 
a_1 & b_1 & c_1 \\ 
a_2 & b_2 & c_2
\end{vmatrix} = \frac{J}{2Z_0^3}.
\end{equation}

By a direct computation in the affine chart $\{Z_0 \neq 0\}$, this determinant happens to be a nonzero multiple of the Jacobian determinant $J$.

The transfer of the logarithmic Wronskian, 
\[
\overline{W}^{\,\prime\prime\prime}
\coloneqq
\begin{vmatrix}
\big(\log\frac{A}{C}\big)^\prime & 
\big(\log\frac{B}{C}\big)^\prime
\\
\big(\log\frac{A}{C}\big)^{\prime\prime} & 
\big(\log\frac{B}{C}\big)^{\prime\prime}
\end{vmatrix}
=
\begin{vmatrix}
\big(\log\frac{a}{c}\big)^\prime & 
\big(\log\frac{b}{c}\big)^\prime
\\
\big(\log\frac{a}{c}\big)^{\prime\prime} & 
\big(\log\frac{b}{c}\big)^{\prime\prime}
\end{vmatrix},
\]
to standard affine jet coordinates is more delicate, and is captured in the following Key Observation.

\begin{KeyObservation}\label{expession of overline W'''}
When expressed in standard affine jet coordinates $(x_1', x_2', x_1'', x_2'')$, the logarithmic Wronskian takes the form:
\begin{equation}\label{transfer formula of log Wronskian}
\begin{aligned}
\overline{W}^{\,\prime\prime\prime}
&=
(\alpha\delta - \beta\gamma)\,
\begin{vmatrix}
x_1' & x_2' \\
x_1'' & x_2''
\end{vmatrix}
+ \kappa\,{x_1'}^3 + \lambda\,{x_1'}^2x_2' + \mu\,x_1'{x_2'}^2 + \nu\,{x_2'}^3,
\end{aligned}
\end{equation}
where the cubic terms in the first-order jets $x_1'$, $x_2'$ have coefficients $\kappa$, $\lambda$, $\mu$, $\nu$ given in terms of the second-order derivatives
\[
a_{11},\ a_{12},\ a_{22},\qquad 
b_{11},\ b_{12},\ b_{22},\qquad 
c_{11},\ c_{12},\ c_{22},
\]
and the coefficients $\alpha$, $\beta$, $\gamma$, $\delta$ as follows:
\[
\begin{aligned}
\kappa &:=
-\gamma\frac{a_{11}}{a}
+ \alpha\frac{b_{11}}{b}
- (\alpha-\gamma)\frac{c_{11}}{c}
+ \alpha^2\gamma - \alpha\gamma^2,
\\[2mm]
\lambda &:=
-\delta\frac{a_{11}}{a}
- 2\gamma\frac{a_{12}}{a}
+ \beta\frac{b_{11}}{b}
+ 2\alpha\frac{b_{12}}{b}
- (\beta-\delta)\frac{c_{11}}{c}
- 2(\alpha-\gamma)\frac{c_{12}}{c}
\\
&\qquad\qquad
+ \alpha^2\delta + 2\alpha\beta\gamma - 2\alpha\gamma\delta - \beta\gamma^2,
\\[2mm]
\mu &:=
-2\delta\frac{a_{12}}{a}
- \gamma\frac{a_{22}}{a}
+ 2\beta\frac{b_{12}}{b}
+ \alpha\frac{b_{22}}{b}
- 2(\beta-\delta)\frac{c_{12}}{c}
- (\alpha-\gamma)\frac{c_{22}}{c}
\\
&\qquad\qquad
+ 2\alpha\beta\delta - \alpha\delta^2 + \beta^2\gamma - 2\beta\gamma\delta,
\\[2mm]
\nu &:=
-\delta\frac{a_{22}}{a}
+ \beta\frac{b_{22}}{b}
- (\beta-\delta)\frac{c_{22}}{c}
\\
&\qquad\qquad
+ \beta^2\delta - \beta\delta^2.
\end{aligned}
\]
\end{KeyObservation}

\proof
This follows by direct computation.
\endproof

According to Key Vanishing Lemma~II, it is enough to consider
$m\leqslant21$.  On the affine chart $(x_1,x_2)\in\C^2\subset\P^2$, away
from the Jacobian zero-locus
\[
\{D=0\} = \{J=0\} \cap \{Z_0 \neq 0\},
\]
the basis in Proposition~\ref{prop4} gives the uniform expansion
\begin{equation}
\label{eq:Jet-ABC-expansion}
\Jet_{ABC}
=\sum_{q=0}^{\floor{m/3}}\ \sum_{i+j=m-3q}
 P_{i,j}^{(q)}
 \Big(\big(\log\frac{a}{c}\big)'\Big)^i
 \Big(\big(\log\frac{b}{c}\big)'\Big)^j
 \big(\overline W^{\,\prime\prime\prime}\big)^q.
\end{equation}
Every coefficient $P_{i,j}^{(q)}$ is holomorphic away from $\{D=0\}$.

The {\sc gaga} principle then ensures that every such coefficient is a
rational expression whose denominator is an integral power of $D$.  Our
preliminary objective is to determine a uniform bound for that power in each
Wronskian layer.

To this end, we consider the transfer formulas from logarithmic jets to standard jets. The first-order jet transfer formulas are already given in \eqref{transfer formulas from logarithmic jets to standard jets}. Using \eqref{transfer formula of log Wronskian}, we have
\begin{equation*}
W''' \coloneqq 
\begin{vmatrix}
x_1' & x_2' \\
x_1'' & x_2''
\end{vmatrix}
= \frac{1}{\alpha\delta - \beta\gamma}\,
\overline{W}^{\,\prime\prime\prime}
- \frac{1}{\alpha\delta - \beta\gamma}
\Bigl( 
\kappa\,{x_1'}^3
+ \lambda\,{x_1'}^2x_2'
+ \mu\,x_1'{x_2'}^2
+ \nu\,{x_2'}^3
\Bigr),
\end{equation*}
where the cubic term $\kappa\, {x_1'}^3 + \lambda\, {x_1'}^2x_2' + \mu\, x_1'{x_2'}^2 + \nu\, {x_2'}^3$ can be expanded using \eqref{transfer formulas from logarithmic jets to standard jets} as:
\begingroup
\setlength{\arraycolsep}{3pt}
\begin{align*}
& \begin{pmatrix}
    {x_1'}^3 & {x_1'}^2x_2' & x_1'{x_2'}^2 & {x_2'}^3
\end{pmatrix}
\begin{pmatrix}
    \kappa \\
    \lambda \\
    \mu \\
    \nu
\end{pmatrix} \\
=& \begin{pmatrix}
    \big(\big(\log\frac{a}{c}\big)^{\prime}\big)^3 
    &
    \big(\big(\log\frac{a}{c}\big)^{\prime}\big)^2 
    \big(\log\frac{b}{c}\big)^{\prime}
    &
    \big(\log\frac{a}{c}\big)^{\prime} 
    \big(\big(\log\frac{b}{c}\big)^{\prime}\big)^2
    &
    \big(\big(\log\frac{b}{c}\big)^{\prime}\big)^3
\end{pmatrix} \\
& \quad \cdot
\frac{1}{(\alpha\delta - \beta\gamma)^3}
\begin{pmatrix}
\delta^3 & -\gamma\delta^2 & \gamma^2\delta & -\gamma^3 \\
-3\beta\delta^2 & \alpha\delta^2+2\beta\gamma\delta & -2\alpha\gamma\delta-\beta\gamma^2 & 3\alpha\gamma^2 \\
3\beta^2\delta & -2\alpha\beta\delta-\beta^2\gamma & \alpha^2\delta+2\alpha\beta\gamma & -3\alpha^2\gamma \\
-\beta^3 & \alpha\beta^2 & -\alpha^2\beta & \alpha^3
\end{pmatrix}
\begin{pmatrix}
    \kappa \\
    \lambda \\
    \mu \\
    \nu
\end{pmatrix}.
\end{align*}
\endgroup
Let us examine this more closely. Define
\[
\begin{pmatrix}
\delta^3 & -\gamma\delta^2 & \gamma^2\delta & -\gamma^3 \\
-3\beta\delta^2 & \alpha\delta^2+2\beta\gamma\delta & -2\alpha\gamma\delta-\beta\gamma^2 & 3\alpha\gamma^2 \\
3\beta^2\delta & -2\alpha\beta\delta-\beta^2\gamma & \alpha^2\delta+2\alpha\beta\gamma & -3\alpha^2\gamma \\
-\beta^3 & \alpha\beta^2 & -\alpha^2\beta & \alpha^3
\end{pmatrix}
\begin{pmatrix}
\kappa \\
\lambda \\
\mu \\
\nu
\end{pmatrix}
\eqqcolon
\begin{pmatrix}
\tau_{3,0} \\
\tau_{2,1} \\
\tau_{1,2} \\
\tau_{0,3}
\end{pmatrix}.
\]

\begin{MiracleFactorizations}
\label{Miracle-k-l-m-n}
The four expressions:
\[
\aligned
\tau_{3,0}
\,:=\,
&\,
\delta^3\,
\kappa
-
\gamma\,\delta^2\,
\lambda
+
\gamma^2\,\delta\,
\mu
-
\gamma^3\,
\nu
\\
\,=\,
&\,
\big(
\alpha\,\delta-\beta\,\gamma
\big)\,
\big(\cdots\big),
\\
\tau_{2,1}
\,:=\,
&\,
-\,
3\,\beta\,\delta^2\,
\kappa
+
\big(\alpha\,\delta^2+2\,\beta\,\gamma\,\delta\big)\,
\lambda
+
\big(-\,2\,\alpha\,\gamma\,\delta-\beta\,\gamma^2\big)\,
\mu
+
3\,\alpha\,\gamma^2\,
\nu
\\
\,=\,
&\,
\big(
\alpha\,\delta-\beta\,\gamma
\big)\,
\big(\cdots\big),
\\
\tau_{1,2}
\,:=\,
&\,
3\,\beta^2\,\delta\,
\kappa
+
\big(-\,2\,\alpha\,\beta\,\delta-\beta^2\,\gamma\big)\,
\lambda
+
\big(\alpha^2\,\delta+2\,\alpha\,\beta\,\gamma\big)\,
\mu
-
3\,\alpha^2\,\gamma\,\nu
\\
\,=\,
&\,
\big(
\alpha\,\delta-\beta\,\gamma
\big)\,
\big(\cdots\big),
\\
\tau_{0,3}
\,:=\,
&\,
-\,\beta^3\,
\kappa
+
\alpha\,\beta^2\,\lambda
-
\alpha^2\,\beta\,\mu
+
\alpha^3\,\nu
\\
\,=\,
&\,
\big(
\alpha\,\delta-\beta\,\gamma
\big)\,
\big(\cdots\big),
\endaligned
\]
are all multiples of $\alpha\delta - \beta\gamma$.
\end{MiracleFactorizations}

\proof
Substitute the displayed formulas for $\kappa,\lambda,\mu,\nu$ and divide
in the polynomial ring generated by $\alpha,\beta,\gamma,\delta$ and the
normalized second derivatives.  Each division has zero remainder.  The four
quotients are displayed and independently checked by exact symbolic
expansion in Appendix~\ref{appendix:miracle-factorizations}.
\endproof

Therefore:
\[
\begin{pmatrix}
\tau_{3,0} \\
\tau_{2,1} \\
\tau_{1,2} \\
\tau_{0,3}
\end{pmatrix}
\eqqcolon
(\alpha\delta - \beta\gamma)
\begin{pmatrix}
\theta_{3,0} \\
\theta_{2,1} \\
\theta_{1,2} \\
\theta_{0,3}
\end{pmatrix},
\]
and finally:
\begin{equation}\label{transfer formula of standard Wronskian}
\begin{vmatrix}
x_1' & x_2' \\
x_1'' & x_2''
\end{vmatrix}
= \frac{1}{\alpha\delta - \beta\gamma}
\begin{vmatrix}
\big(\log\frac{a}{c}\big)^\prime & 
\big(\log\frac{b}{c}\big)^\prime \\
\big(\log\frac{a}{c}\big)^{\prime\prime} & 
\big(\log\frac{b}{c}\big)^{\prime\prime}
\end{vmatrix}
- \frac{1}{(\alpha\delta - \beta\gamma)^3}
\sum_{i+j=3} \theta_{i,j}\, 
\big(\big(\log\frac{a}{c}\big)^{\prime}\big)^i 
\big(\big(\log\frac{b}{c}\big)^{\prime}\big)^j.
\end{equation}

Recall that our main assumption is that, in the standard affine coordinates $(x_1,x_2) \in \C^2$, the logarithmic $2$-jet differential $\Jet_{ABC}$ is holomorphic over $\C^2 \setminus \{abc=0\}$. Therefore, it can be expressed as 
\begin{equation}
\label{eq:Jet-ABC-standard-expansion}
\Jet_{ABC}
=\sum_{q=0}^{\floor{m/3}}\ \sum_{r+s=m-3q}
 \overline P_{r,s}^{(q)}\,{x_1'}^r{x_2'}^s\big(W'''\big)^q,
\end{equation}
where every $\overline P_{r,s}^{(q)}$ is holomorphic on
$\C^2\setminus\{abc=0\}$.  By the {\sc gaga} principle it is a rational
function whose denominator may contain only $a$, $b$, and $c$.

We now compare \eqref{eq:Jet-ABC-expansion} with
\eqref{eq:Jet-ABC-standard-expansion}.  Substitution of
\eqref{transfer formulas from logarithmic jets to standard jets} and
\eqref{transfer formula of standard Wronskian} shows that the coefficients
$\{P_{i,j}^{(q)}\}_{i+j=m-3q}$ can be written as rational expressions whose
denominators consist of factors $a$, $b$, $c$ (which are coprime to $D$) and
$D$ (see~\eqref{key point about D}).  In Wronskian layer $q$, the maximal
possible pole order with respect to $D$ is $m-2q$.

Indeed, an ordinary first derivative contributes at most one inverse power of
$D$, while each occurrence of the lower cubic term in
\eqref{transfer formula of standard Wronskian} trades three first-derivative
weights for two additional inverse powers of $D$.  For a monomial of total
weight $m$ in layer $q$, this gives precisely $D^{-(m-2q)}$.  The argument is
independent of $q$ and therefore covers every $0\leqslant q\leqslant
\floor{m/3}$.

Consequently the invariant logarithmic $2$-jet differential $\Jet_{ABC}$ in
\eqref{eq:Jet-ABC-expansion}, of weighted order $m\leqslant21$, has the
following uniform all-layer form in homogeneous coordinates:
\begin{equation}
\label{final formula J}
\Jet_{ABC}
=\sum_{q=0}^{\floor{m/3}}\ \sum_{i+j=m-3q}
 \frac{P_{i,j}^{(q)}(Z_0,Z_1,Z_2)}
 {J(Z_0,Z_1,Z_2)^{m-2q}}
 \Big(\big(\log\frac{A}{C}\big)'\Big)^i
 \Big(\big(\log\frac{B}{C}\big)'\Big)^j
 \big(\overline W^{\,\prime\prime\prime}\big)^q,
\end{equation}
where
\begin{equation}\label{eq:all-layer-degrees}
n_q:=m-2q,
\qquad
P_{i,j}^{(q)}\in\C[Z_0,Z_1,Z_2]_{d_q},
\qquad
d_q:=3n_q-t.
\end{equation}
Indeed, homogeneity gives $d_q-3(m-2q)=-t$ in every layer.  Thus
\eqref{final formula J} is not a truncation at a fixed Wronskian power: it is
the complete finite shape, and it is precisely the shape enumerated by the
C++ verifier.  It represents a \emph{meromorphic} (rational) section of
\[
E_{2,m}T_{\P^2}^\ast\big(\log\mathcal{C}\big)
\otimes
\mathcal{O}_{\P^2}(-t).
\]

For $\Jet_{ABC}$ to be a genuine \emph{holomorphic} logarithmic jet
differential, the coefficients of all polynomials $P_{i,j}^{(q)}$ must
satisfy a finite homogeneous linear
system.  The system records that, after the complete coupled transition to
ordinary affine jets, every coefficient has the divisibility required to
cancel the appropriate power of the Jacobian $J$.  Cancellations are taken
only after summing the contributions of all Wronskian layers to the same
ordinary-jet monomial; testing the layers separately would give only a
diagnostic, not the desired kernel computation.

The number of homogeneous coefficient variables is
\begin{equation}\label{eq:KVL-II-variable-count}
 N(m,t)=
 \sum_{q=0}^{\floor{m/3}}
 (m-3q+1)\binom{3(m-2q)-t+2}{2},
\end{equation}
where a summand is omitted when $3(m-2q)-t<0$.  In particular,
$N(13,9)=12{,}550$, whereas the largest strengthened target
$(m,t)=(21,11)$ has $N(21,11)=81{,}840$ variables.

To streamline explicit computer-assisted computations, we adopt a strategic simplification by choosing specific coefficients $A_{\bullet}, B_{\bullet}, C_{\bullet}$ such that the three conics $\mathcal{C}$ are of Fermat type; for example:
\begin{align}
\label{Fermat conics}
    A(Z) = 2 Z_0^2 + Z_1^2 + Z_2^2, \quad
    B(Z) = Z_0^2 + 2Z_1^2 + Z_2^2, \quad
    C(Z) = Z_0^2 + Z_1^2 + 2Z_2^2.
\end{align}
Each of the three conics in \eqref{Fermat conics} is smooth.  Any two meet
transversally, and the three have no common point: subtracting their equations
gives $Z_0^2=Z_1^2=Z_2^2$, after which any one equation forces all three
coordinates to vanish.  Hence their union is a simple normal crossing
divisor.  This choice also gives
$J=32Z_0Z_1Z_2$, which makes the divisibility test particularly transparent.

Then we test extension on both
$\{Z_0ABC\neq0\}$ and $\{Z_2ABC\neq0\}$.  On
$\P^2\setminus\mathcal C$, their union contains $\{J=0\}$ except for the
single point $[0:1:0]$.  That omitted set has codimension two, so Hartogs'
theorem completes the extension there.
% In $\{Z_0 ABC \neq 0\}$, all the numerators in the expansion of $\Jet_{ABC}$ to the affine jet coordinates $(x_1, x_2, x_1', x_2', x_1'', x_2'')$ must be divisible by $(32 x_1 x_2)^m$; meanwhile, in $\{Z_1 ABC \neq 0\}$, with coordinates $(y_0,y_2) = (Z_0/Z_1,Z_2/Z_1)$, all the numerators in the expansion of $\Jet_{ABC}$ to the affine jet coordinates $(y_0, y_2, y_0', y_2', y_0'', y_2'')$ must be divisible by $(32 y_0 y_2)^m$.

More precisely, we work in the two affine charts:
\[
\{Z_0 \neq 0\} \qquad \text{and} \qquad \{Z_2 \neq 0\},
\]
with the coordinate substitutions:
\[
(Z_0,Z_1,Z_2) \longmapsto (1,x_1,x_2)
\qquad \text{and} \qquad
(Z_0,Z_1,Z_2) \longmapsto (y_0,y_1,1).
\]
In these coordinates, $J$ reads as $32x_1x_2$ or $32y_0y_1$.
Thus the required cancellation is equivalent to divisibility by
$(x_1x_2)^{m-2k}$ or $(y_0y_1)^{m-2k}$ in the coefficient of the ordinary
Wronskian power $W^k$.  These conditions form the finite homogeneous linear
system computed below.

\medskip\noindent\textbf{The exact matrix computed by the C++ verifier.}
For $\nu\in\{0,2\}$, let $U_\nu=\{Z_\nu\neq0\}$ and write
\[
 (u_1,u_2)=
 \begin{cases}
   (x_1,x_2),&\nu=0,\\
   (y_0,y_1),&\nu=2.
 \end{cases}
\]
Let $a,b,c$ denote the dehomogenizations of the homogeneous quadrics
$A,B,C$ on the current chart.  We suppress the chart index below, and
likewise use the same symbols $J$ and $P_{i,j}^{(q)}$ for their
dehomogenizations.  Put
\[
 H:=abc,
 \qquad
 \widehat L_1:=b(ca'-ac'),
 \qquad
 \widehat L_2:=a(cb'-bc'),
\]
and define
\[
 D:=
 \det
 \begin{pmatrix}
   a&b&c\\
   \partial_{u_1}a&\partial_{u_1}b&\partial_{u_1}c\\
   \partial_{u_2}a&\partial_{u_2}b&\partial_{u_2}c
 \end{pmatrix},
 \qquad
 W:=u_1'u_2''-u_1''u_2'.
\]
Then
\[
 H\bigl(\log(a/c)\bigr)'=\widehat L_1,
 \qquad
 H\bigl(\log(b/c)\bigr)'=\widehat L_2.
\]
Define the cubic first-jet polynomial $R_3$ by
\begin{equation}\label{eq:cpp-R3-definition}
 \widehat L_1(\widehat L_2)'
 -\widehat L_2(\widehat L_1)'
 =HDW+R_3.
\end{equation}
Equivalently,
\begin{equation}\label{eq:cpp-cleared-log-Wronskian}
 H^3\overline W^{\,\prime\prime\prime}
 =H^2DW+HR_3.
\end{equation}
The program constructs $R_3$ from \eqref{eq:cpp-R3-definition}; it does not
enter a precomputed formula for it.

Fix a source layer $q$ in \eqref{final formula J}, and put
\[
 n_q=m-2q,
 \qquad s_q=m-3q,
 \qquad i+j=s_q,
 \qquad \deg P_{i,j}^{(q)}=3n_q-t.
\]
After substituting \eqref{eq:cpp-cleared-log-Wronskian}, the contribution of
$P_{i,j}^{(q)}$ to the coefficient of $W^k$, $0\leqslant k\leqslant q$,
is placed over the common denominator
$H^{m-2k}J^{m-2k}$.  Let $T_{i,j,q}$ denote the $(i,j,q)$-summand in
\eqref{final formula J}.  Using
\[
 \left(\log\frac AC\right)'=\frac{\widehat L_1}{H},
 \qquad
 \left(\log\frac BC\right)'=\frac{\widehat L_2}{H},
 \qquad
 \overline W^{\,\prime\prime\prime}
 =\frac{HDW+R_3}{H^2},
\]
and $i+j=m-3q$, we obtain, for $0\leqslant k\leqslant q$,
\[
 [W^k]T_{i,j,q}
 =
 \frac{
 \binom qk P_{i,j}^{(q)}
 \widehat L_1^i\widehat L_2^jD^k
 \big(HJ^2R_3\big)^{q-k}
 }{
 H^{m-2k}J^{m-2k}
 }.
\]
Thus the polynomial in the numerator
\begin{equation}\label{eq:cpp-column-contribution} 
\binom qk P_{i,j}^{(q)}\widehat L_1^i\widehat L_2^jD^k \big(HJ^2R_3\big)^{q-k} 
\end{equation}
is the contribution of the source
term $T_{i,j,q}$ to the common numerator of the $W^k$-coefficient.
For fixed $k$ and a fixed velocity monomial, the program first sums
\eqref{eq:cpp-column-contribution} over \emph{all} $q\geqslant k$ and only
then tests divisibility by $J^{m-2k}$.  This order matters: cancellation
between different Wronskian layers is part of the extension problem.  A
layer-by-layer rank test would therefore compute a different, generally
stronger system.

Each scalar coefficient of each homogeneous polynomial $P_{i,j}^{(q)}$ is a
column.  A row is indexed by a chart, the power $k$, a velocity monomial, and
a base monomial that is not divisible by the required power of the monomial
Jacobian.  The verifier generates these sparse columns on demand and performs
exact column elimination over $\mathbb F_5$; it never materializes a dense
matrix. 

Before a rank computation, the self-test verifies the determinant identity
for $D$, checks \eqref{eq:cpp-cleared-log-Wronskian} directly on nonsingular
jets in both charts, and checks the expected variable counts.  The transition
matrix has integer entries.  Hence full column rank after reduction modulo
$5$ exhibits an integer maximal minor that is nonzero modulo $5$; that minor
is nonzero over $\mathbb Q$, and the original matrix has zero kernel over
$\mathbb C$.

All eighteen matrices in Key Vanishing Lemma~II have full column rank.  For
the largest case $(m,t)=(21,11)$, the program processes $81{,}840$ columns,
finds rank $81{,}840$, and returns nullity zero.\footnote{The complete C++
evidence archive---source, mathematical explanation, certificates, hashes,
and full logs---is available from \url{https://xiesongyan.github.io/}.}

\begin{rmk}
\label{m=5 for three Fermat-type conics}
For the explicit Fermat-type configuration~\eqref{Fermat conics}, our exact
finite computations give the vanishings for
$(m,t)=(1,1),(2,1),(3,1),(4,1),(5,2)$.  For $(m,t)=(5,1)$, by contrast,
they give
\[
\dim\, H^0\big(\mathbb{P}^2 , E_{2,5}T^*_{\mathbb{P}^2}(\log \mathcal{C}) \otimes \mathcal{O}_{\mathbb{P}^2}(-1)\big) = 9.
\]
 
Let $\eta_1,\ldots,\eta_9$ be a basis of this space.  The vanishings in
weighted degrees $1,2,3,4$ exclude a common invariant jet generator of
smaller positive weight.  The standard minimal-weight argument therefore
shows that the $9$-dimensional space contains two algebraically independent
members, say $\eta$ and $\theta$.  If either pullback is nonzero,
Theorem~\ref{smt-form-logarithmic-diff-jet} gives the desired estimate with
coefficient $m/t=5$.  If both pullbacks vanish, the common zero locus of
$\eta$ and $\theta$ has proper image in the first Semple stage; hence
$f_{[1]}$ is algebraically degenerate and
Theorem~\ref{McQuillan's result} applies.  Thus the nine sections give a
Second Main Theorem with constant $5$ for the Fermat-type
configuration~\eqref{Fermat conics}.  This special conclusion is not needed
in the proof of Theorem~\ref{main theorem 3 conics}.
\end{rmk}

\section*{\bf Declaration of AI Use}

AI tools were used in the preparation of this article as follows.
The proof of Proposition~\ref{pro:one-jet-vanishing} was provided
by ChatGPT (model: GPT-5.6 Sol). The authors have carefully
checked the proof and verified its mathematical correctness.

The C++ code used for the computations in this article was
produced by Codex by translating and adapting Maple code
written by the authors. The authors have checked the resulting
C++ code and verified the computational results used in the article.

AI tools were also used to assist with writing and to improve
the clarity and language of this article. The authors have
reviewed the resulting text and take full responsibility for
the mathematical content, the code and computational results,
and the final wording of this article.

\section*{\bf Acknowledgements}

We are grateful to Julien Duval for introducing us a decade ago to the problem of establishing a Second Main Theorem for three conics in $\mathbb{P}^2$. We thank Junjiro Noguchi for sharing his private notes and providing helpful comments. We thank Duc Viet Vu for offering valuable insights from a pluripotential perspective.

S.-Y. Xie would like to thank Fr\'ed\'eric Campana, Yunling Chen, Tao Cui, Hanlong Fang, Jie Liu, Mihai P\u{a}un, Erwan Rousseau, Yum-Tong Siu, Ruiran Sun, and Kang Zuo for inspiring discussions.
S.-Y. Xie thanks Daniel Barlet and Renjie Lv for discussions concerning the one-jet vanishing problem.
S.-Y. Xie and L. Hou also wish to thank  Tianyuan Mathematics Research Center for its excellent  environment. 

We also thank Yunling Chen, Alexandre Eremenko, and Yi C. Huang for their  suggestions on an earlier version.
We thank Steven Lu and Katsutoshi Yamanoi for helpful comments. 

We are grateful to Eric Riedl for pointing out an error in
Proposition~\ref{pro:the higher direct images vanishing} in an earlier
version and for providing a concrete counterexample. This observation helped
us correct the relevant argument.

We thank Pengchao Wang for pointing out two important errors in an earlier version.

The results of this article were first presented by S.-Y. Xie at {\sl The Seminar on Several Complex Variables and Complex Geometry} (held from July 27 to August 2, 2025), which was organized by Xiangyu Zhou, Yum-Tong Siu, and John Erik Fornæss.

\section*{\bf Funding}

S.-Y. Xie acknowledges partial support from the National Key R\&D Program of China under Grants No. 2023YFA1010500 and No. 2021YFA1003100, and from the National Natural Science Foundation of China under Grants No. 12288201 and No. 12471081, as well as support from the 
Xiaomi Young Talents Program. D.T. Huynh is supported by the Vietnam National Foundation for Science and Technology Development (NAFOSTED) under the grant number 101.04-2025.19.

\medskip
\appendix

\section{Explicit quotients in the four factorization identities}
\label{appendix:miracle-factorizations}

We complete the calculation in
Miracle Factorizations~\ref{Miracle-k-l-m-n}.  Set
\[
 \widehat a_{ij}:=\frac{a_{ij}}a,
 \qquad
 \widehat b_{ij}:=\frac{b_{ij}}b,
 \qquad
 \widehat c_{ij}:=\frac{c_{ij}}c,
 \qquad
 \Delta:=\alpha\delta-\beta\gamma.
\]
After substituting the formulas for
$\kappa,\lambda,\mu,\nu$, direct expansion gives
\[
 \tau_{3,0}=\Delta\theta_{3,0},\qquad
 \tau_{2,1}=\Delta\theta_{2,1},\qquad
 \tau_{1,2}=\Delta\theta_{1,2},\qquad
 \tau_{0,3}=\Delta\theta_{0,3},
\]
where
\begin{align*}
\theta_{3,0}={}&
 \delta^2(\widehat b_{11}-\widehat c_{11})
 -2\delta\gamma(\widehat b_{12}-\widehat c_{12})
 +\gamma^2(\widehat b_{22}-\widehat c_{22}),
\\[2mm]
\theta_{2,1}={}&
 -\delta^2\widehat a_{11}
 +2\delta\gamma\widehat a_{12}
 -\gamma^2\widehat a_{22}
 -2\beta\delta\widehat b_{11}
 +2(\alpha\delta+\beta\gamma)\widehat b_{12}
 -2\alpha\gamma\widehat b_{22}
\\
& +(2\beta\delta+\delta^2)\widehat c_{11}
 -2(\alpha\delta+\beta\gamma+\delta\gamma)\widehat c_{12}
 +(2\alpha\gamma+\gamma^2)\widehat c_{22}
 +\Delta^2,
\\[2mm]
\theta_{1,2}={}&
 2\beta\delta\widehat a_{11}
 -2(\alpha\delta+\beta\gamma)\widehat a_{12}
 +2\alpha\gamma\widehat a_{22}
 +\beta^2\widehat b_{11}
 -2\alpha\beta\widehat b_{12}
 +\alpha^2\widehat b_{22}
\\
& -(\beta^2+2\beta\delta)\widehat c_{11}
 +2(\alpha\beta+\alpha\delta+\beta\gamma)\widehat c_{12}
 -(\alpha^2+2\alpha\gamma)\widehat c_{22}
 -\Delta^2,
\\[2mm]
\theta_{0,3}={}&
 -\beta^2(\widehat a_{11}-\widehat c_{11})
 +2\alpha\beta(\widehat a_{12}-\widehat c_{12})
 -\alpha^2(\widehat a_{22}-\widehat c_{22}).
\end{align*}

For completeness, we also verified these identities by exact symbolic
computation.  The calculation was performed with SymPy~1.14.0 in the
polynomial ring
\[
\mathbb Z[\alpha,\beta,\gamma,\delta,
 \widehat a_{11},\widehat a_{12},\widehat a_{22},
 \widehat b_{11},\widehat b_{12},\widehat b_{22},
 \widehat c_{11},\widehat c_{12},\widehat c_{22}].
\]
For each of the four expressions, exact division by $\Delta$ returned the
quotient displayed above and remainder zero.  A second expansion of
$\tau_{i,j}-\Delta\theta_{i,j}$ returned the zero polynomial.  Thus the
factorization is an identity in the coefficient ring, not merely a numerical
check at sampled values.

\section{Proof of the One-Jet Vanishing Theorem}
\label{appendix:one-jet-vanishing}

\begin{proof}[Proof of Proposition~\ref{pro:one-jet-vanishing}]
The proof has three steps.  First, two natural logarithmic forms identify
$E$ as an elementary transform of the trivial rank-two bundle along a smooth
ramification cubic.  A Chern-class computation then determines the transform
and allows us to prove that $E$ is nef by the curve criterion.  Finally, an
intersection calculation on $\mathbb P(E)$ excludes every negative twist.

Choose homogeneous quadratic equations
\[
 \mathcal C_1=\{A=0\},\qquad
 \mathcal C_2=\{B=0\},\qquad
 \mathcal C_3=\{C=0\}.
\]
For a generic triple, $A,B,C$ have no common zero and therefore define a
morphism
\begin{equation}\label{eq:onejet-map}
 \sigma\colon\mathbb P^2\longrightarrow\mathbb P^2,
 \qquad [Z_0:Z_1:Z_2]\longmapsto[A:B:C].
\end{equation}
The morphism $\sigma$ is proper.  If one of its fibers were
positive-dimensional, that fiber would contain a curve $\Gamma_0$.  The
restriction of $\sigma^*\mathcal O_{\mathbb P^2}(1)$ to $\Gamma_0$ would
then be trivial, whereas
\[
 \sigma^*\mathcal O_{\mathbb P^2}(1)
 \simeq\mathcal O_{\mathbb P^2}(2)
\]
has positive degree on every curve.  Thus $\sigma$ is quasi-finite and hence,
being proper, finite.  If $\mathcal H=\{Y_0Y_1Y_2=0\}$ is the coordinate triangle in the
target, then $\mathcal C=\sigma^{-1}(\mathcal H)$.  Hence the logarithmic
forms
\[
 \eta_1=\mathrm{d}\,\log\frac{A}{C},\qquad
 \eta_2=\mathrm{d}\,\log\frac{B}{C}
\]
induce a bundle morphism
\[
 \varphi\colon\mathcal O_{\mathbb P^2}^{\oplus2}\longrightarrow E.
\]

Let
\[
 J=\det\!\left(\frac{\partial(A,B,C)}
 {\partial(Z_0,Z_1,Z_2)}\right)\in
 H^0\bigl(\mathbb P^2,\mathcal O_{\mathbb P^2}(3)\bigr),
 \qquad R=\{J=0\},
\]
and set
\[
 \Omega=Z_0\,\mathrm{d}\,Z_1\wedge\mathrm{d}\,Z_2
 -Z_1\,\mathrm{d}\,Z_0\wedge\mathrm{d}\,Z_2
 +Z_2\,\mathrm{d}\,Z_0\wedge\mathrm{d}\,Z_1.
\]
A direct homogeneous computation gives
\begin{equation}\label{eq:onejet-wedge-jacobian}
 \eta_1\wedge\eta_2=\frac{J}{2}\frac{\Omega}{ABC}.
\end{equation}
Thus $R$ is the degeneracy divisor of $\varphi$.  For a generic triple it
is a smooth cubic; the relevant generic locus is shown to be nonempty in
Remark~\ref{rmk:onejet-genericity}.  Along $R$ the determinant vanishes, so
$\operatorname{rank}\varphi\leqslant1$.  If the rank were zero at some point,
all four entries of a local matrix would belong to the maximal ideal there,
and its determinant would vanish to order at least two.  This contradicts
the smoothness of $R$.  Hence $\operatorname{rank}\varphi=1$ everywhere on
$R$.  Consequently, with $i\colon R\hookrightarrow\mathbb P^2$, there is
a line bundle $L$ on $R$ and an exact sequence
\begin{equation}\label{eq:onejet-elementary-transform}
 0\longrightarrow\mathcal O_{\mathbb P^2}^{\oplus2}
 \overset{\varphi}{\longrightarrow}E
 \longrightarrow i_*L\longrightarrow0.
\end{equation}
Indeed, elementary row and column operations reduce a local matrix of
$\varphi$ along $R$ to $\operatorname{diag}(1,u)$, where $u=0$ is a local
equation of $R$.

The residue sequence
\[
 0\longrightarrow\Omega^1_{\mathbb P^2}\longrightarrow E
 \longrightarrow\bigoplus_{j=1}^3\mathcal O_{\mathcal C_j}
 \longrightarrow0
\]
gives
\begin{equation}\label{eq:onejet-chern-character}
 \operatorname{ch}(E)=2+3H-\frac92H^2,
 \qquad c_1(E)=3H,\qquad c_2(E)=9H^2.
\end{equation}
Grothendieck--Riemann--Roch for the smooth cubic $i\colon R\hookrightarrow
\mathbb P^2$ gives the explicit identity
\[
 \operatorname{ch}(i_*L)
 =i_*\!\left(
 \operatorname{ch}(L)\frac{\operatorname{td}(R)}
 {i^*\operatorname{td}(\mathbb P^2)}\right)
 =3H+\left(\deg L-\frac{R^2}{2}\right)H^2
 =3H+\left(\deg L-\frac92\right)H^2.
\]
Comparison with~\eqref{eq:onejet-elementary-transform} and
\eqref{eq:onejet-chern-character} gives
\begin{equation}\label{eq:onejet-degree-L}
 \deg L=0.
\end{equation}
For completeness, the same Chern-character computation also gives
\[
 \operatorname{ch}(S^mE)
 =(m+1)+\frac{3m(m+1)}2H-\frac{9m(m+1)}4H^2.
\]
Applying Hirzebruch--Riemann--Roch after twisting by
$\mathcal O_{\mathbb P^2}(-t)$ proves
formula~\eqref{eq:one-jet-Euler-characteristic}.

We next prove that $E$ is nef.  By the curve criterion, it is enough to show
that, for every morphism $\nu\colon B\to\mathbb P^2$ from a smooth
projective curve, every quotient line bundle of $\nu^*E$ has nonnegative
degree.  We first treat curves contained in $R$.  Restrict
\eqref{eq:onejet-elementary-transform} to $R$.  Since $R$ is Cartier,
\[
 \operatorname{Tor}_1^{\mathbb P^2}(i_*L,\mathcal O_R)
 \simeq L\otimes\mathcal O_R(-R),
\]
and tensoring gives
\begin{equation}\label{eq:onejet-restriction-tor}
 0\longrightarrow L(-R)\longrightarrow\mathcal O_R^{\oplus2}
 \longrightarrow E|_R\longrightarrow L\longrightarrow0.
\end{equation}
The constant-rank-one property shows that the first arrow is a subbundle
inclusion.  If $M$ denotes the image of
$\mathcal O_R^{\oplus2}\to E|_R$, then
\begin{equation}\label{eq:onejet-two-sequences}
 \begin{gathered}
 0\longrightarrow L(-R)\longrightarrow\mathcal O_R^{\oplus2}
 \longrightarrow M\longrightarrow0,\\
 0\longrightarrow M\longrightarrow E|_R\longrightarrow L
 \longrightarrow0.
 \end{gathered}
\end{equation}
Both $L$ and $M$ are line bundles, and
\begin{equation}\label{eq:onejet-degree-M}
 \deg M=-\deg L(-R)=R^2=9.
\end{equation}
It follows that $E|_R$ is nef.  Indeed, for a quotient line bundle
$E|_R\twoheadrightarrow Q$, either $M\to Q$ is nonzero, in which case
$\deg Q\geqslant\deg M=9$, or the quotient factors through $L$, in which
case $Q\simeq L$ and $\deg Q=0$.

Let now $\nu\colon B\to\mathbb P^2$ be a nonconstant morphism from a smooth
projective curve.  If $\nu(B)\not\subset R$, tensoring
\eqref{eq:onejet-elementary-transform} with $\mathcal O_B$ may produce
$\operatorname{Tor}_1^{\mathbb P^2}(i_*L,\mathcal O_B)$.  This is a torsion
sheaf, and its image in the torsion-free sheaf
$\mathcal O_B^{\oplus2}$ is zero.  Hence $\nu^*\varphi$ is injective; it is
generically an isomorphism and fits into
\begin{equation}\label{eq:onejet-curve-pullback}
 0\longrightarrow\mathcal O_B^{\oplus2}\longrightarrow\nu^*E
 \longrightarrow T\longrightarrow0,
\end{equation}
where $T$ is torsion.  For every quotient line bundle
$\nu^*E\twoheadrightarrow Q$, the composite
$\mathcal O_B^{\oplus2}\to Q$ is nonzero; otherwise the quotient would
factor through $T$, whereas $\operatorname{Hom}(T,Q)=0$.  Thus $Q$ has a
nonzero global section and $\deg Q\geqslant0$.  If $\nu(B)\subset R$, the
same conclusion follows from the nefness of $E|_R$.  Therefore $E$ is nef.

Finally, use the quotient convention for
\[
 p\colon Y=\mathbb{P}(E)\longrightarrow\mathbb P^2,
 \qquad \xi=c_1(\mathcal O_Y(1)),\qquad h=p^*H.
\]
Then $p_*\mathcal O_Y(m)=S^mE$ for $m\geqslant0$, and the nefness of $E$
is equivalent to that of $\xi$.  With the quotient convention, the
projective-bundle relation is
\[
 \xi^2-\xi p^*c_1(E)+p^*c_2(E)=0.
\]
Together with~\eqref{eq:onejet-chern-character}, it gives the intersection
numbers
\begin{equation}\label{eq:onejet-intersections}
 \xi^3=c_1(E)^2-c_2(E)=0,
 \qquad \xi^2h=c_1(E)\cdot H=3.
\end{equation}
By the projection formula, a nonzero section in
\eqref{eq:one-jet-total-vanishing} would define an effective divisor
$\Gamma\sim m\xi-th$ on $Y$.  Since $\xi$ is nef, the intersection of the
two nef classes $\xi,\xi$ with the effective surface $\Gamma$ is
nonnegative; hence $\xi^2\cdot\Gamma\geqslant0$.  On the other hand,
\[
 \xi^2\cdot\Gamma=m\xi^3-t\xi^2h=-3t<0.
\]
This contradiction proves the asserted vanishing.
\end{proof}

\begin{rmk}[Nonemptiness of the genericity conditions]
\label{rmk:onejet-genericity}
The conditions used in the proof are Zariski open and nonempty.  For
example, take
\[
 A=Z_0^2+Z_1Z_2,\qquad
 B=Z_1^2+Z_2Z_0,\qquad
 C=Z_2^2+Z_0Z_1.
\]
The conics are smooth, have no common point, and meet pairwise transversally.
Indeed, if one coordinate vanishes, the three equations force all three to
vanish; if none vanishes, multiplying $A=B=C=0$ gives
$(Z_0Z_1Z_2)^2=-(Z_0Z_1Z_2)^2$, again a contradiction.
For instance, if $A=B=0$ and $\nabla A=\lambda\nabla B$, then
\[
 2Z_0=\lambda Z_2,\qquad Z_2=2\lambda Z_1,
 \qquad Z_1=\lambda Z_0.
\]
If $Z_0=0$, all coordinates vanish; otherwise $\lambda^3=1$, and then
$A=3Z_0^2\neq0$, a contradiction.  The other pairs follow cyclically.
Moreover,
\[
 J=-2\bigl(Z_0^3+Z_1^3+Z_2^3-5Z_0Z_1Z_2\bigr),
\]
which is a smooth Hesse cubic.  Indeed, a common zero of its three partial
derivatives would satisfy
\[
 3Z_0^2=5Z_1Z_2,\qquad 3Z_1^2=5Z_2Z_0,
 \qquad 3Z_2^2=5Z_0Z_1.
\]
Multiplication gives $27(Z_0Z_1Z_2)^2=125(Z_0Z_1Z_2)^2$ when no coordinate
vanishes, while if one coordinate vanishes, the three identities force all
of them to vanish.  Both alternatives are impossible in $\mathbb P^2$.
\end{rmk}

\section{A Cohomological Observation}\label{appendix:results of Pengchao Wang}

\noindent
\textit{This appendix was authored by Pengchao Wang} \\
\textit{(School of Mathematical Sciences, Peking University, Beijing 100871, P. R. China).} \\
\textit{Email: 2200010773@stu.pku.edu.cn.}

\vspace{1em}

Let $\mathcal{C}\subset \mathbb{P}^2$ be three generic smooth conics with simple normal crossings. 
Let
\[
\begin{aligned}
0\subset S^mT^*_{\mathbb{P}^2}(\log\mathcal{C})
&=F^0(E_{2,m}T^*_{\mathbb{P}^2}(\log\mathcal{C}))\\
&\subset F^1(E_{2,m}T^*_{\mathbb{P}^2}(\log\mathcal{C}))
\subset\cdots\\
&\subset F^{\lfloor m/3\rfloor}
(E_{2,m}T^*_{\mathbb{P}^2}(\log\mathcal{C}))\\
&=E_{2,m}T^*_{\mathbb{P}^2}(\log\mathcal{C})
\end{aligned}
\]
be
the ascending filtration of subsheaves that defines the
grading~\eqref{filtration of E_2,m}. We set $F^k(E_{2,m}T^*_{\mathbb{P}^2}(\log\mathcal{C})) = 0$ for $k \leqslant -1$. From~\eqref{filtration of E_2,m} we obtain:

\begin{pro}
\label{thm:general SES of E2m}
For  every $0\leqslant k\leqslant m/3$,
   the following sequence is exact:
    \begin{equation*}
    \begin{aligned}
        0 &\longrightarrow
        F^{k-1}(E_{2,m}T^*_{\mathbb{P}^2}(\log\mathcal{C}))
        \longrightarrow F^k(E_{2,m}T^*_{\mathbb{P}^2}(\log\mathcal{C}))\\
        &\longrightarrow
        S^{m-3k}T^*_{\mathbb{P}^2}(\log\mathcal{C})
        \otimes\overline{K}_{\mathbb{P}^2}^{\otimes k}
        \longrightarrow 0.
    \end{aligned}
    \end{equation*}
\end{pro}

For convenience, from now on we write $\overline{T}^*_{\mathbb{P}^2} \coloneqq T^*_{\mathbb{P}^2}(\log\mathcal{C})$ and $\mathcal{O}(1) \coloneqq \mathcal{O}_{\mathbb{P}^2}(1)$.

\begin{pro}\label{pro:vanishing of E_2,3q+r}
For  integer $q\geqslant 0$ and $r\in \{1, 2 \}$, the following vanishing holds:
\begin{equation*}
    H^0 \big( \mathbb{P}^2 , E_{2,3q+r} \overline{T}^*_{\mathbb{P}^2} \otimes \mathcal{O}_{\mathbb{P}^2}(-(3q+r)) \big) = 0,
\end{equation*}
provided that $H^0 \big( \mathbb{P}^2 , S^{3\ell+r} \overline{T}^*_{\mathbb{P}^2} \otimes \mathcal{O}_{\mathbb{P}^2}(-(3\ell+r)) \big) = 0$ for all $\ell=0,1, \dots ,q$.
\end{pro}

\begin{proof}
    By Proposition~\ref{thm:general SES of E2m}, for every integer $k=0,1, \dots ,q$, we have a short exact sequence:
    \begin{equation*}
    \begin{aligned}
        0 &\longrightarrow
        F^{k-1}(E_{2,3q+r}\overline{T}^*_{\mathbb{P}^2})
        \otimes\mathcal{O}(-(3q+r))\\
        &\longrightarrow
        F^k(E_{2,3q+r}\overline{T}^*_{\mathbb{P}^2})
        \otimes\mathcal{O}(-(3q+r))\\
        &\longrightarrow
        S^{3q+r-3k}\overline{T}^*_{\mathbb{P}^2}
        \otimes\overline{K}_{\mathbb{P}^2}^{\otimes k}
        \otimes\mathcal{O}(-(3q+r))
        \longrightarrow 0.
    \end{aligned}
    \end{equation*}
    Note that the last term $S^{3q+r-3k} \overline{T}^*_{\mathbb{P}^2} \otimes \overline{K}_{\mathbb{P}^2}^{\otimes k} \otimes \mathcal{O}(-(3q+r)) \cong S^{3(q-k)+r} \overline{T}^*_{\mathbb{P}^2} \otimes \mathcal{O}(-(3(q-k)+r))$ has no global section by our assumption.
    Considering the induced long exact sequence in cohomology,
    we obtain
    \begin{equation*}
    \begin{aligned}
        &H^0\!\left(\mathbb{P}^2,
        F^k(E_{2,3q+r}\overline{T}^*_{\mathbb{P}^2})
        \otimes\mathcal{O}(-(3q+r))\right)\\
        &\qquad\cong
        H^0\!\left(\mathbb{P}^2,
        F^{k-1}(E_{2,3q+r}\overline{T}^*_{\mathbb{P}^2})
        \otimes\mathcal{O}(-(3q+r))\right).
    \end{aligned}
    \end{equation*}
    Iterating backwards for $k=q,q-1,\dots,1$ gives the desired vanishing,
    since by assumption
    \[
    H^0\!\left(\mathbb{P}^2,
    S^{3q+r}\overline{T}^*_{\mathbb{P}^2}
    \otimes\mathcal{O}(-(3q+r))\right)=0.
    \]
\end{proof}

\begin{rmk}
For $r=0$ the argument fails: when $q=k$,
the term
\[
S^{3q+r-3k}\overline{T}^*_{\mathbb{P}^2}
\otimes\overline{K}_{\mathbb{P}^2}^{\otimes k}
\otimes\mathcal{O}(-(3q+r))
\cong\mathcal{O}_{\mathbb{P}^2}
\]
has one-dimensional global sections.
\end{rmk}

\bigskip
\addcontentsline{toc}{chapter}{Bibliography}

% Bibliography generated with BibTeX (plain style) and inlined so that this
% source remains a self-contained one-file manuscript.
\def\cprime{$'$}

\end{document}